\documentclass[preprint,review,3p,12pt]{elsarticle}

\usepackage{amssymb,dsfont}
\usepackage{lineno}
\usepackage{graphicx}
\usepackage{epstopdf}

\newtheorem{definition}{Definition}
\newtheorem{principle}{Principle}
\biboptions{comma,sort}

\journal{Journal of Computational Physics}

\begin{document}

\begin{frontmatter}
  \title{PDEs on moving surfaces via the closest point method  and a modified grid based particle method}
  \author[rvt]{A.~Petras\corref{cor1}}
  \ead{apetras@sfu.ca}
  \author[rvt]{S.J.~Ruuth}
  \ead{sruuth@sfu.ca}
  \cortext[cor1]{Corresponding author}
  \address[rvt]{Department of Mathematics, Simon Fraser University, Burnaby, British Columbia, Canada V5A1S6}

  \begin{abstract}
Partial differential equations (PDEs) on surfaces arise in a wide range of applications.    The closest point method
(Ruuth and Merriman, J. Comput. Phys. 227(3):1943–-1961, [2008]) is a recent embedding method that has been
used to solve a variety of PDEs on smooth surfaces using a closest point representation of the surface and
standard Cartesian grid methods in the embedding space.
The original closest point method (CPM) was designed for problems posed on static surfaces, however the
solution of PDEs on moving surfaces is of considerable interest as well.
Here we propose solving PDEs on moving surfaces using a combination of the CPM and
a modification of the grid based particle method (Leung and Zhao, J. Comput. Phys. 228(8):2993--3024, [2009]).
The grid based particle method (GBPM)
represents and tracks surfaces using meshless particles and an Eulerian reference grid.
Our modification of the GBPM introduces a reconstruction step into the original method
to ensure that all the grid points within a computational tube surrounding the surface are active.
We present a number of examples to illustrate the numerical convergence properties of our combined method.
Experiments for advection-diffusion equations  that are strongly coupled to the velocity of the surface are also presented.
  \end{abstract}

\end{frontmatter}

\section{Introduction}
Partial differential equations (PDEs) on moving surfaces arise throughout the natural and applied sciences. Biological applications include the modelling of biomembranes \cite{elliott2010modeling}, cell motility and chemotaxis \cite{elliott2012modelling}, and pattern formation \cite{barreira2011surface,venkataraman2011modeling}. In material science, PDEs on moving surfaces arise in dealloying \cite{eilks2008numerical}, while in fluid dynamics they appear in the modelling of the mass conservation/transportation of surfactants in two phase flows \cite{adalsteinsson2003transport,james2004surfactant}. A wide range of applications also appear in computer graphics. See \cite{auer2013semi} for some examples.

A variety of methods have been used to solve PDEs on moving surfaces.
For example, Dziuk and Elliott \cite{dziuk2007finite} used
moving triangle meshes as part of a finite element method for conservation laws.
Their method is an example of a surface finite element method.  (See \cite{dziuk2013finite} for further details and
references related to the numerical analysis and applications of surface finite elements.)
While finite element methods are commonly used with moving triangle meshes, other discretizations have also been proposed.
For example, Nemadjieu \cite{nemadjieu2012finite} used a finite volume discretization to solve advection-diffusion equations on moving triangle meshes.

Particle and parametrization methods have also been used to represent and move surfaces.
Of particular interest to us is the work of Leung and Zhao \cite{leung2009grid}, where these concepts are combined into the grid based particle method.
The grid based particle method was used in combination with a local parametrization method to solve  advection-diffusion equations on moving surfaces in \cite{leung2011grid}.
It was also used with the Gaussian beam summation method to approximate the high frequency asymptotic solution of the Helmholtz equation on moving surfaces in \cite{leung2010gaussian}.
We provide a detailed description of the grid based particle method in Section~\ref{GBPMSection}.

Level set methods \cite{sethian1999level,osher2006level} are another major class of methods for representing and moving surfaces.
In the level set approach, the PDE-on-surface problem is solved on a narrow computational domain surrounding the surface.
Methods of this type include the  semi-implicit finite difference scheme of Xu and Zhao \cite{xu2003eulerian} and the
 finite element discretization of Dziuk and Elliott \cite{dziuk2010eulerian}.   Related to the level set representation of a surface, is the closest point representation
of a surface, whereby grid nodes store the closest point in Euclidean distance to the surface.
Closest point representations constructed from moving triangulated surfaces
have been used in a closest point method for the Navier-Stokes and wave equations in \cite{auer2013semi}.
Closest point representations constructed from level sets have also been used; see  \cite{kim2013closest} for a closest point method
 for irrotational flow on moving surfaces.
We provide a detailed description of the closest point method \cite{ruuth2008simple} in Section~\ref{CPMSection}.

In this paper, we combine the closest point method  (CPM) and the grid based particle method (GBPM) to
solve PDEs on moving surfaces.
The closest point method is a simple embedding method for solving rather general PDEs on smooth stationary surfaces.
It decouples surface geometry and PDE evolution via a closest point extension step, yielding a method that involves standard Cartesian grid methods in the embedding space.
The original method only considered stationary surfaces, however.   Complementing the CPM is the GBPM.
In the GBPM, a closest point representation is formed as part of an algorithm to move surfaces according to curvature-dependent (and higher-order \cite{leung2011grid}) motions.
In this paper, we combine the CPM and a (small) modification of the GBPM to solve PDEs on moving surfaces.
Our proposed method is a highly modular combination of two well-known, tested methods.
Indeed,  in our approach, the component methods are applied in an alternating fashion, yielding a particularly simple
method to solve PDEs on moving surfaces.


The paper is organized as follows. In Section~\ref{ProblemSection}, we present the PDE model on the evolving surface and we review the CPM and the GBPM. Section~\ref{NewMethodDescriptionSection} introduces our  modification of the GBPM. Numerical experiments for a number of geometric motion laws are provided to verify the correctness of the algorithm. In Section~\ref{CoupledMethodSection}, we couple the modified GBPM with the CPM to obtain a method for PDEs on evolving surfaces. A variety of numerical experiments are provided to illustrate the method. Finally, in Section~\ref{SummarySection} we draw conclusions and we discuss possible future work.

\section{Mathematical formulation and numerical methods review}\label{ProblemSection}

In this section, we introduce notation and present the PDE-on-surface model.  Also, we briefly review the closest point method (CPM) and the grid based particle method (GBPM).

\subsection{Notation and formulation of PDE}\label{PDEsSection}

Advection-diffusion is a fundamental PDE model. On an evolving surface $\Gamma(t)$, the advection-diffusion PDE takes the form
\begin{equation}\label{AdvDifPDE}
  u_t+\mathbf{v}\cdot\nabla u+u\nabla_\Gamma\cdot\mathbf{v}-\nabla_\Gamma\cdot \mathbf{q} = f\qquad \textrm{on } \Gamma(t)
\end{equation}
where $\mathbf{v}$ is the velocity of the surface, $\mathbf{q}$ is the flux for $u$, $f$ is a source term and $u$ is a scalar quantity. By splitting the velocity into normal and tangential components, $\mathbf{v}=V\mathbf{n}+\mathbf{T}$, an alternative formulation of (\ref{AdvDifPDE}) can be derived \cite{dziuk2013finite}. This has the form
\begin{equation}\label{AdvDifPDE2}
  u_t+V\frac{\partial u}{\partial n}-V\kappa u+\nabla_\Gamma\cdot(u\mathbf{T}) -\nabla_\Gamma\cdot \mathbf{q} = f\qquad \textrm{on } \Gamma(t)
\end{equation}
where $\kappa$ is the mean curvature of the surface.

Throughout this paper, we will consider (\ref{AdvDifPDE2}) with the diffusion flux,
\begin{equation}\label{DiffusionFlux}
\mathbf{q}=\mathcal{D}\nabla_\Gamma u
\end{equation}
where $\mathcal{D}>0$ is the diffusivity parameter. Other important cases include the Cahn-Hilliard flux
\begin{equation}\label{CahnHilliardFlux}
\mathbf{q} = b(u)\nabla_{\Gamma} w
\end{equation}
where $b(u)$ is the mobility, $w = -c\Delta_\Gamma u+g(u)$ is the chemical potential and $g$ is the first derivative of a double-well free energy. More information on these flux choices can be found in \cite{dziuk2013finite}.

To complete the initial value problem, initial conditions $u(\mathbf{x},0)$ are imposed. The surfaces considered throughout this work are closed, therefore boundary conditions are not required.

\subsection{The closest point method}\label{CPMSection}
The closest point method \cite{ruuth2008simple} is a numerical method for approximating PDEs on stationary surfaces. In the original formulation \cite{ruuth2008simple}, a uniform finite difference grid is formed around the surface. Grid nodes store the coordinates of the closest point in Euclidean distance to the surface:
\begin{definition}
  Let $\mathbf{x}$ be some point in the embedding space $\mathds{R}^d$. Then,
  $$cp_\Gamma(\mathbf{x})=\arg\min_{\mathbf{z}\in \Gamma}\|\mathbf{x}-\mathbf{z}\|_2$$
  is the closest point of $\mathbf{x}$ to the surface $\Gamma$.
\end{definition}
This gives a closest point representation of the surface in the embedding space.

The evolution of the surface PDE is carried out using an embedding PDE. To form the embedding PDE, derivatives intrinsic to the surface in the original PDE are replaced with a combination of closest point operators and derivatives expressed in standard Cartesian coordinates in the embedding space. Two principles are used to define the embedding PDE. The first of these relies on the fact that if surface values are extended off the surface in the normal direction, then, at the surface, the standard gradient in the embedding space agrees with the gradient intrinsic to the surface. This fundamental idea is expressed in Principle~\ref{equivalence of gradients}, \emph{equivalence of gradients}:
\begin{principle}\label{equivalence of gradients}
  Let $v$ be any function on $\mathds{R}^d$ that is constant along normal directions of $\Gamma$. Then, at the surface, intrinsic gradients are equivalent to standard gradients, $\nabla_\Gamma v=\nabla v$.
\end{principle}
Similarly, the extension of surface vectors into the embedding space leads to Principle~\ref{equivalence of divergence}, \emph{equivalence of divergence}:
\begin{principle}\label{equivalence of divergence}
  Let $\mathbf{v}$ be any vector field on $\mathds{R}^d$ that is tangent to $\Gamma$ and tangent to all surfaces displaced by a fixed distance from $\Gamma$. Then, at the surface, $\nabla_\Gamma\cdot\mathbf{v}=\nabla\cdot\mathbf{v}$.
\end{principle}
Combinations of these two principles yield embeddings for more general differential operators \cite{ruuth2008simple,macdonald2009implicit}. Of particular interest for our advection-diffusion problem (\ref{AdvDifPDE2}) is the \emph{equivalence of diffusion}. Specifically, by combining Principles~\ref{equivalence of gradients} and \ref{equivalence of divergence}, we may replace the Laplace-Beltrami operator with the standard Laplacian after a constant normal extension of surface values $u$.

To initialize the CPM, we set up the closest point representation of the surface on a tubular computational domain surrounding the surface, and extend the initial values of $u$ onto the computational domain $\Omega_c$. Time-stepping for the explicit CPM is carried out by alternating two steps:
\begin{enumerate}
  \item \textbf{Closest point extension}.  The solution on $\Gamma$ is extended to the computational domain by replacing $u$ with $u(cp_\Gamma)$ for all $\mathbf{x}\in\Omega_c$.
  \item \textbf{Evolution}.  The embedding PDE is solved on the tubular computational domain $\Omega_c$ for one time step (or one stage of a Runge-Kutta method).
\end{enumerate}

Note that the closest point extension defined in the first step will require interpolation since $cp_\Gamma(\mathbf{x})$ is not necessarily a grid point in $\Omega_c$. Following \cite{ruuth2008simple}, interpolation is carried out using barycentric Lagrange interpolation with polynomial degree $p=q+r-1$, where $q$ is the order of the differencing scheme and $r$ is the derivative order. The radius of the computational tube depends on the width of the differencing and interpolation stencils. In our computations for the advection-diffusion equation (\ref{AdvDifPDE2}), second-order centered differences are used. This leads to a computational tube radius of
\begin{equation}\label{Bandwidth}
  \gamma_{CPM}=\sqrt{(d-1)\left(\frac{p+1}{2}\right)^2+\left(1+\frac{p+1}{2}\right)^2}\Delta x
\end{equation}
in the $d$-dimensional embedding space, using polynomial interpolation of degree $p$ (our examples take $p=3$) \cite{ruuth2008simple}.

Throughout this paper, we consider the explicit CPM \cite{ruuth2008simple}, however, implicit time-stepping could alternatively be used. See \cite{macdonald2009implicit,Macdonald20117944} for details.

\subsection{The grid based particle method}\label{GBPMSection}
We now review Leung and Zhao's grid based particle method (GBPM) for capturing the motion of the surface \cite{leung2009grid}.
Notably, it constructs a closest point representation of the surface for each time step,
which makes the study of its combination with the CPM particularly interesting.
Our interest will be on curvature-dependent flows on closed surfaces,
however the method has also been applied to open surfaces \cite{leung2009gridopen}
and higher-order geometric motions \cite{leung2011grid}.

To initialize the GBPM, a grid is constructed that contains the surface. Grid points that are within a Euclidean distance $\gamma_{GBPM}$ (the tube radius) of the surface are identified: These are the \emph{active grid points}. As part of the construction of the computational tube, we compute the closest point on the surface for each active grid point. The nodal values are called \emph{footpoints} in the GBPM and give the surface representation.

After initialization, the system is evolved in time.  Each time step of size $\Delta t$ consists of three steps:
\begin{enumerate}
   \item \textbf{Motion}: The footpoints are moved according to the desired motion law.
  \item \textbf{Resampling}: For each active grid point, the closest point to the surface (as defined by the footpoints) is computed. This gives the updated footpoints.
  \item \textbf{Update of the Computational Tube}: This step consists of two stages. The first stage activates all the grid points that have neighboring active grid points and applies the resampling step to find their corresponding footpoints. The second stage deactivates all the grid points that are far from the surface (i.e., the distance between the grid point and its footpoint is larger than the tube radius $\gamma_{GBPM}$).
\end{enumerate}

\begin{figure}
    \centering
    \includegraphics[width=0.49\textwidth]{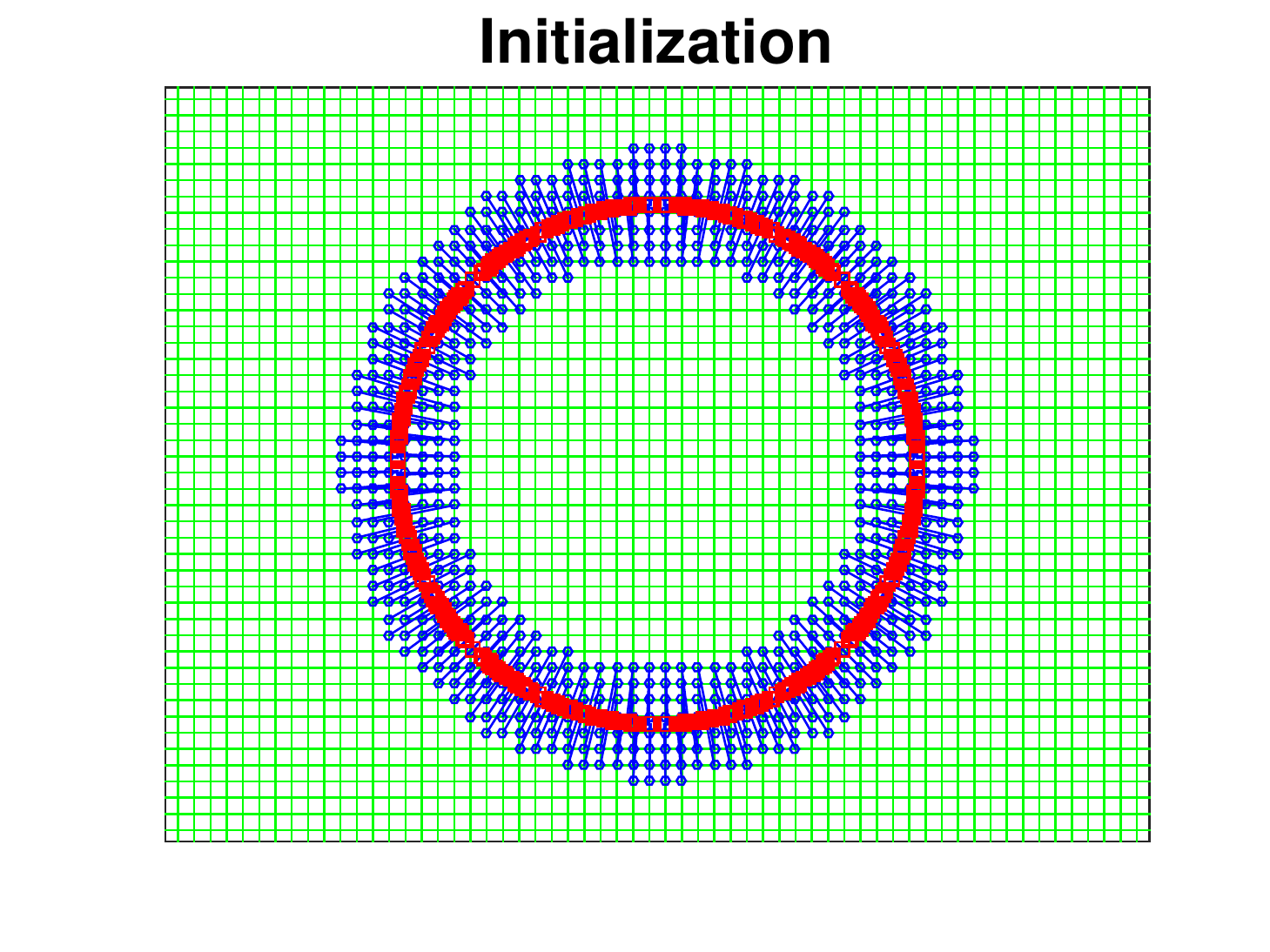}\\

    \includegraphics[width=0.49\textwidth]{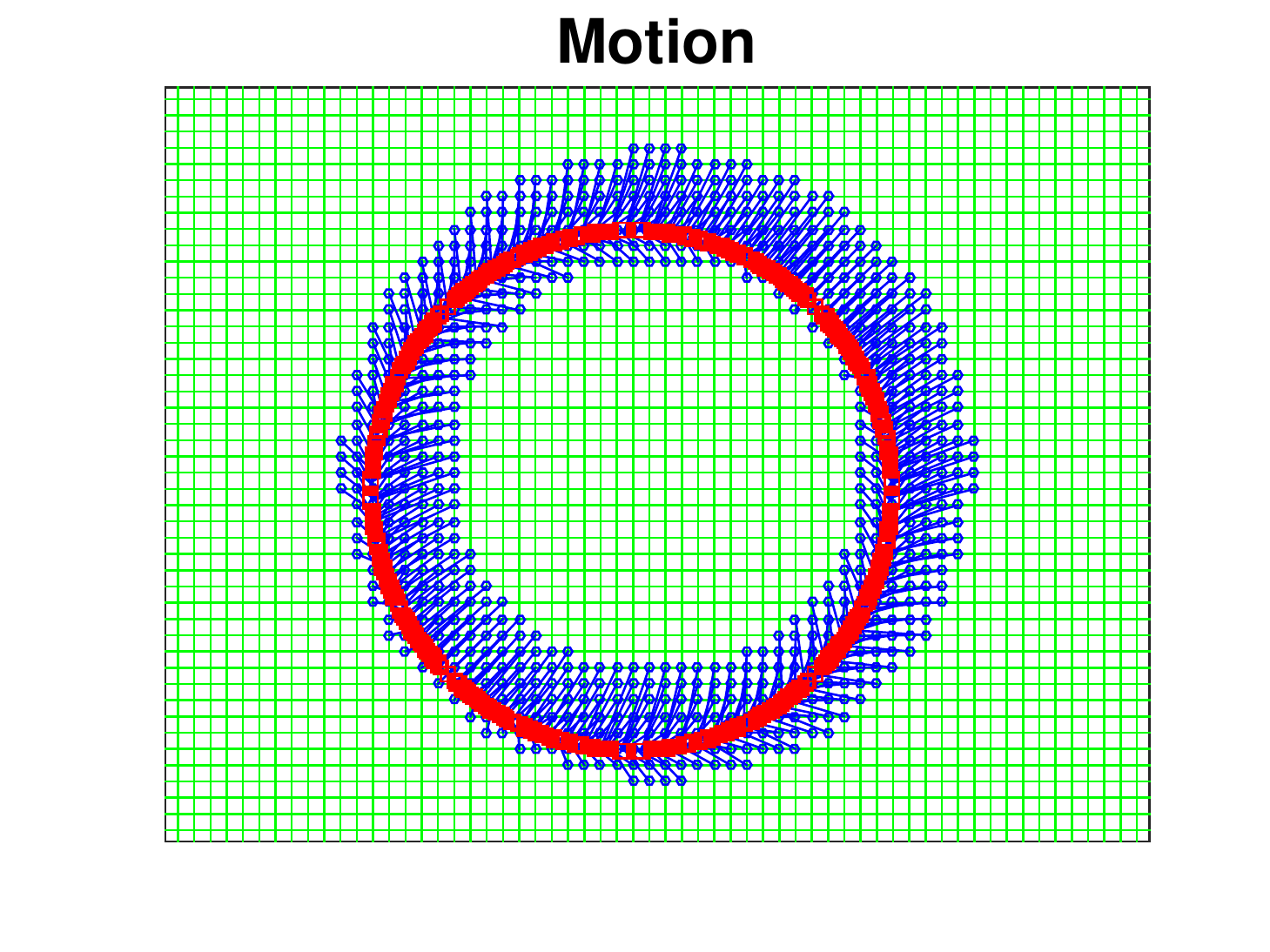}
    \includegraphics[width=0.49\textwidth]{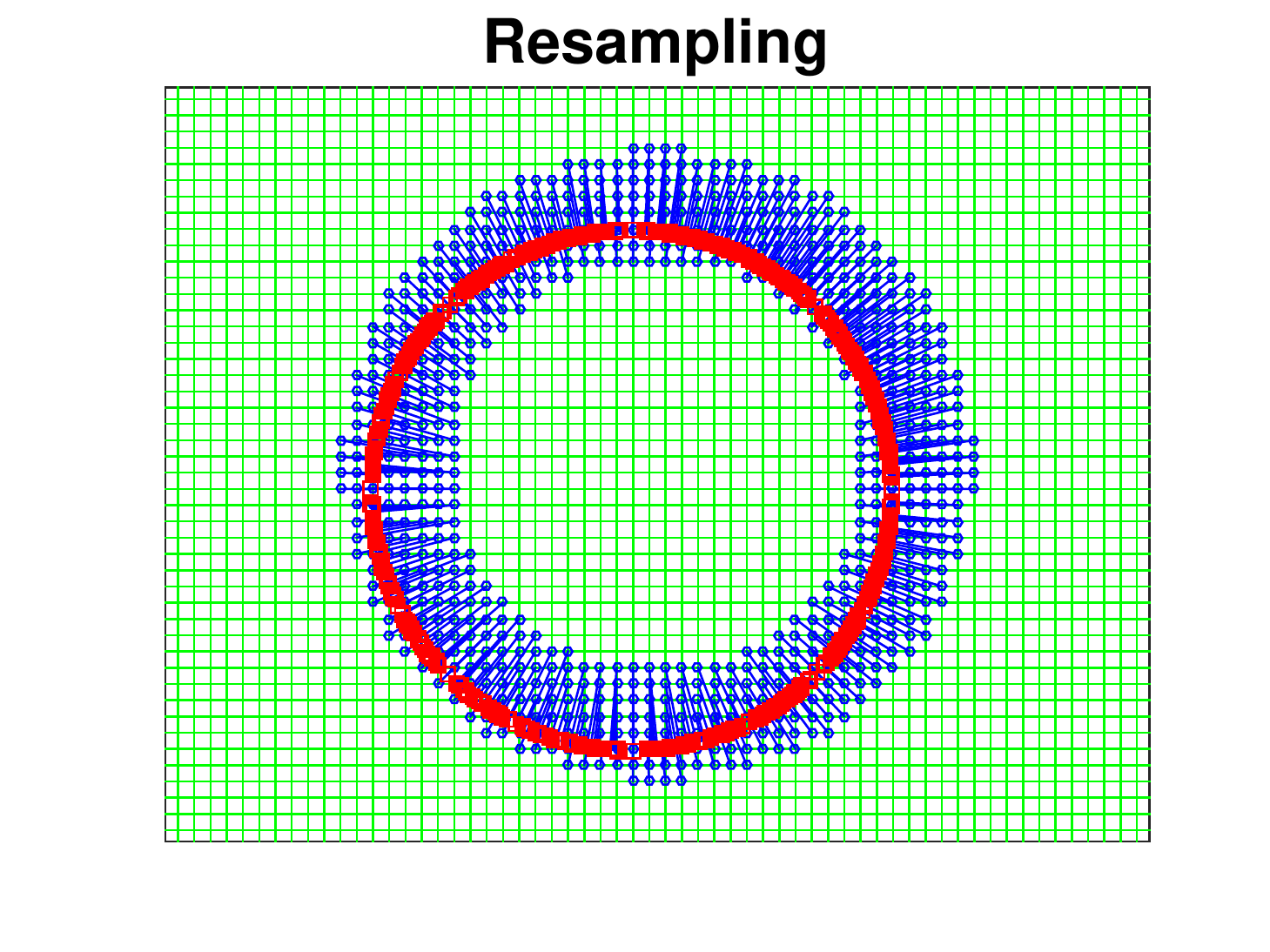}\\

    \includegraphics[width=0.49\textwidth]{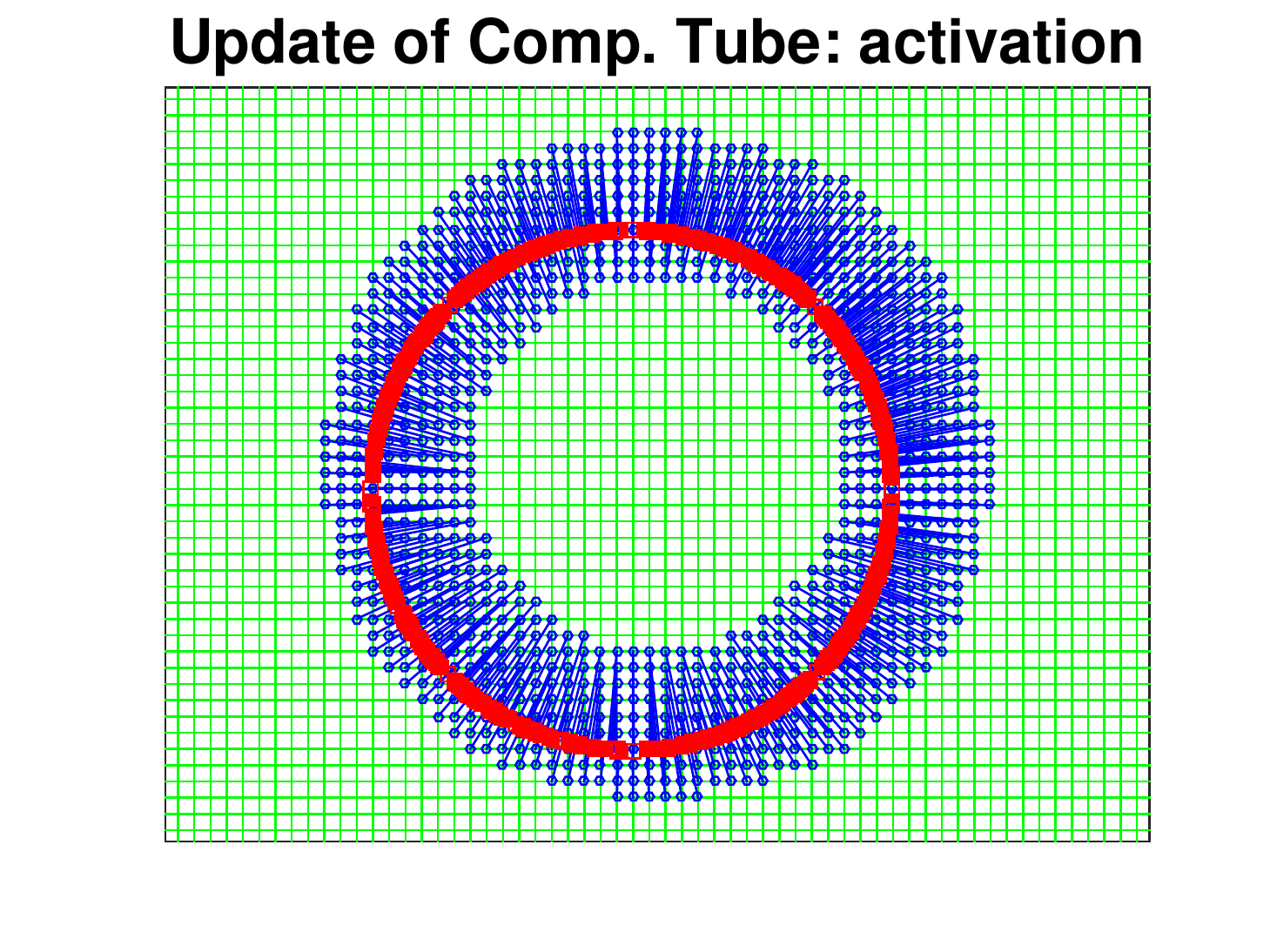}
    \includegraphics[width=0.49\textwidth]{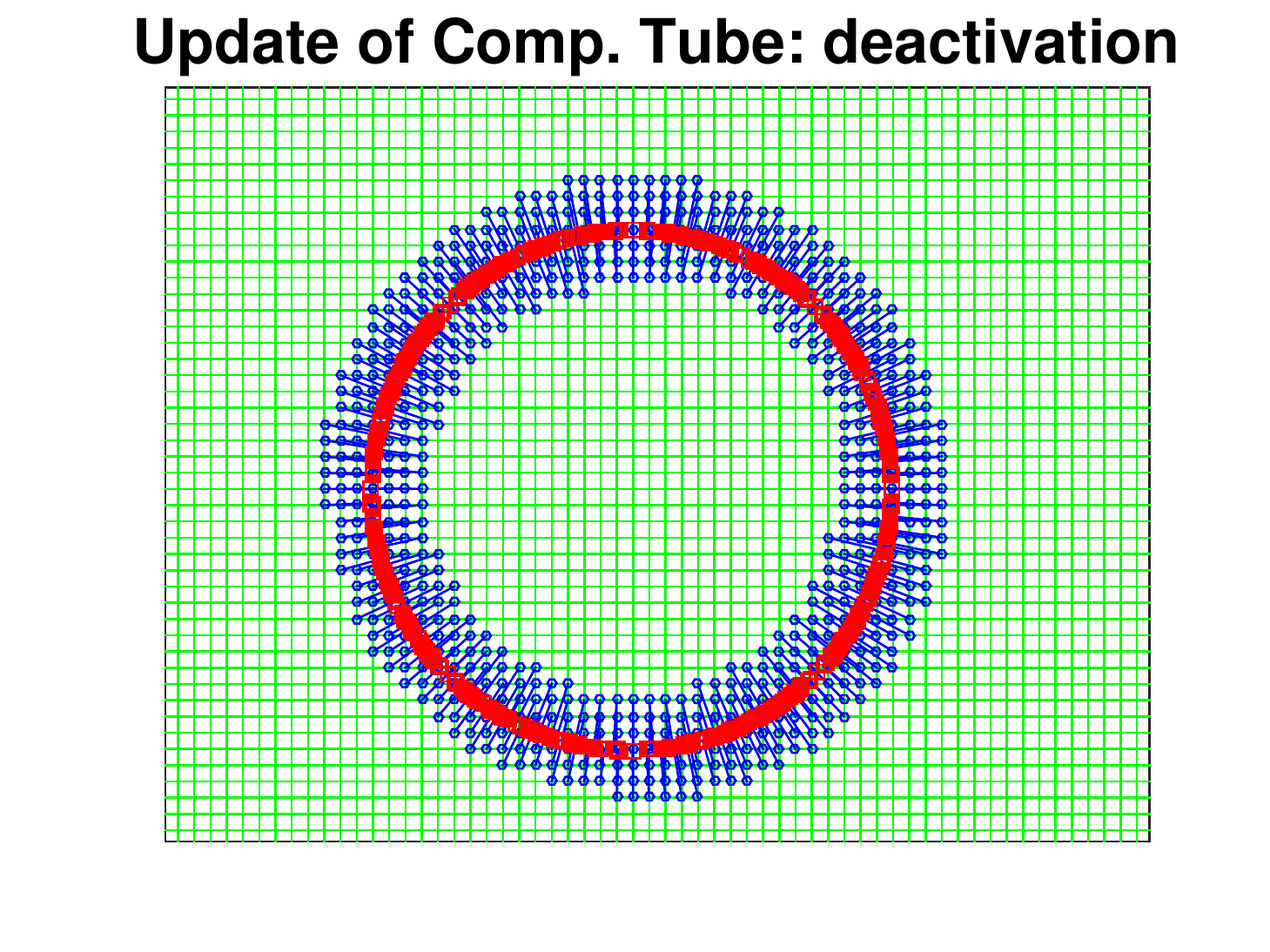}
    \caption{An illustration of the main steps of the GBPM (from top left to bottom right). Active grid points (blue dots) are connected to their footpoints (red dots) with blue lines. The green lines correspond to the grid.}
    \label{GBPM}
\end{figure}

An illustration of the main steps of the GBPM is provided in Figure~\ref{GBPM}. The initialization sets all the grid points that are within a distance $\gamma_{GBPM}$ of the surface to active, and sets the footpoints to be the corresponding closest points. This yields a closest point representation over a computational tube of radius $\gamma_{GBPM}$. The motion step evolves the surface by moving the footpoints; the footpoints no longer give a closest point representation of the surface. In the resampling step, the footpoint of each active grid point is re-assigned to be the closest point on the surface. Finally, we update the computational tube by adding points that are neighboring the active grid points and by deleting points outside the tube radius. This restores the tubular closest point representation.

\subsubsection{The Resampling Step} \label{GBPMSubsection}
We now provide a more detailed review of the resampling step of the GBPM since this step will be a focus for our modified algorithm described in Section \ref{RemarksSection}.
(For more information and implementation details on the GBPM, see \cite{leung2009grid,leung2009gridopen,leung2011grid}.)

Resampling begins by  collecting the footpoints that will be used.   To do this,
we identify neighboring points and their corresponding footpoints for each active grid point $\mathbf{p}$.
From these, the $m$ closest footpoints to $\mathbf{p}$ are determined.
(These points need to be distinct and in practice we impose that they are at least some minimum distance $\delta$ apart.)

Next, a local coordinate system is defined by the normal and tangential vector(s) of the closest footpoint to $\mathbf{p}$, and an interpolating function is constructed by the least squares method. Using Newton's method, we obtain the point on the local reconstruction function that minimizes the Euclidean distance from the point $\mathbf{p}$.
Some conditions are imposed on the minimizing points obtained in the resampling step.
Specifically, the curvature cannot be too large (we require $\kappa<1/\Delta x$),
and the minimizer must lie between the footpoints defining the reconstruction.
If these conditions fail, then the grid point and its footpoint are deactivated.
Note that any needed geometric information, such as the local normal, or mean curvature $\kappa$,
is computed and stored as part of this step.

We further note that in cases where two segments of the surface are close to one another,
a condition on the maximum angle  between unit normal vectors is imposed to ensure that the chosen footpoints lie on the same segment.
Finally, in cases of topological change, the footpoints that contain Lagrangian inconsistencies in their local neighborhood are removed.

\section{A modified GBPM for geometric surface motion}\label{NewMethodDescriptionSection}
In this section,
a modified GBPM is proposed which is amenable to coupling with the CPM. Numerical experiments on a variety of geometric motions are provided to verify the correctness of the new algorithm.

\subsection{A modified GBPM}\label{RemarksSection}
The GBPM constructs a closest point representation of the surface at every time step of a surface evolution, making it a natural candidate to consider for coupling with the CPM.   However, a challenge must be overcome before combining the methods:
In the CPM,  all nodes within the computational tube must have valid closest point values,
otherwise a standard implementation of the method will fail.   This condition is not met with the GBPM, since the method deactivates grid points
within the computational tube in certain situations (see Section~\ref{GBPMSection}).

Note that grid point deactivation is expected to occur with greater frequency in a coupled method than in the standard GBPM.
The reason for this is that a wider computational tube is required in the CPM than in the standard GBPM,
leading to the use of a wider tube in a coupled method as well.
To illustrate the difference, consider a CPM discretization of diffusion using a second-order finite difference scheme and cubic polynomial interpolation.
In the CPM, this leads to an approximate tube radius of $3.6\Delta x$ in two dimensions and $4.1\Delta x$ in three dimensions.  In contrast,
a standard GBPM implementation in 2 or 3D would use a relatively narrow tube radius of up to $1.5\Delta x$.

To obtain a method for geometric surface motion that is compatible with the CPM, we introduce a (slightly) modified GBPM.   The {\it modified GBPM}
applies the original GBPM for the {\bf initialization}, {\bf motion} (step 1),  and {\bf update of the computational tube} (step 3).   For  {\bf resampling} (step 2),
a change is introduced: Whenever a footpoint deactivation occurs in the original GBPM resampling step, the modified method constructs
an osculating circle/sphere to locally approximate the curve/surface.  The closest point on the reconstruction is accepted as the updated footpoint.

\begin{figure}
    \centering
    \includegraphics[width=0.49\textwidth]{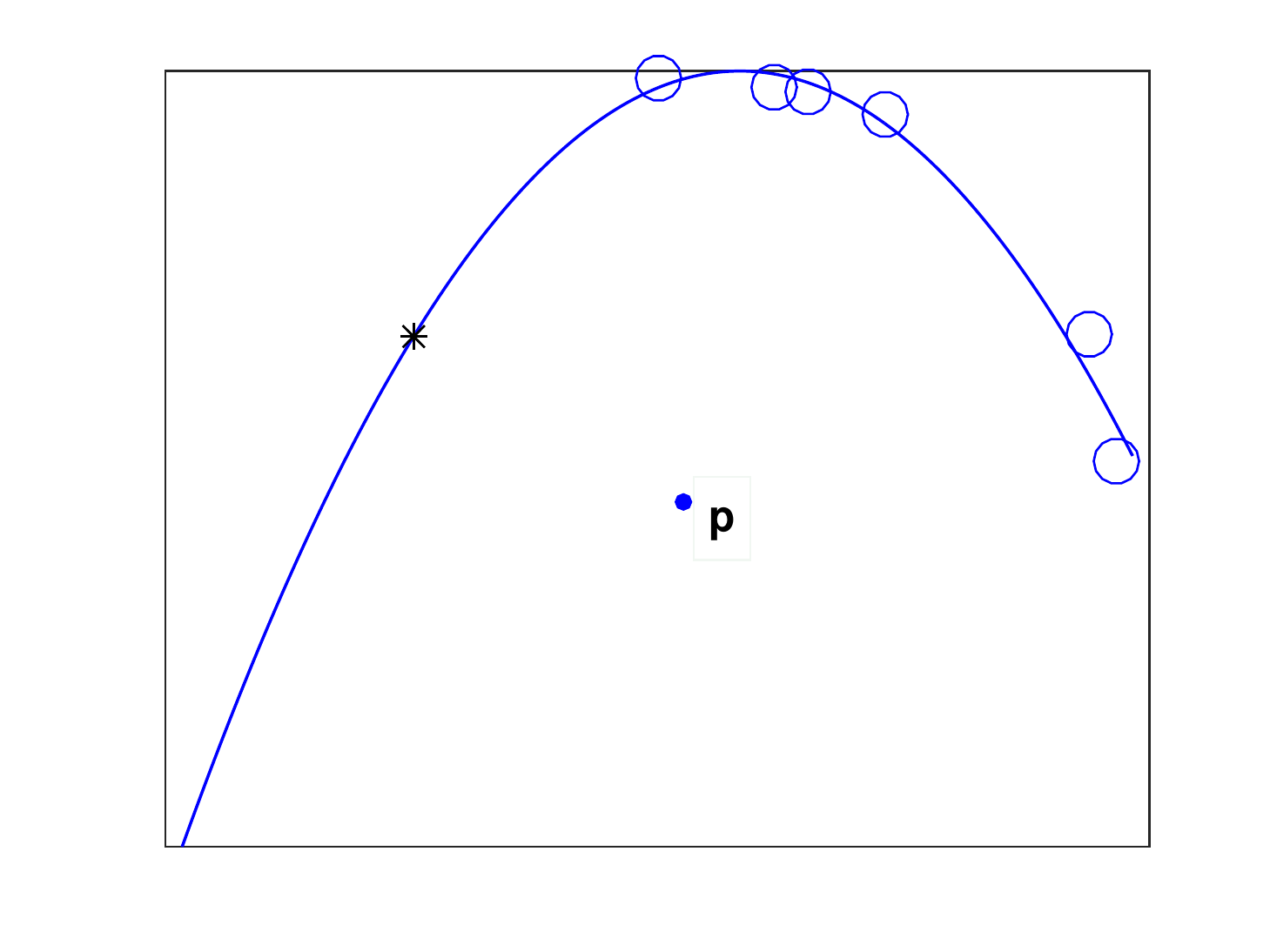}
    \caption{The local GBPM reconstruction (blue line) using $m=6$ footpoints (blue circles).  The minimizer of distance (black star) lies outside the interpolating points. This leads to the deactivation of grid point $\mathbf{p}$ (blue point). }
    \label{GBPMReconstruction}
\end{figure}

\begin{figure}
    \includegraphics[width=0.49\textwidth]{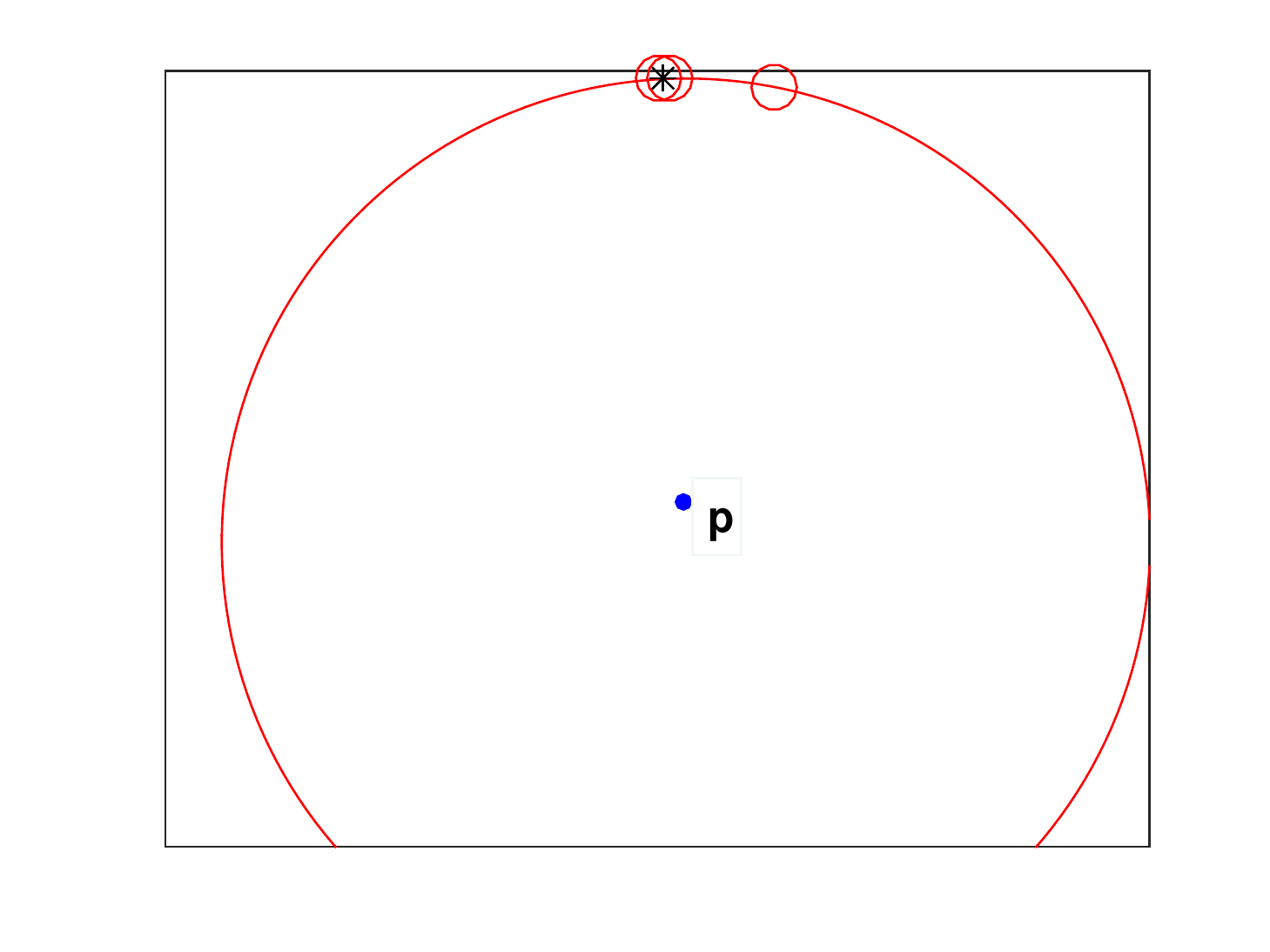}
    \includegraphics[width=0.49\textwidth]{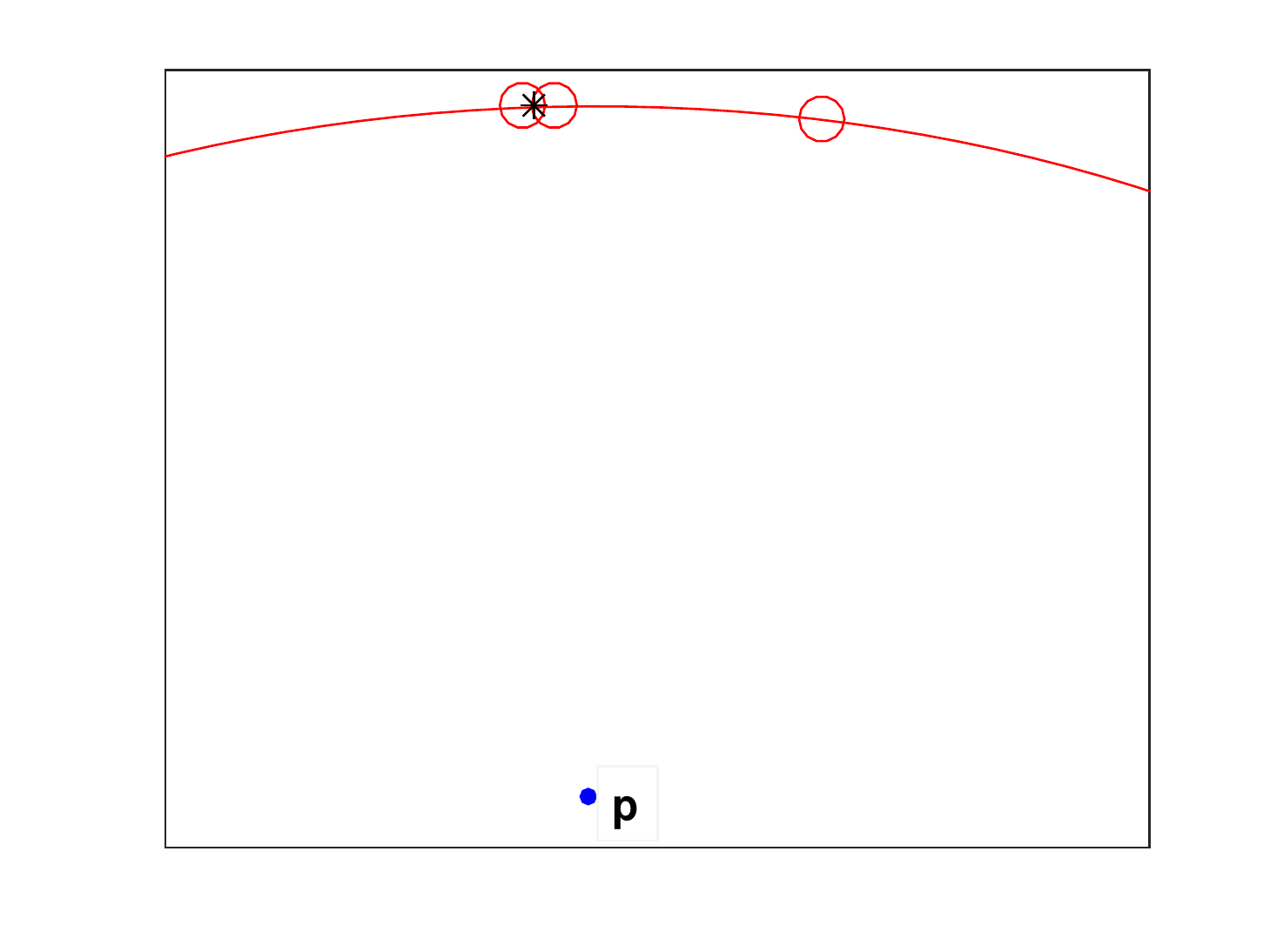}
    \caption{Left: The osculating circle reconstruction (red line) from three points (red circles) and the minimizer of distance (black star). The minimizer is chosen as the new footpoint. Right: A zoom-in of the osculating circle reconstruction.}
    \label{CircleReconstruction}
\end{figure}

We now detail this additional reconstruction step.
Consider first the reconstruction for the case of curves in 2D.
When a footpoint deactivation is flagged by the GBPM, we select three nonlinear points that are closest to the grid point $\mathbf{p}$ and
 fit the unique osculating circle through these points.
The update to the footpoint is taken to be the local minimizer of distance to $\mathbf{p}$ on the osculating circle.

A 2D illustration of the GBPM and the proposed modification appears in Figures~\ref{GBPMReconstruction} and \ref{CircleReconstruction}. Recall that for a grid point $\mathbf{p}$, the GBPM collects the $m$ closest footpoints and locally reconstructs the surface. The new footpoint is determined by minimizing the $L^2$-distance between the grid point $\mathbf{p}$ and the local reconstruction. However, in some cases, the GBPM rejects the new footpoint, and flags the point $\mathbf{p}$ for deactivation. Figure~\ref{GBPMReconstruction} gives an example where the grid point is flagged for deactivation because the minimizer (black star) lies outside the $m=6$ interpolating points. In our modified algorithm, an osculating circle based on the three points closest to $\mathbf{p}$ is constructed.  The closest point (black star) on the circle to {\bf p} gives the updated footpoint; see Figure~\ref{CircleReconstruction}.

For the case of surfaces, a similar strategy is employed. Specifically, whenever a grid point is flagged for deactivation by the GBPM, we reconstruct the surface by fitting a sphere to four non-coplanar footpoints. Spherical shapes have been shown to be effective for reconstructing surfaces moving by curvature-dependent speeds in Hon et al. \cite{hon2014cell}.  See  \cite{hon2014cell} for a method that uses circles and spheres to locally reconstruct curves and surfaces in the least squares sense.

For both cases (circles and spheres), it is necessary to impose the condition that the new footpoint lies at the same side as the interpolation points. Two candidate points are considered: the minimizer and its antipodal point on the circle/sphere. The point closest to the closest point of the reference grid point is the new footpoint.

With this proposed modification, all grid points within the computational tube remain active. We emphasize that the osculating circle and sphere reconstruction should not be used as a replacement for the GBPM. The GBPM is an effective, tested method for reconstructing a closest point representation. Only rarely, when grid point deactivation occurs, do we apply the osculating circle/spherical reconstruction.

\subsection{Numerical experiments for geometric motion}\label{RemarksTestSection}
In this section, we apply the modified GBPM to a number of examples to illustrate the performance and versatility of the method.
In view of our objective of coupling the method with the CPM, we select a  tube radius of $\gamma_{CPM}$ (see Equation (\ref{Bandwidth})).
Quadratic polynomials are used to locally reconstruct the surfaces in the least squares sense.
In two dimensions, $m=6$ footpoints with a minimum distance of $\delta=dx/4$ between one another are used,
while $m=20$ footpoints with $\delta=dx/2$ are used in three dimensions. The motion step is carried out by evolving the ODE system
$$\frac{d\mathbf{x}}{dt}=\mathbf{v}$$
with a step of forward Euler.

\subsubsection{Mean curvature motion on a circle and a sphere}

We begin by considering mean curvature motion on a circle and a sphere. The velocity is given by
 $$\mathbf{v} = -\kappa\mathbf{n}$$
 where $\kappa$ is the mean curvature and $\mathbf{n}$ is the outward unit normal vector. The analytical solution may be found by integrating
 $$\dot{R}=\kappa$$
 with $\kappa=1/R$ for the circle and $\kappa=2/R$ for the sphere. For all times $t\geq 0$, this gives
 $$R_{exact}(t) = \sqrt{R_0^2-2t}$$
 for the case of a circle and
 $$R_{exact}(t) = \sqrt{R_0^2-4t}$$
 for the sphere.

\begin{figure}
    \centering
    \includegraphics[width=0.49\textwidth]{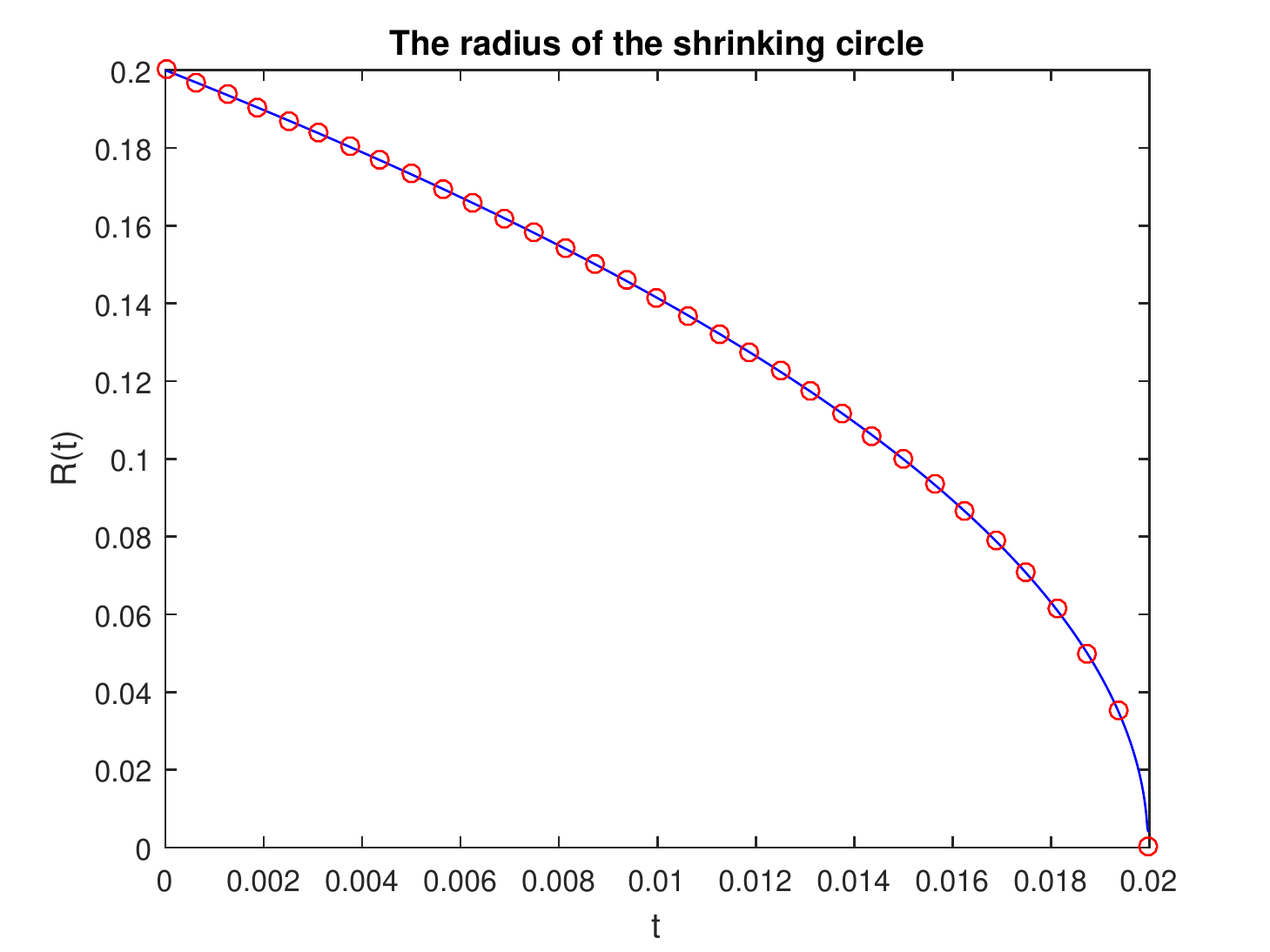}
    \includegraphics[width=0.49\textwidth]{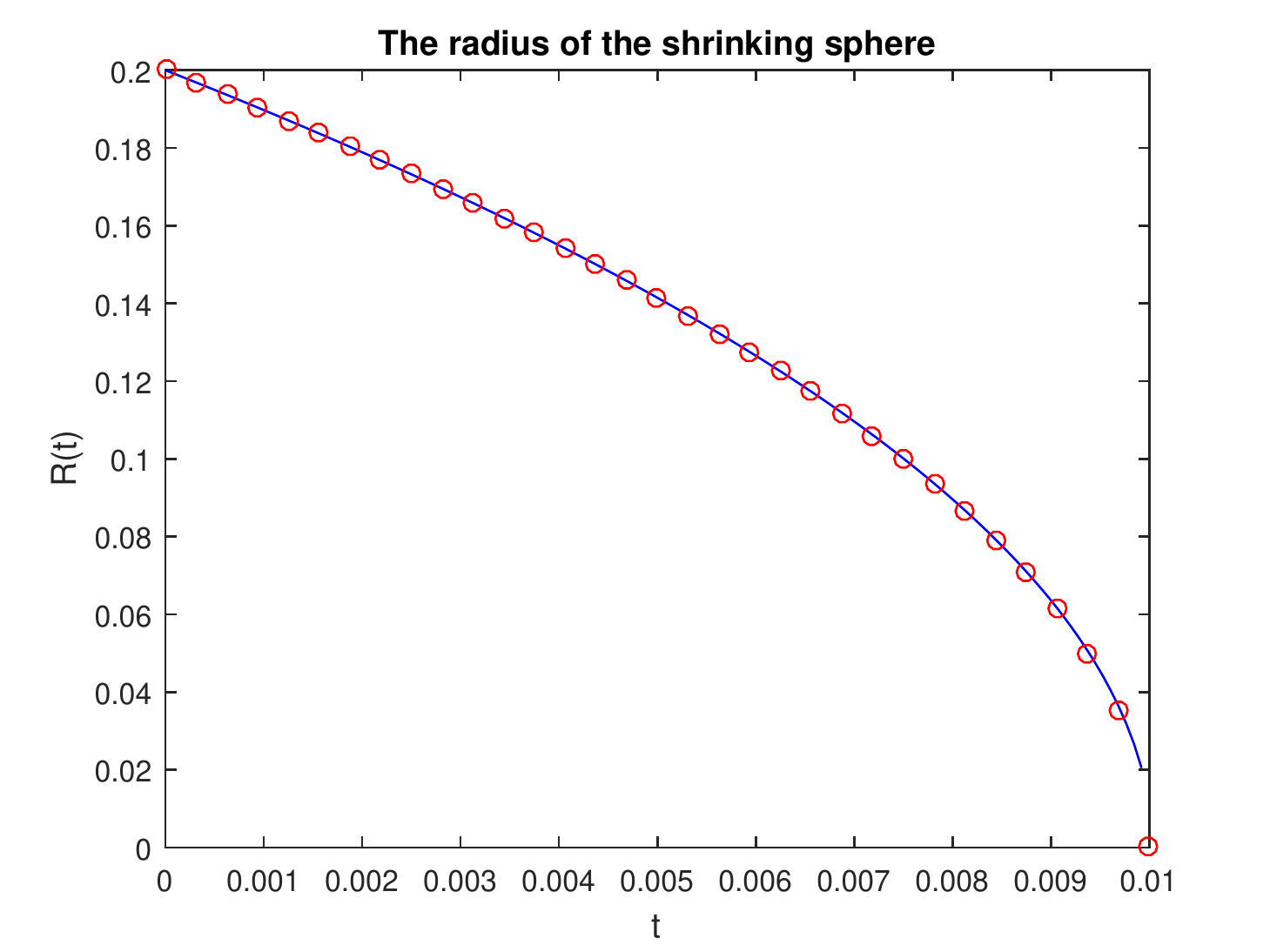}
    \caption{Blue lines give the computed radii as a function of time for the circle (left) and the sphere (right). Red dots mark the exact solution at selected times.}
    \label{Radius}
\end{figure}

 In our discretization, we select a spatial step-size of $\Delta x=0.00625$ for the circle and $\Delta x=0.0125$ for the sphere. A forward Euler time discretization is chosen with step-size $\Delta t=0.5\Delta x^2$. In both cases, the computed results give a good match to the exact solution. See Figure~\ref{Radius} for details.

\subsubsection{Vortex flow with rewind}

To test area conservation, we follow \cite{leveque1996high} and apply a vortex flow motion followed by a reversal of the velocity field. Starting from a circle centered at $(0.5,0.75)$ with radius $0.15$, we apply the velocity $\mathbf{v} = (v_1,v_2)$, where
$$\begin{array}{ll}
  v_1 = & -\sin(\pi x)^2\sin(2\pi y)\cos(\frac{\pi t}{t_{final}}),\\
  v_2 = & \sin(2\pi x)\sin(\pi y)^2\cos(\frac{\pi t}{t_{final}}),\\
\end{array}$$
and $t_{final}=4$. This \emph{vortex flow with rewind motion} yields the original surface at the prescribed final time $t_{final}$. In this example, different segments of the surface get close to one another, so a check on the consistency of the Lagrangian information is applied. Specifically, we require that
$$\mathbf{n}_0\cdot\mathbf{n}>\cos\Big(\frac{\pi}{2}\Big)$$
for every footpoint added in the collection of $m$ footpoints, where $\mathbf{n}_0$ is the unit normal vector of the closest footpoint to the reference grid point and $\mathbf{n}$ is the unit normal of the candidate footpoint.

\begin{figure}
    \centering
    \includegraphics[width=0.32\textwidth]{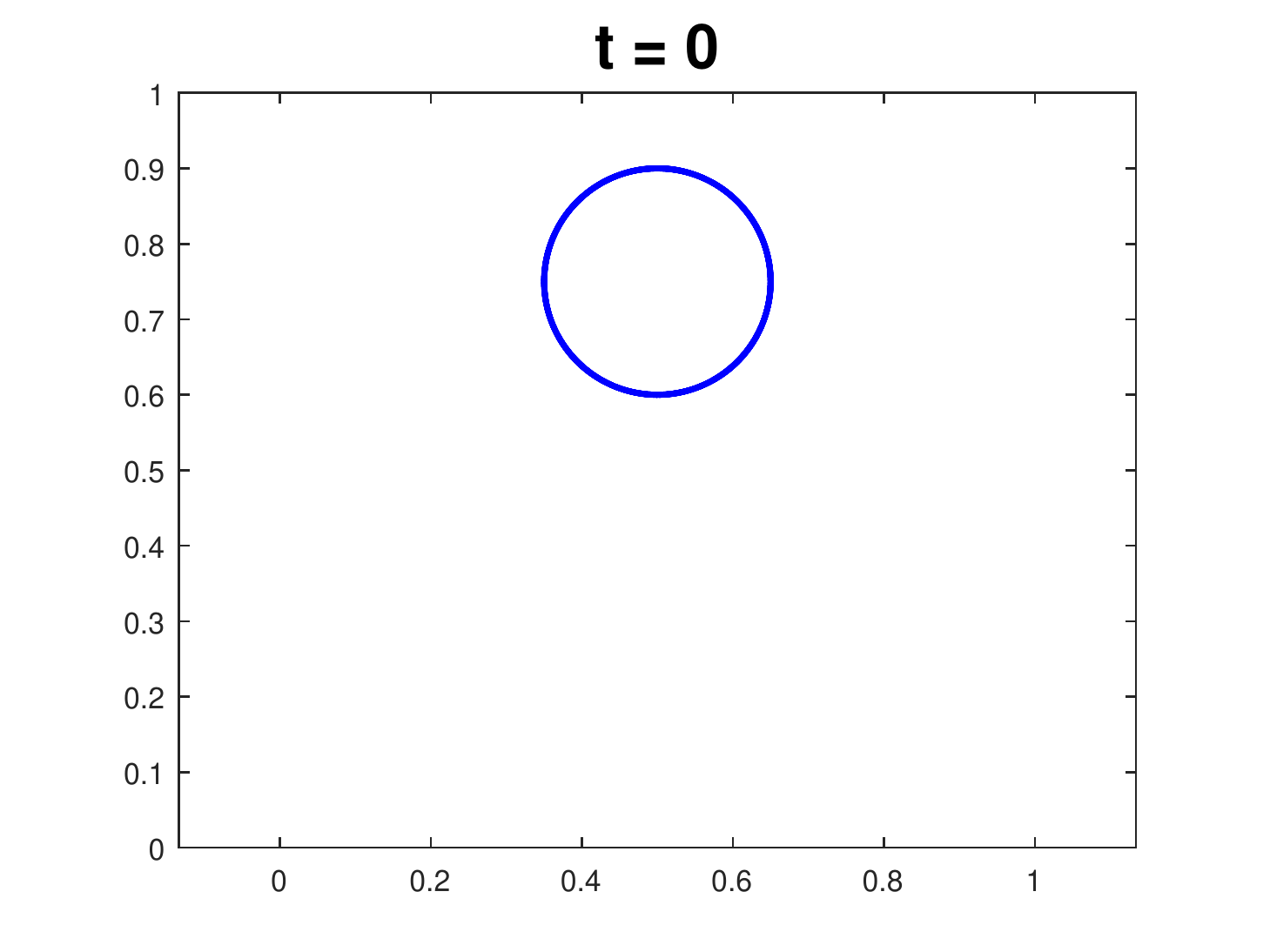}
    \includegraphics[width=0.32\textwidth]{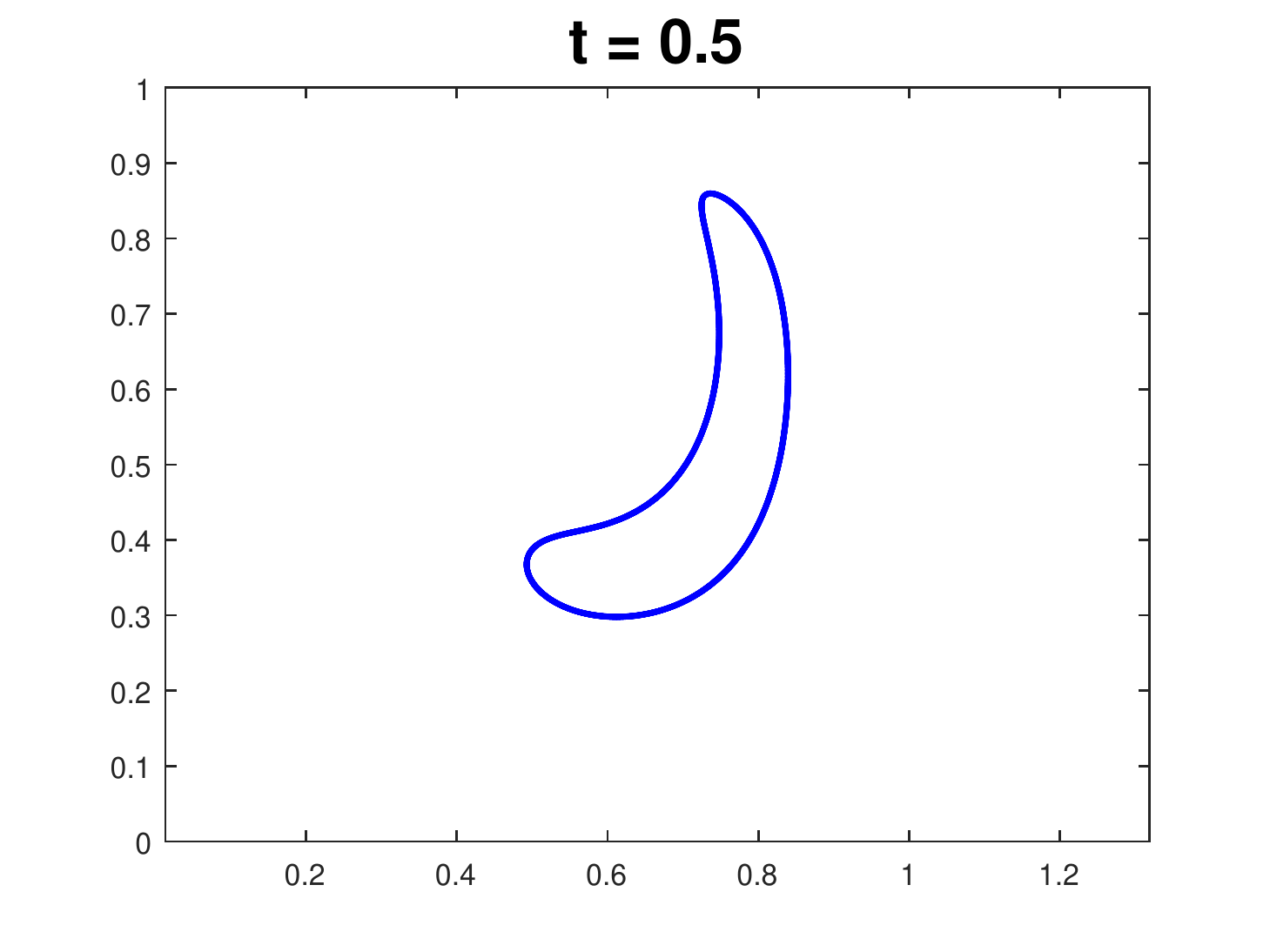}
    \includegraphics[width=0.32\textwidth]{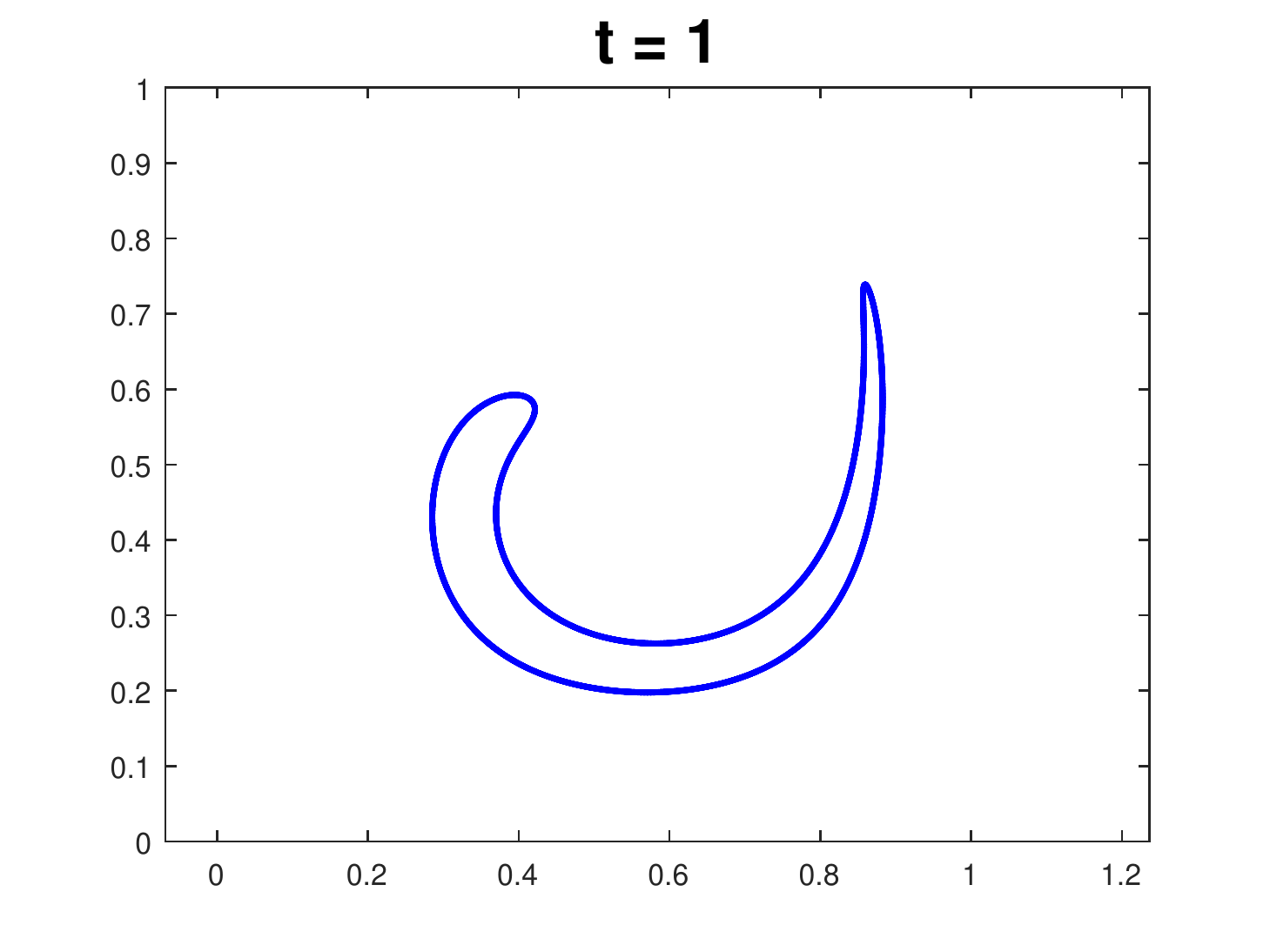}

    \includegraphics[width=0.32\textwidth]{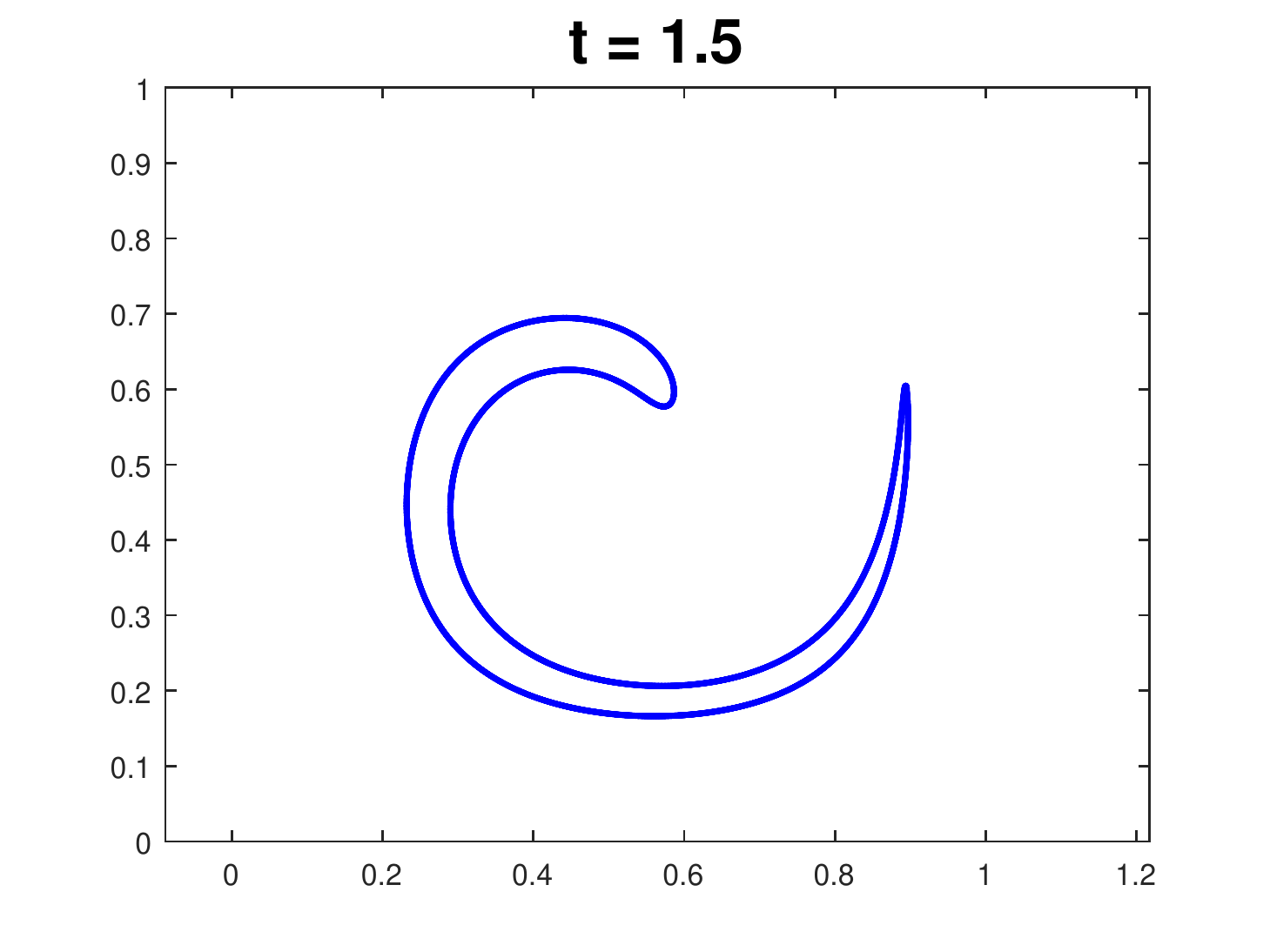}
    \includegraphics[width=0.32\textwidth]{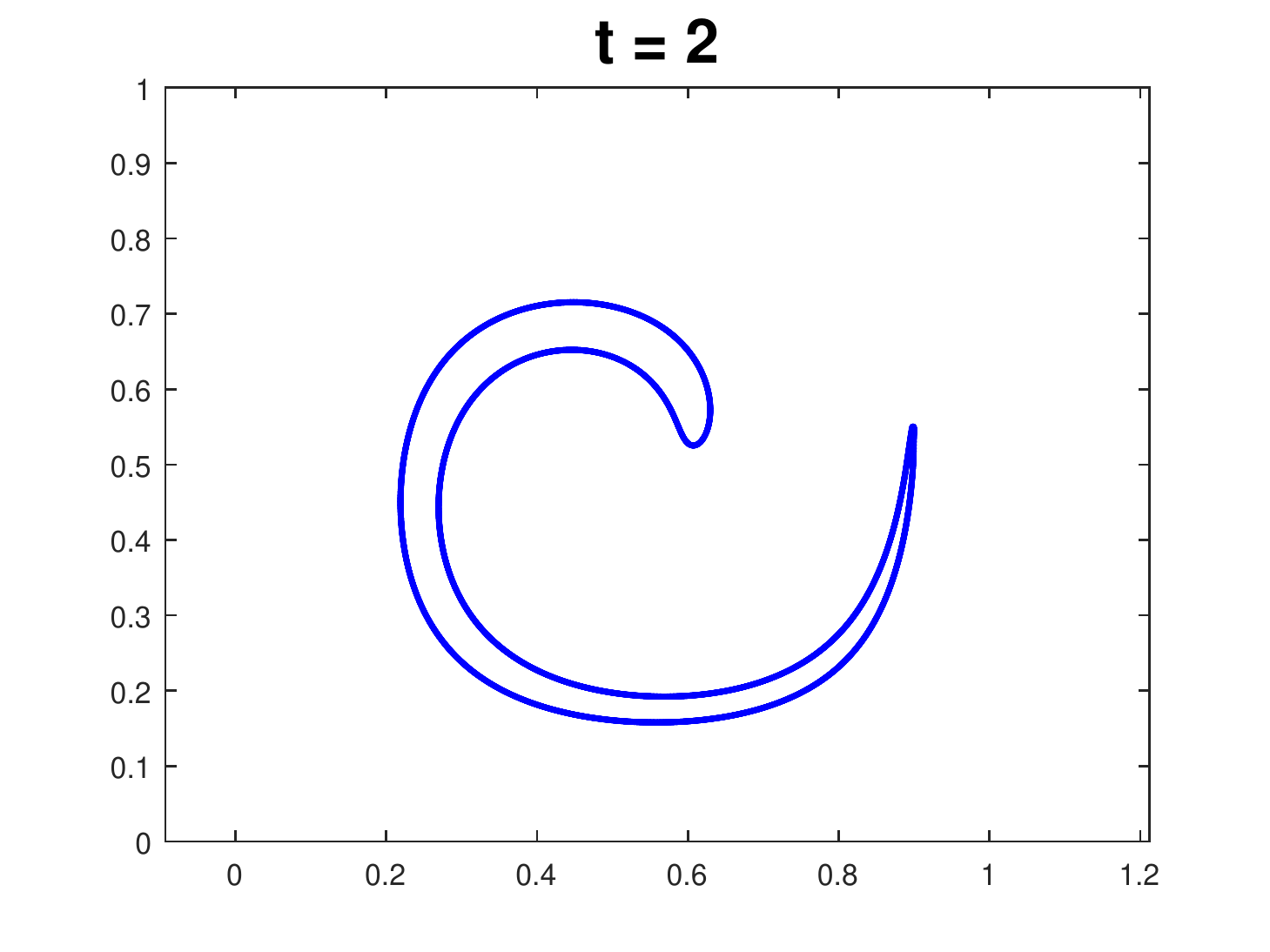}
    \includegraphics[width=0.32\textwidth]{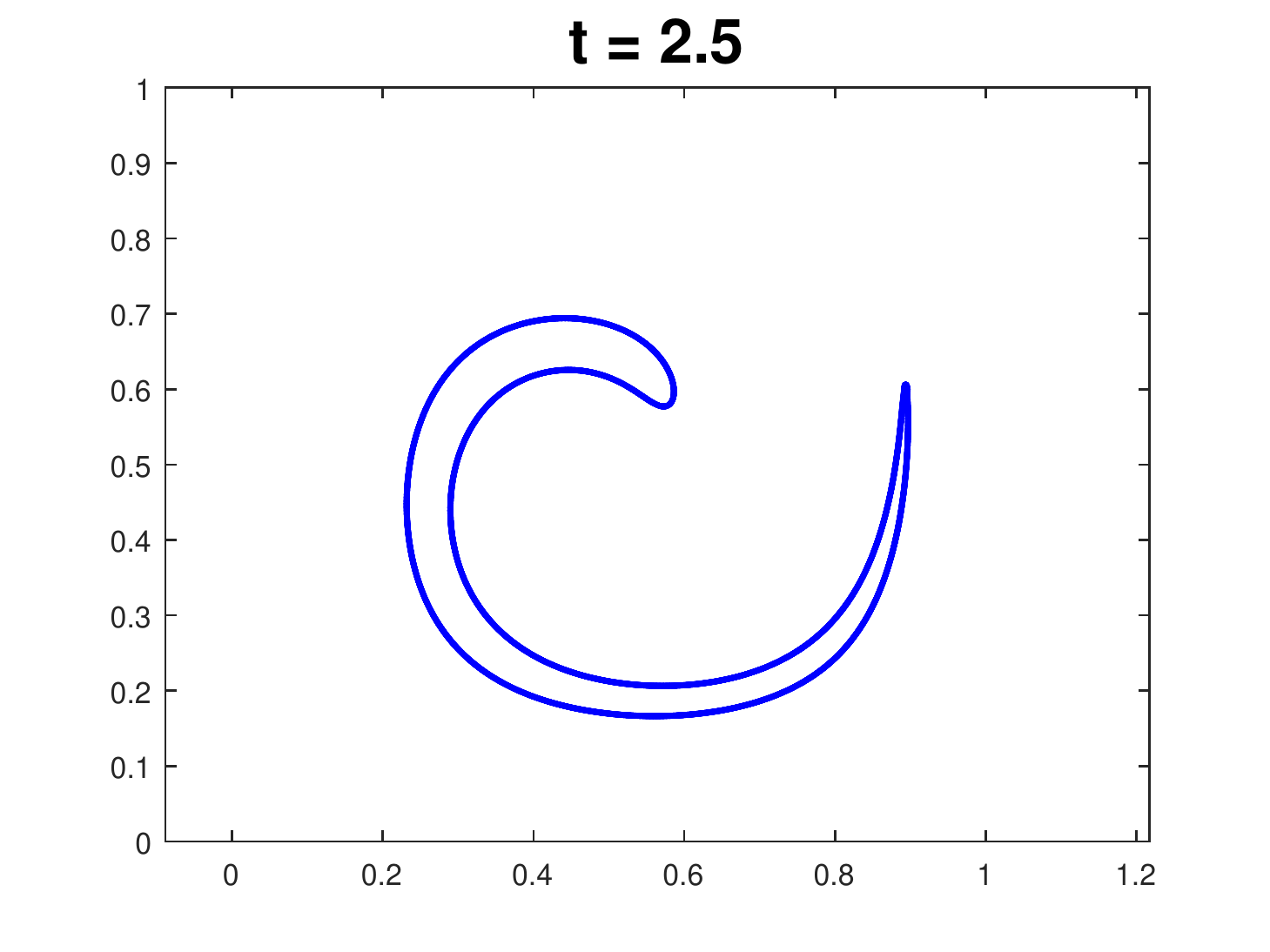}

    \includegraphics[width=0.32\textwidth]{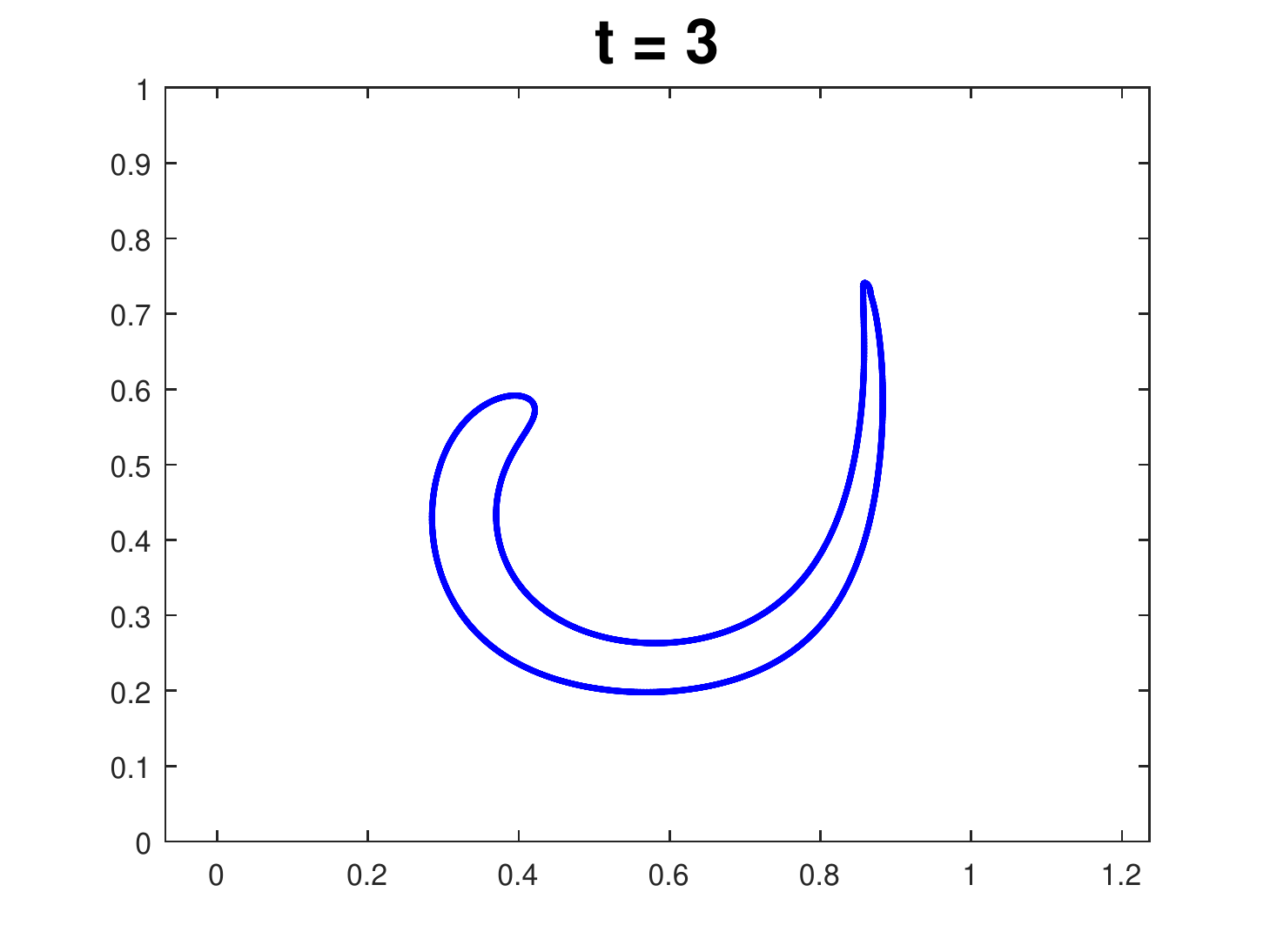}
    \includegraphics[width=0.32\textwidth]{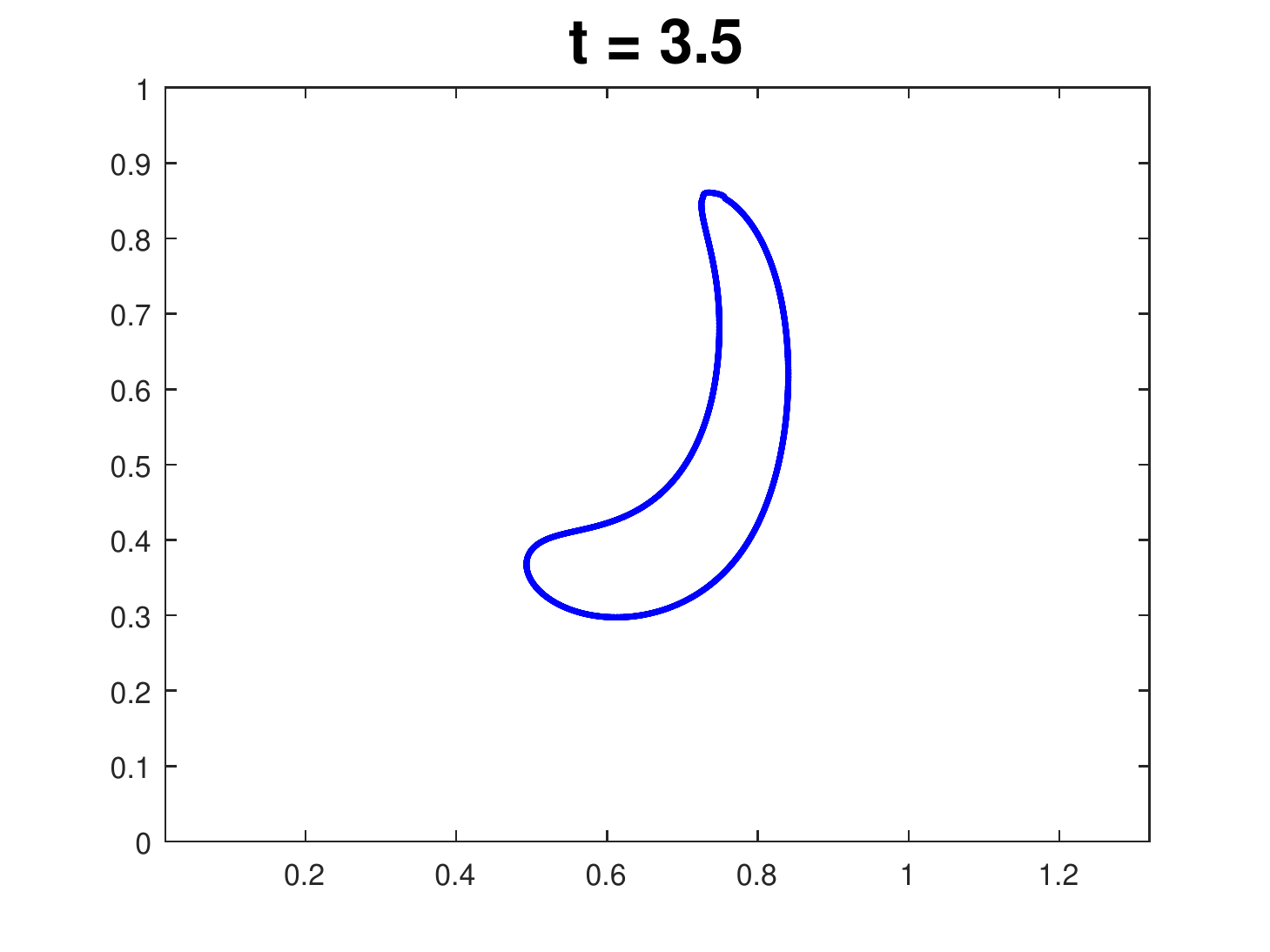}
    \includegraphics[width=0.32\textwidth]{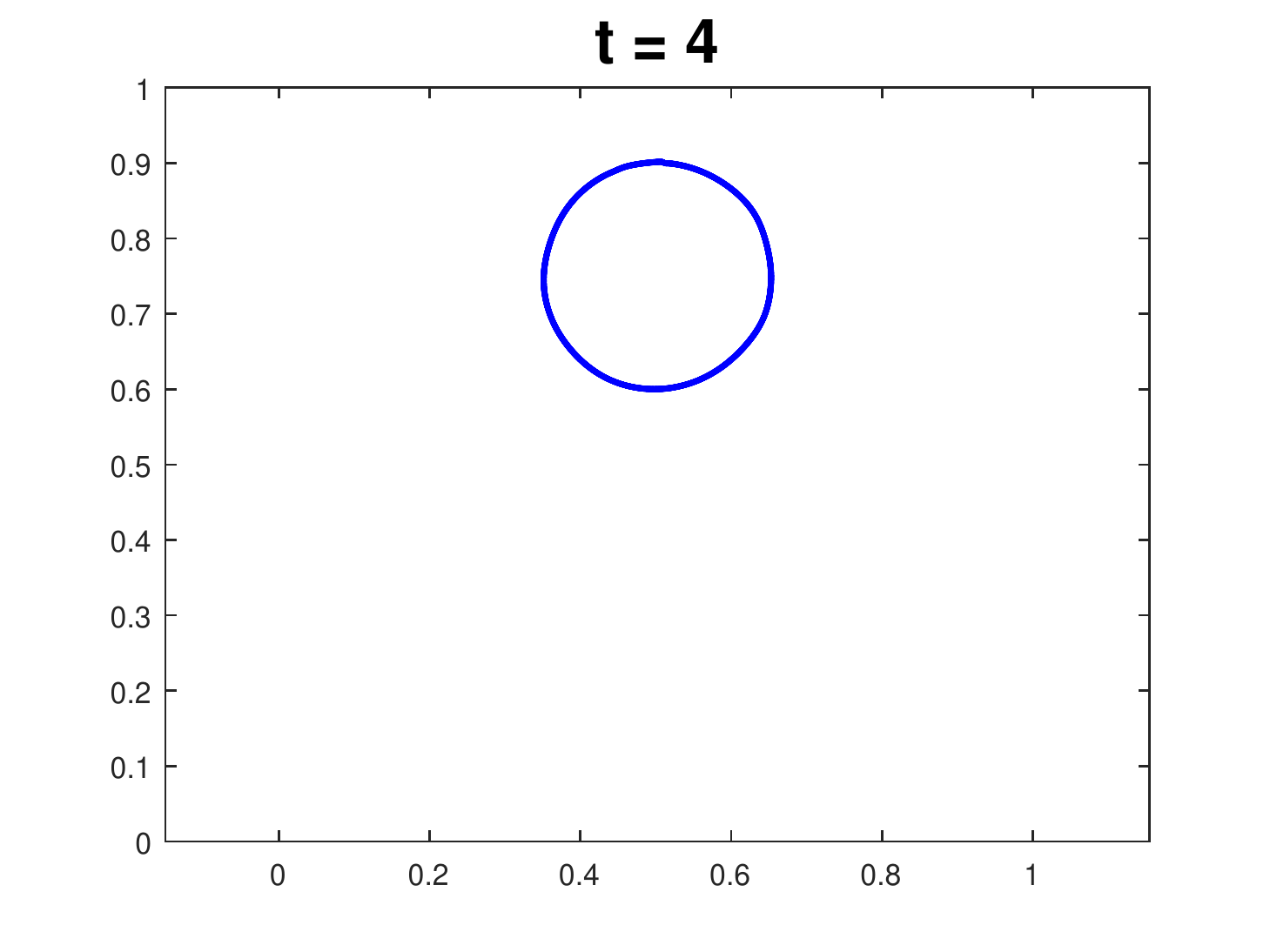}
    \caption{The evolution of a circle under vortex flow with rewind at various times $t$ (from top left to bottom right).}
    \label{VortexFlow}
\end{figure}

We compute to the final time $t_{final}=4$ using a grid spacing $\Delta x=0.0015625$ and a time step-size $\Delta t=0.8\Delta x$. Results at various times are displayed in Figure~\ref{VortexFlow}. A good visual match is observed between the initial and final contours. Indeed, the mean radius of the final contour is $R=0.1505$, which is a $0.33\%$ change from the initial circle.

\subsubsection{Merging circles}

We now consider an example that includes topological change. Starting from two circles with radius 0.15 centered at $(0.4,0.4)$ and $(0.6,0.6)$, we evolve outwards in the normal direction with constant unit speed. Under this flow, the two circles intersect and merge to form a single closed curve.

\begin{figure}
    \centering
    \includegraphics[width=0.49\textwidth]{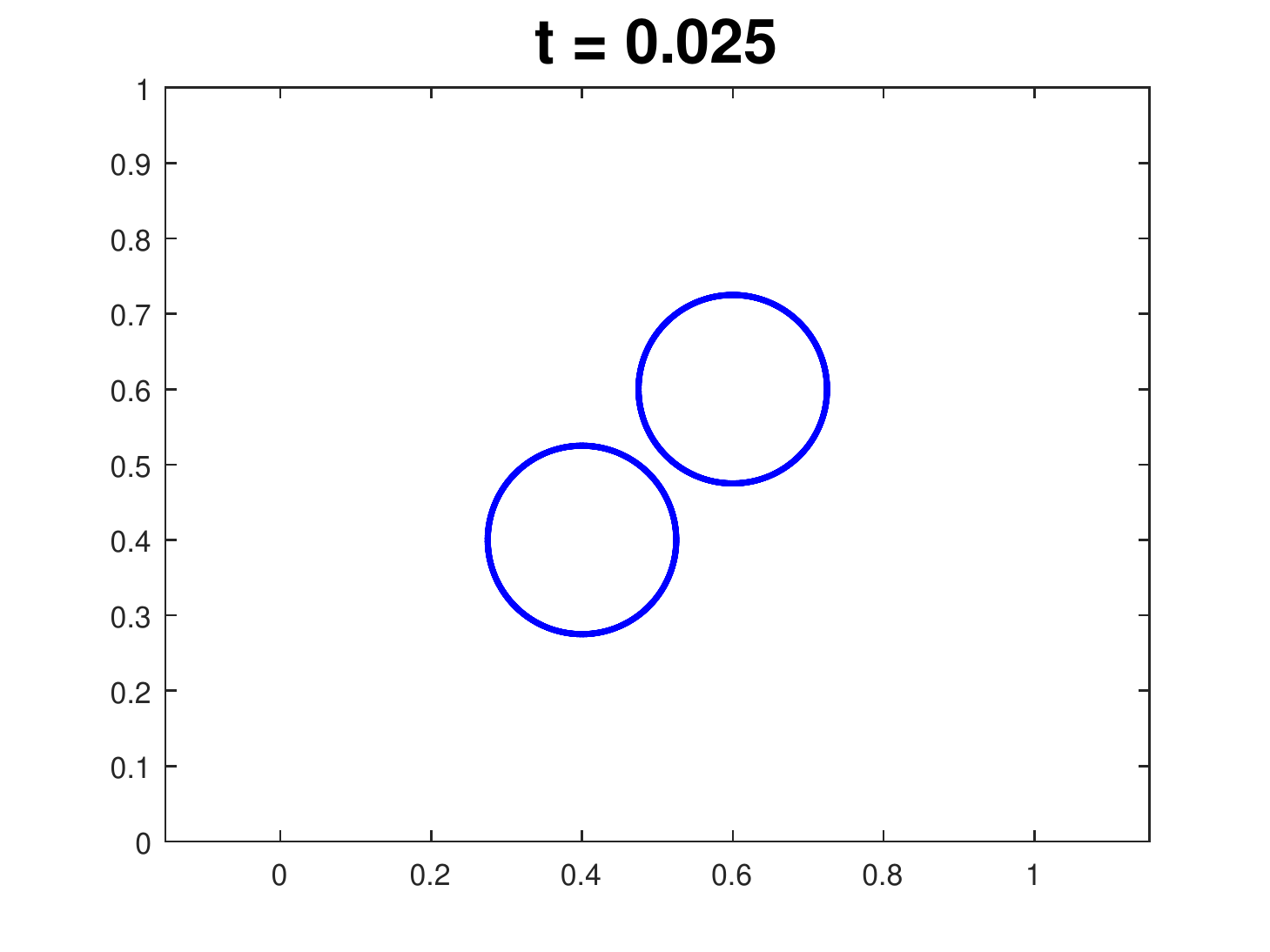}
    \includegraphics[width=0.49\textwidth]{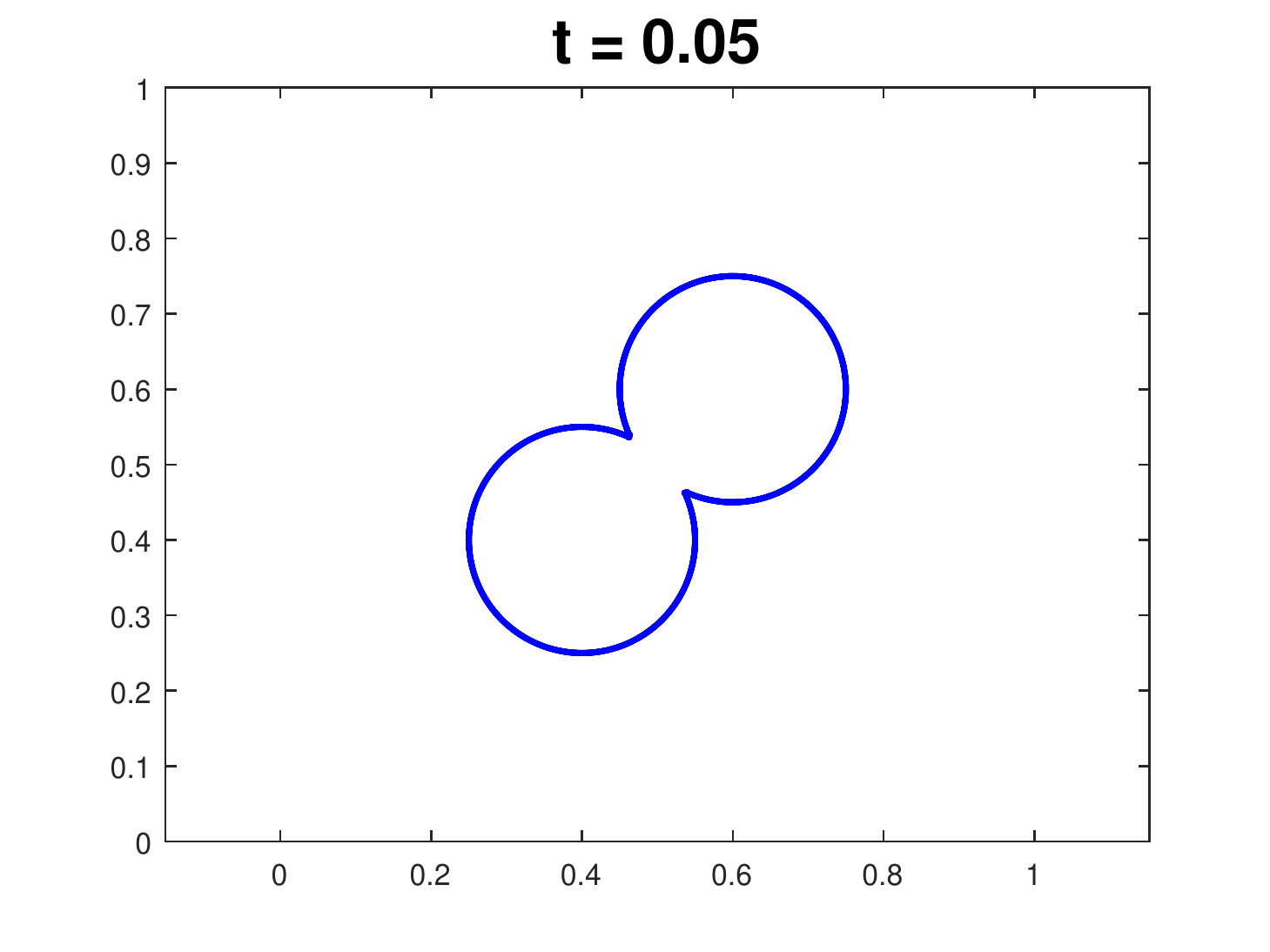}
    \includegraphics[width=0.49\textwidth]{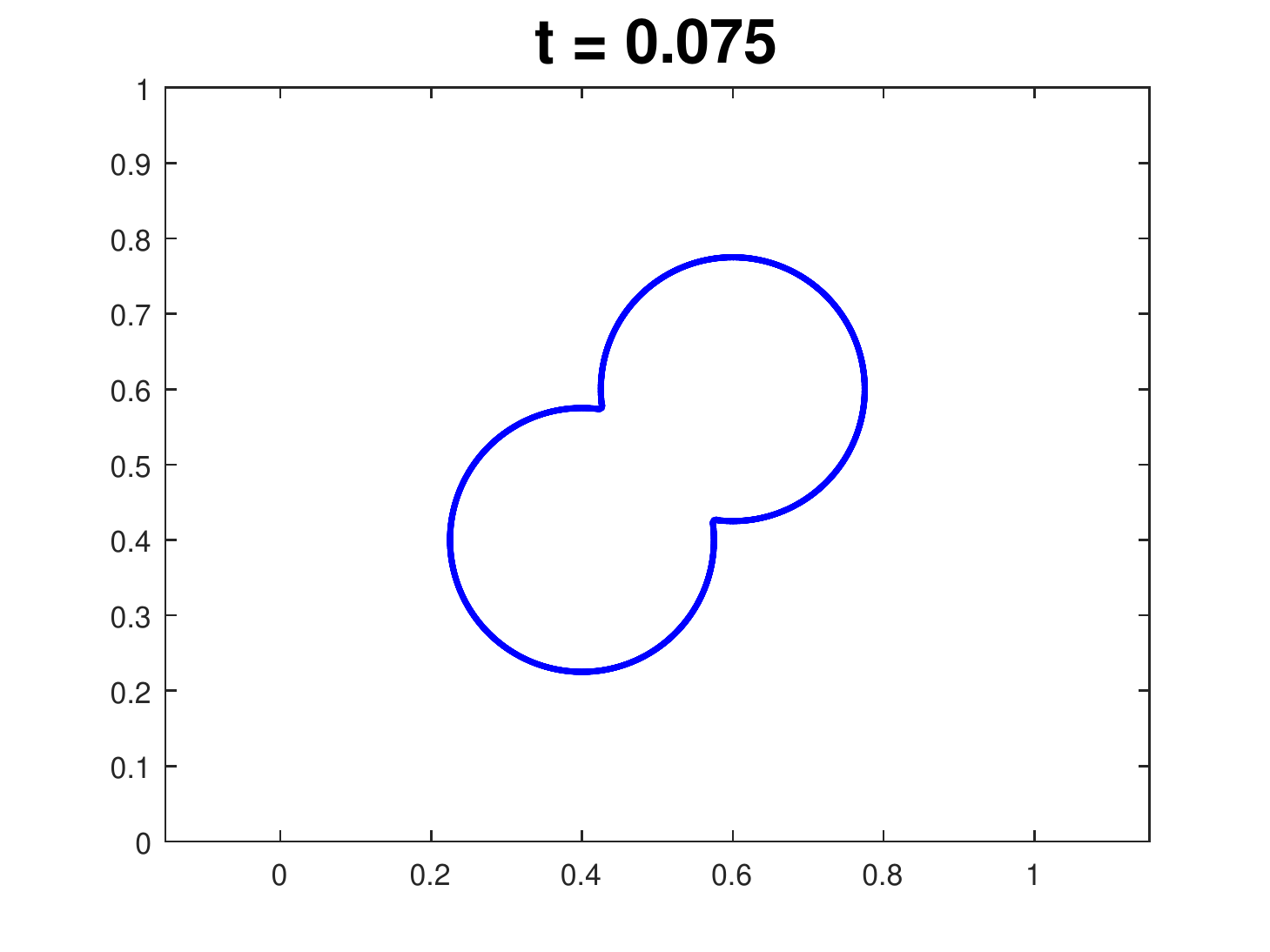}
    \includegraphics[width=0.49\textwidth]{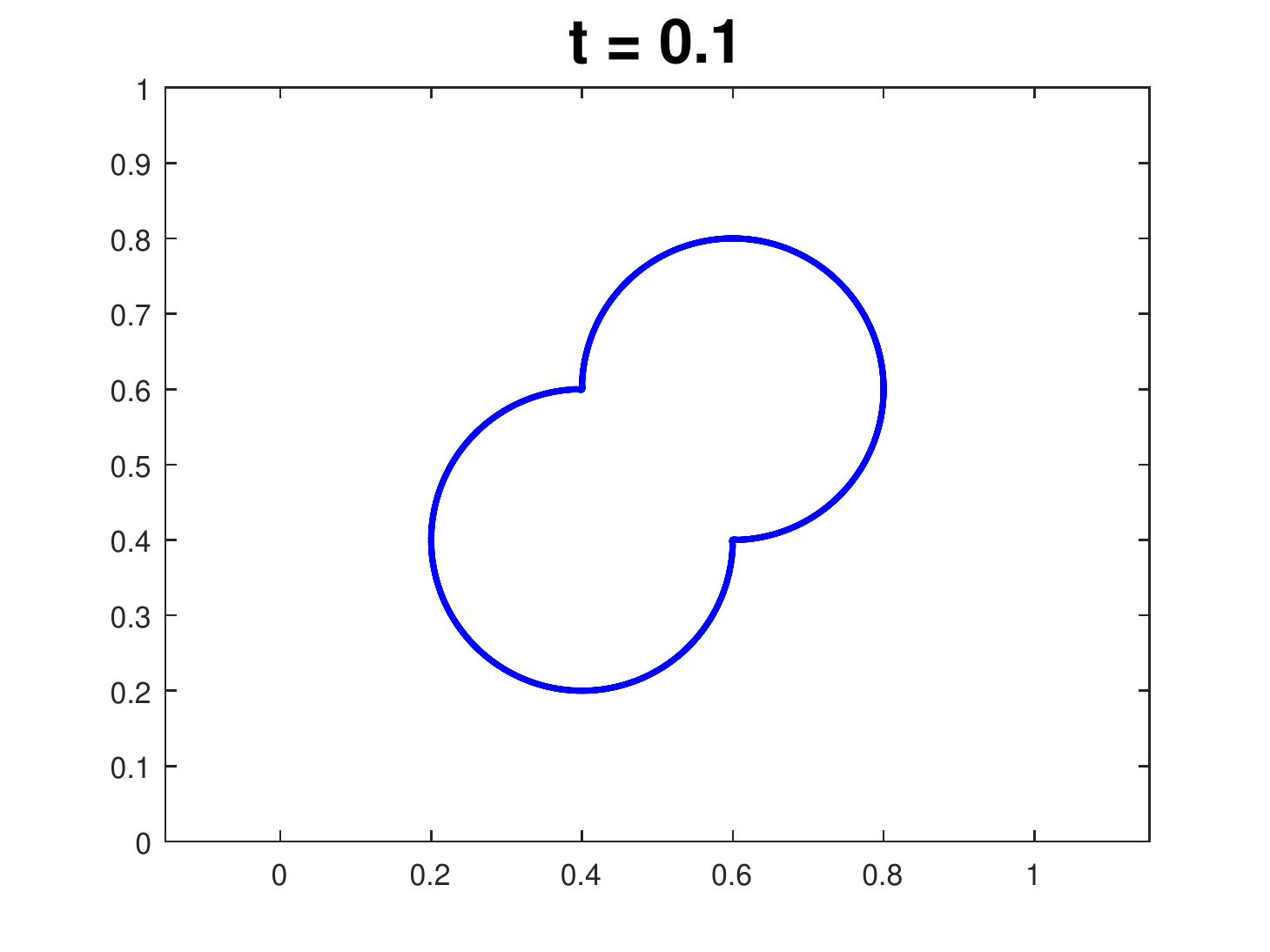}
    \caption{Two circles moving outwards with constant unit normal speed. A merger occurs.}
    \label{MergingCircles}
\end{figure}

Evolving to time $t_{final}=0.1$ using the modified GBPM with a mesh spacing $\Delta x=0.0015625$ and a time step-size $\Delta t=0.5\Delta x$ yields the results displayed in Figure~\ref{MergingCircles}. In agreement with \cite{leung2009grid}, we observe that the surfaces intersect to form a single curve which continues its outward normal expansion. In this example, it was necessary to deactivate grid points and their corresponding footpoints to ensure the merging of the surface. Specifically, if the selected $m$ footpoints have inconsistent Lagrangian information (i.e., $\mathbf{n_0}\cdot\mathbf{n}<\cos\Big(\frac{3\pi}{4}\Big)$), we deactivate the reference grid point and its footpoint; see Section~\ref{GBPMSection} for further details.

\subsubsection{Vortex flow with rewind in 3D}

Volume-conserving flows are also of considerable interest. We consider next vortex flow with rewind applied to an initial sphere with radius 0.15 centered at $(0.35,0.35,0.35)$. The velocity is given by $\mathbf{v}=(v_1,v_2,v_3)$, where
$$\begin{array}{ll}
  v_1 = & 2\sin(\pi x)^2\sin(2\pi y)\sin(2\pi z)\cos\Big(\frac{\pi t}{t_{final}}\Big),\\
  v_2 = & -\sin(2\pi x)\sin(\pi y)^2\sin(2\pi z)\cos\Big(\frac{\pi t}{t_{final}}\Big),\\
  v_3 = & -\sin(2\pi x)\sin(2\pi y)\sin(\pi z)^2\cos\Big(\frac{\pi t}{t_{final}}\Big),\\
\end{array}$$
and $t_{final}=1.5$. In this example, different segments of the surface get close to one another so Lagrangian consistency information is used.
Similar to vortex flow with rewind in 2D, we insist that
$$\mathbf{n}_0\cdot\mathbf{n}>\cos\Big(\frac{\pi}{3}\Big)$$
for every footpoint added in the collection of $m$ footpoints, where $\mathbf{n}_0$ is the unit normal vector of the closest footpoint to the reference grid point and $\mathbf{n}$ is the unit normal of the candidate footpoint.

\begin{figure}
    \centering
    \includegraphics[width=0.32\textwidth]{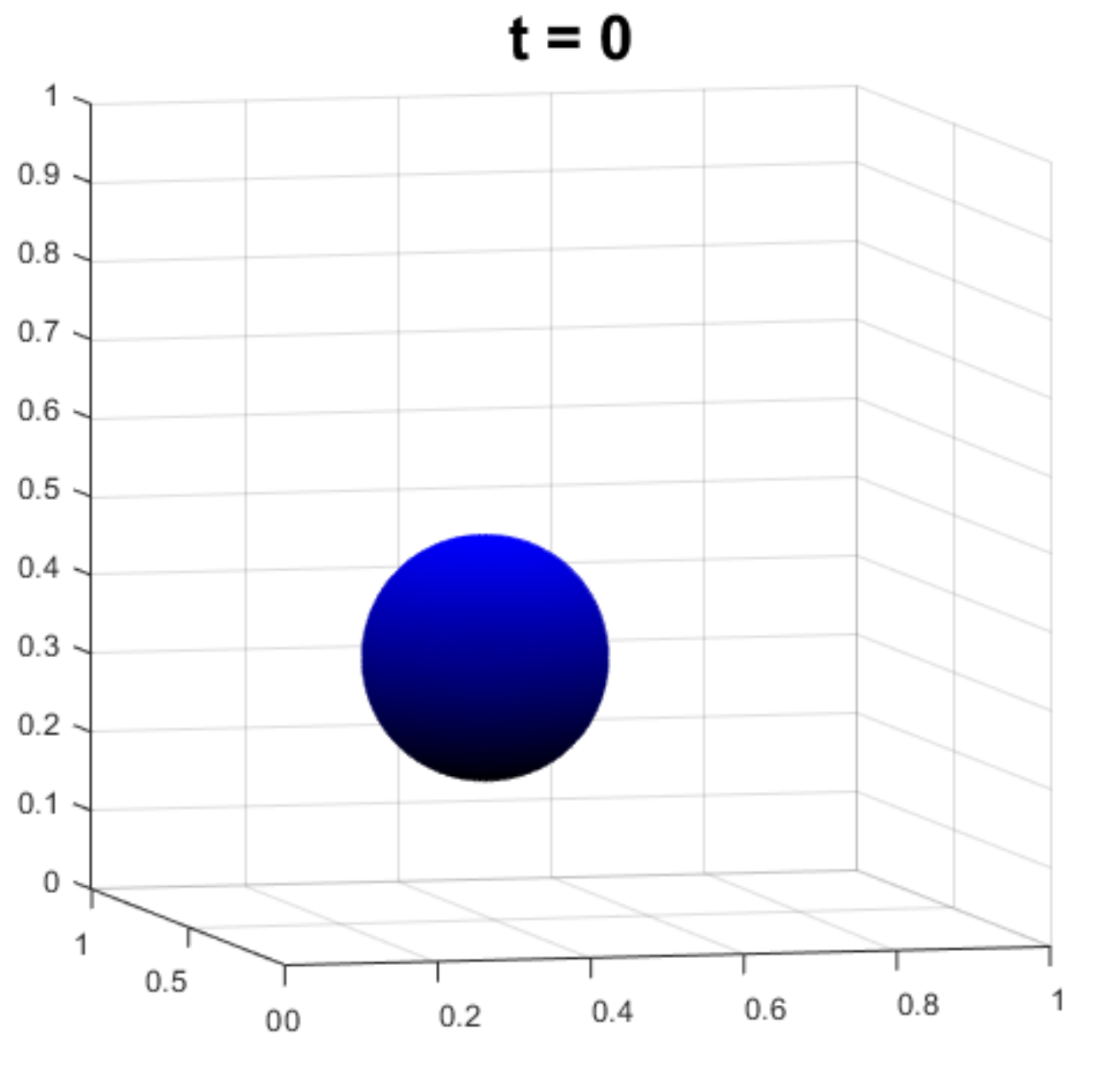}
    \includegraphics[width=0.32\textwidth]{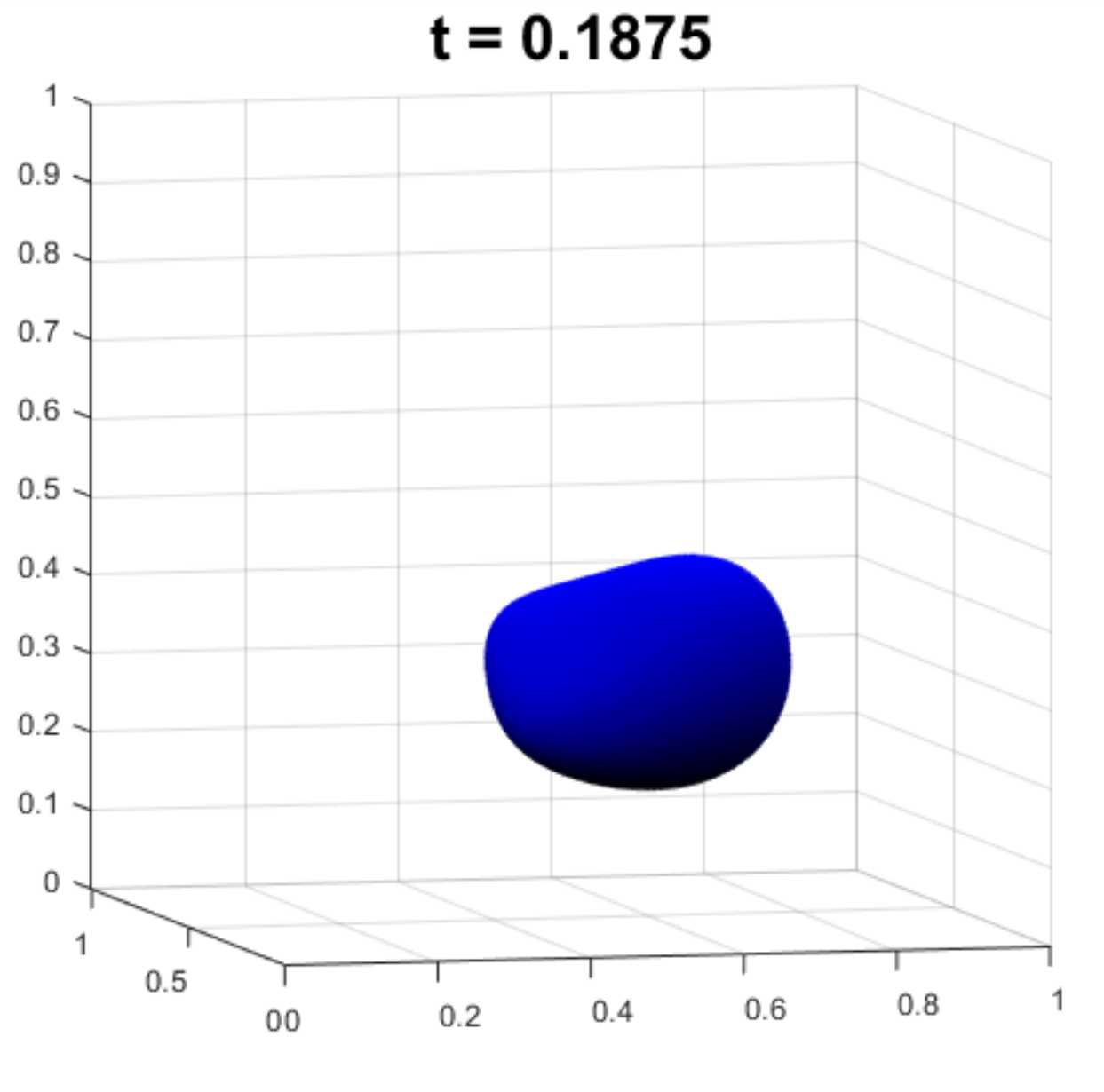}
    \includegraphics[width=0.32\textwidth]{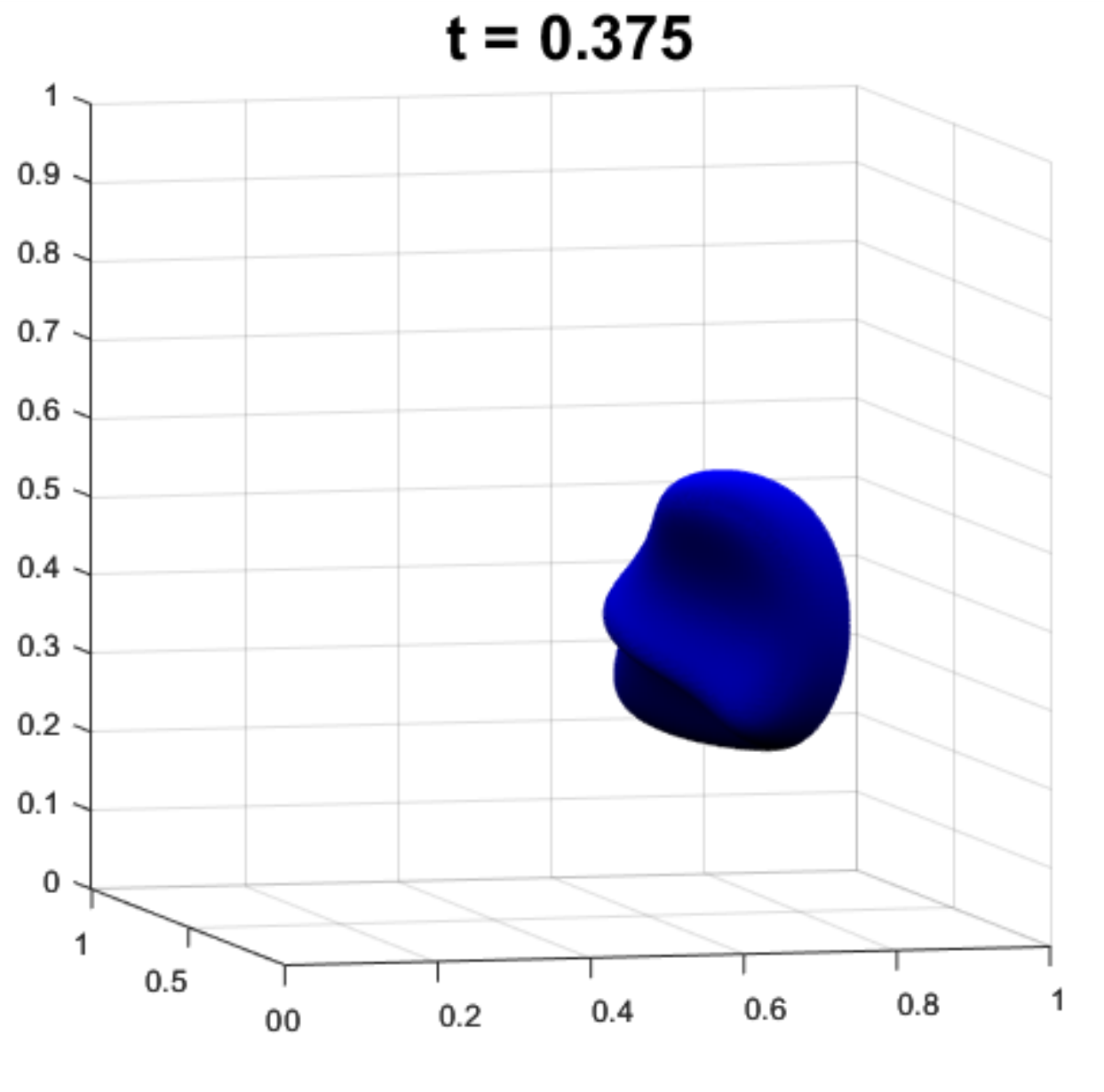}

    \includegraphics[width=0.32\textwidth]{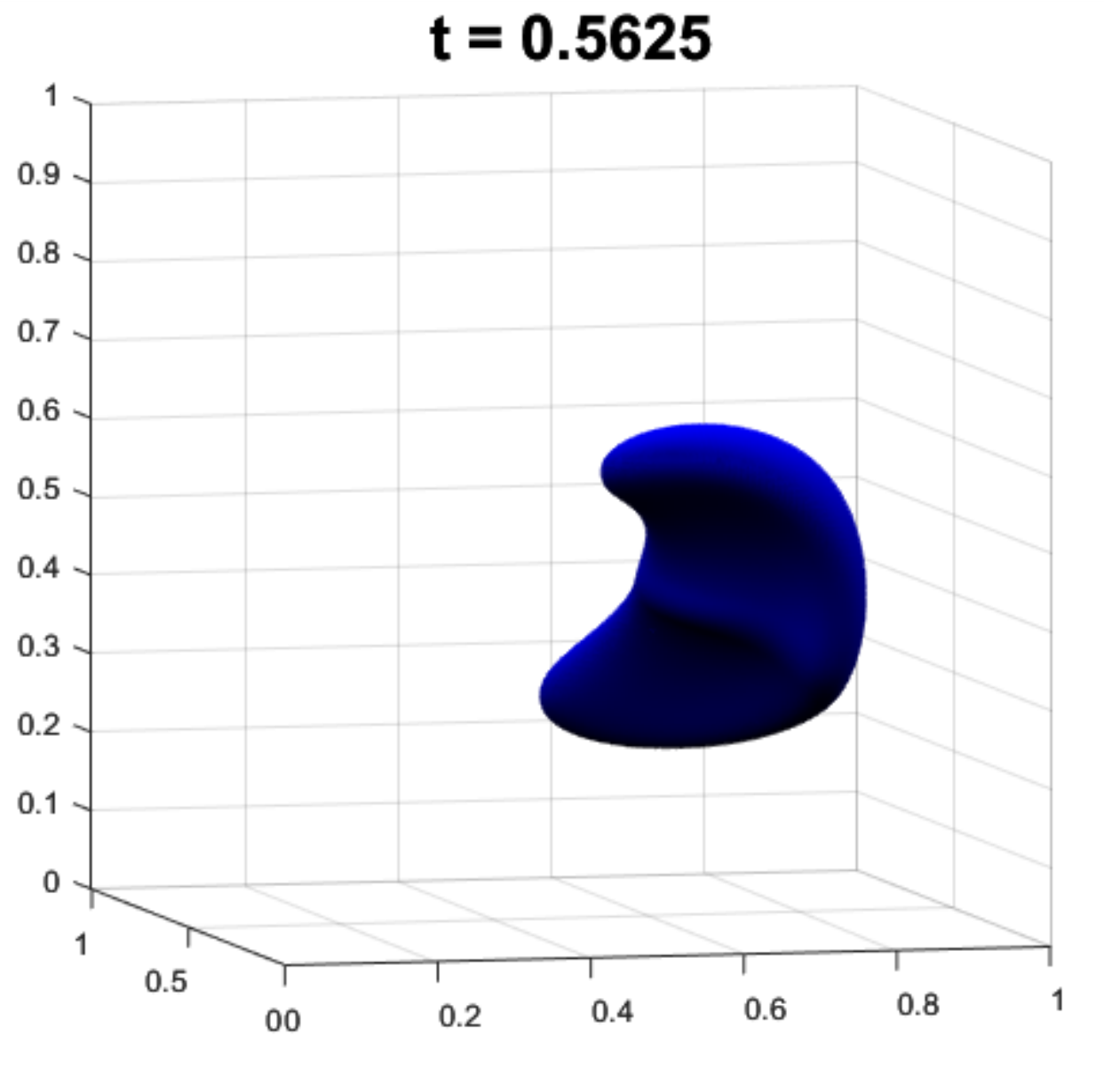}
    \includegraphics[width=0.32\textwidth]{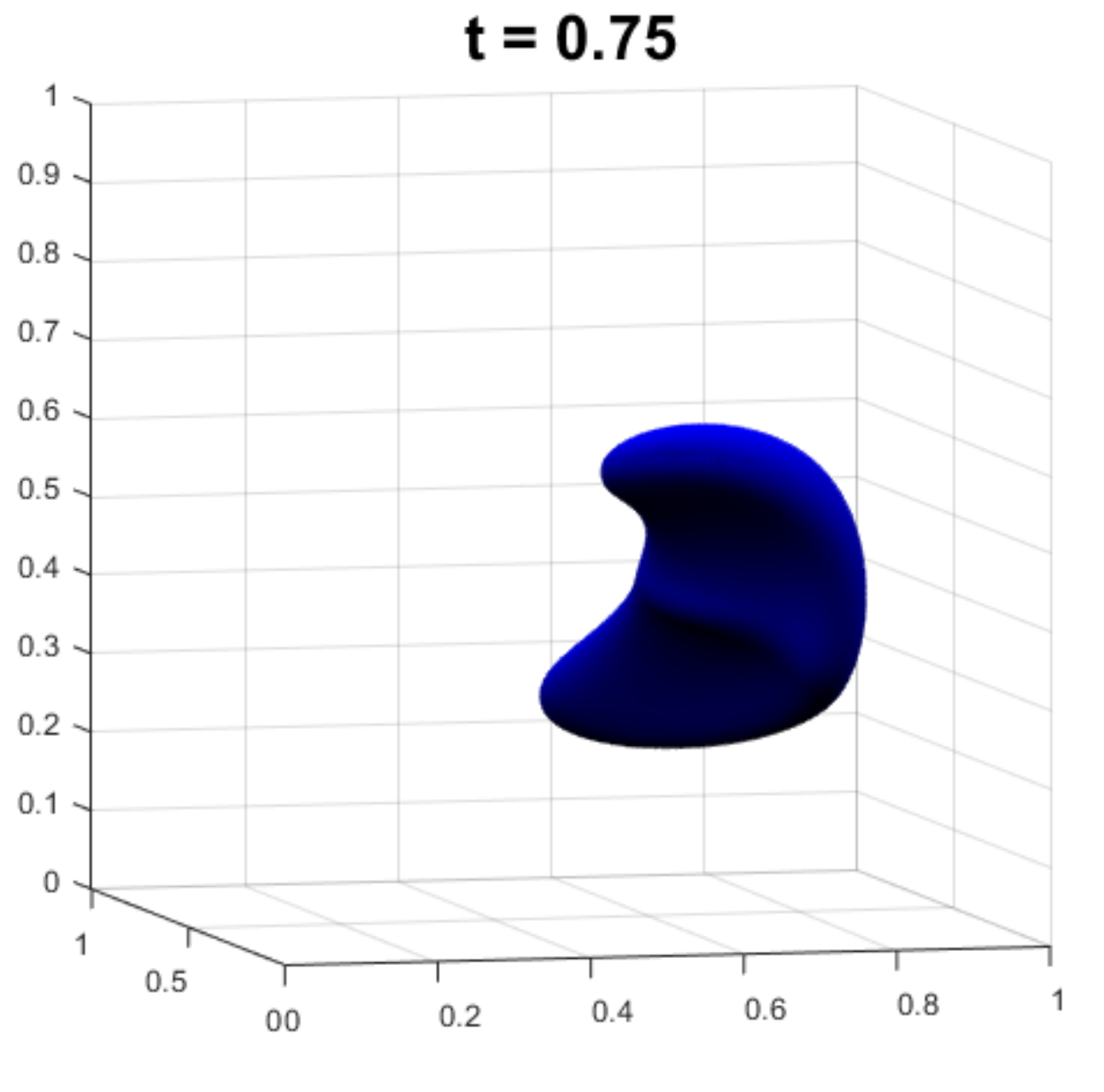}
    \includegraphics[width=0.32\textwidth]{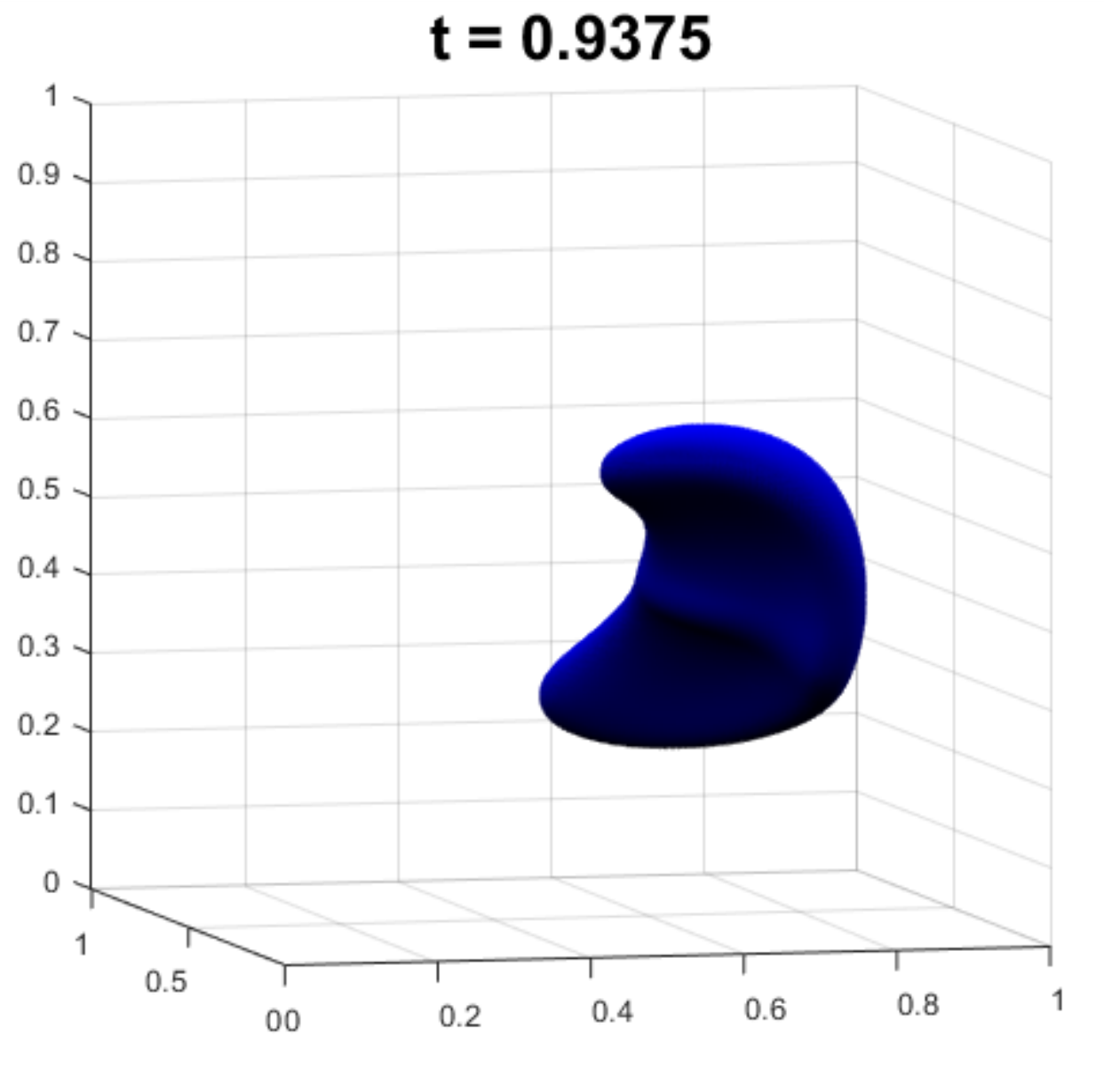}

    \includegraphics[width=0.32\textwidth]{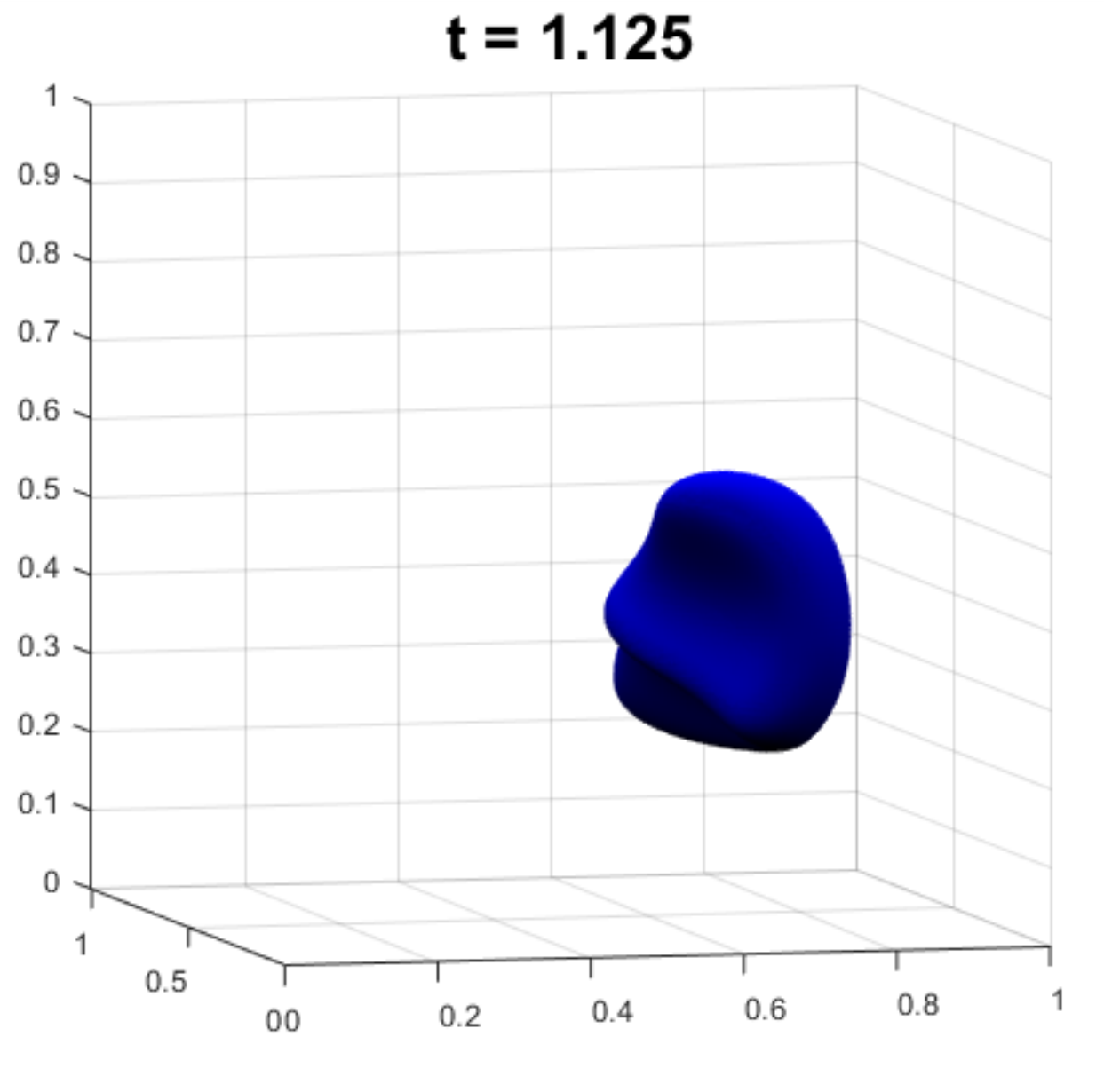}
    \includegraphics[width=0.32\textwidth]{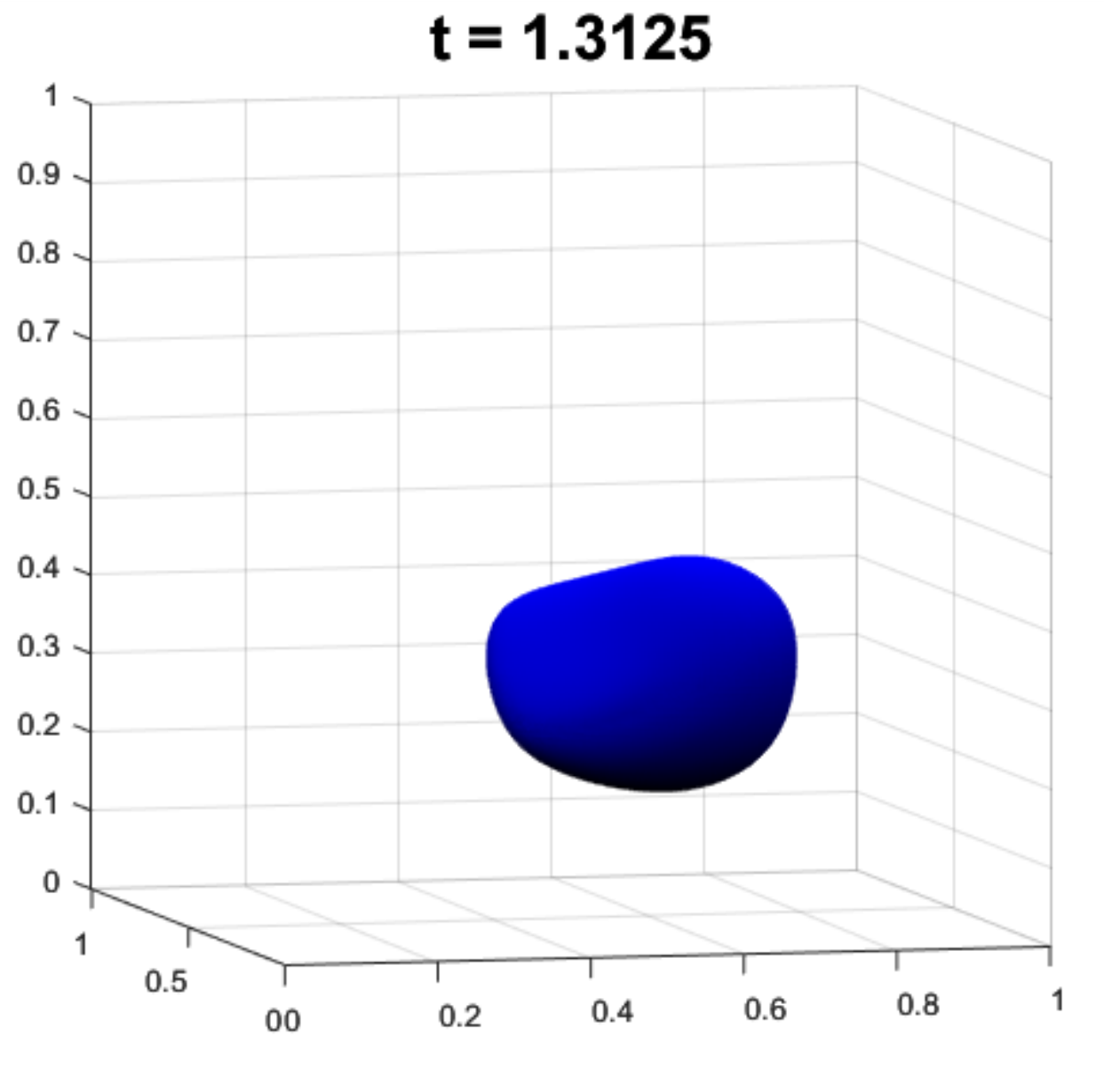}
    \includegraphics[width=0.32\textwidth]{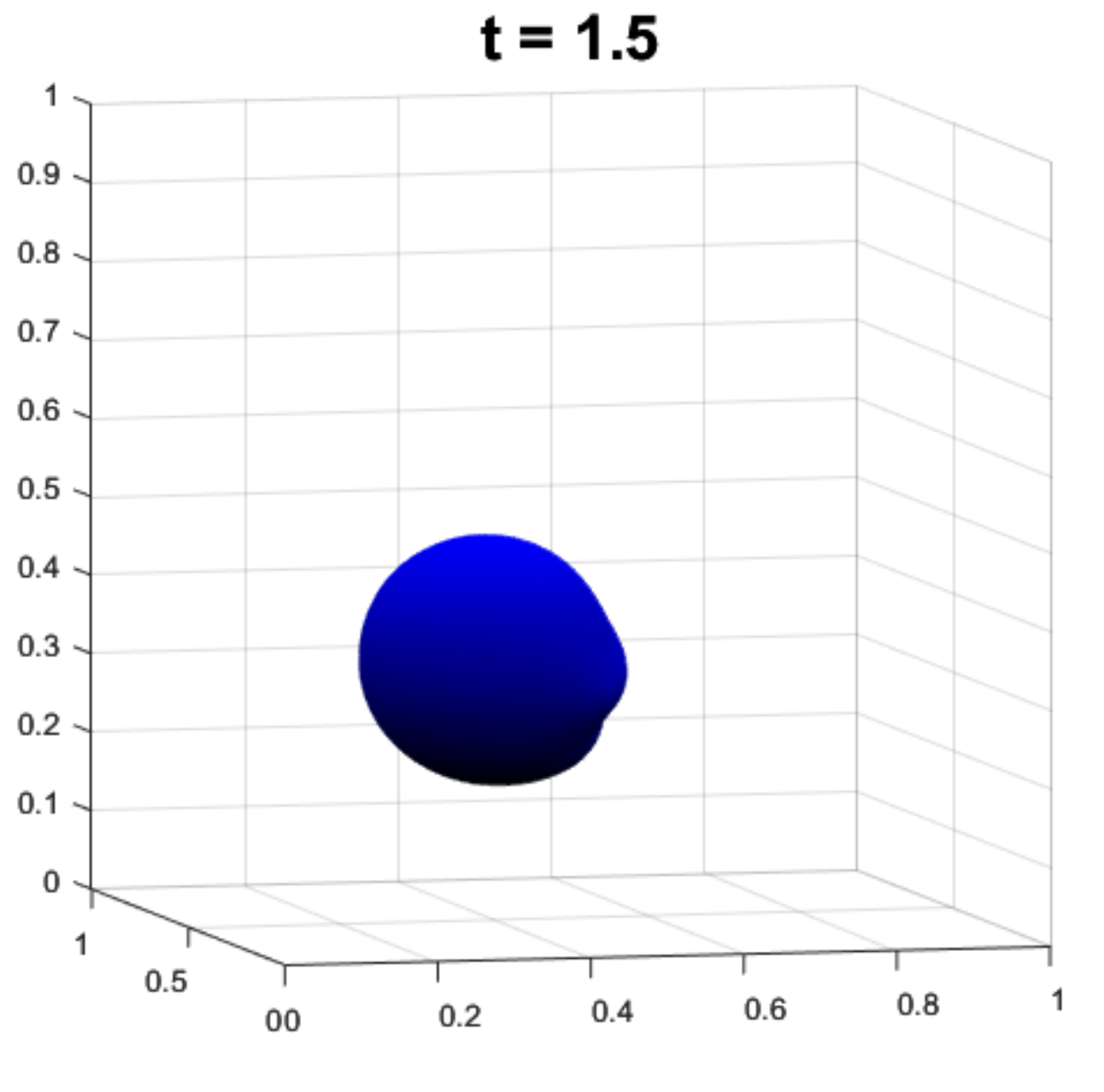}
    \caption{The evolution of a sphere under vortex flow with rewind at selected times $t$. The final time is $t_{final}=1.5$.}
    \label{VortexFlow3D}
\end{figure}

Selecting a mesh spacing $\Delta x=0.0083$ and a time step-size $\Delta t=0.8\Delta x$, we obtain the evolution displayed in Figure~\ref{VortexFlow3D}. We find that our final shape ($t=1.5$) is in good agreement with the initial sphere ($t=0$). The mean radius of the final shape is $R=0.1534$ with a standard deviation of $0.0048$.

\subsubsection{Topological change in 3D}

We conclude our tests for geometric motion with an example exhibiting topological change in three dimensions.

\begin{figure}
    \centering
    \includegraphics[width=0.49\textwidth]{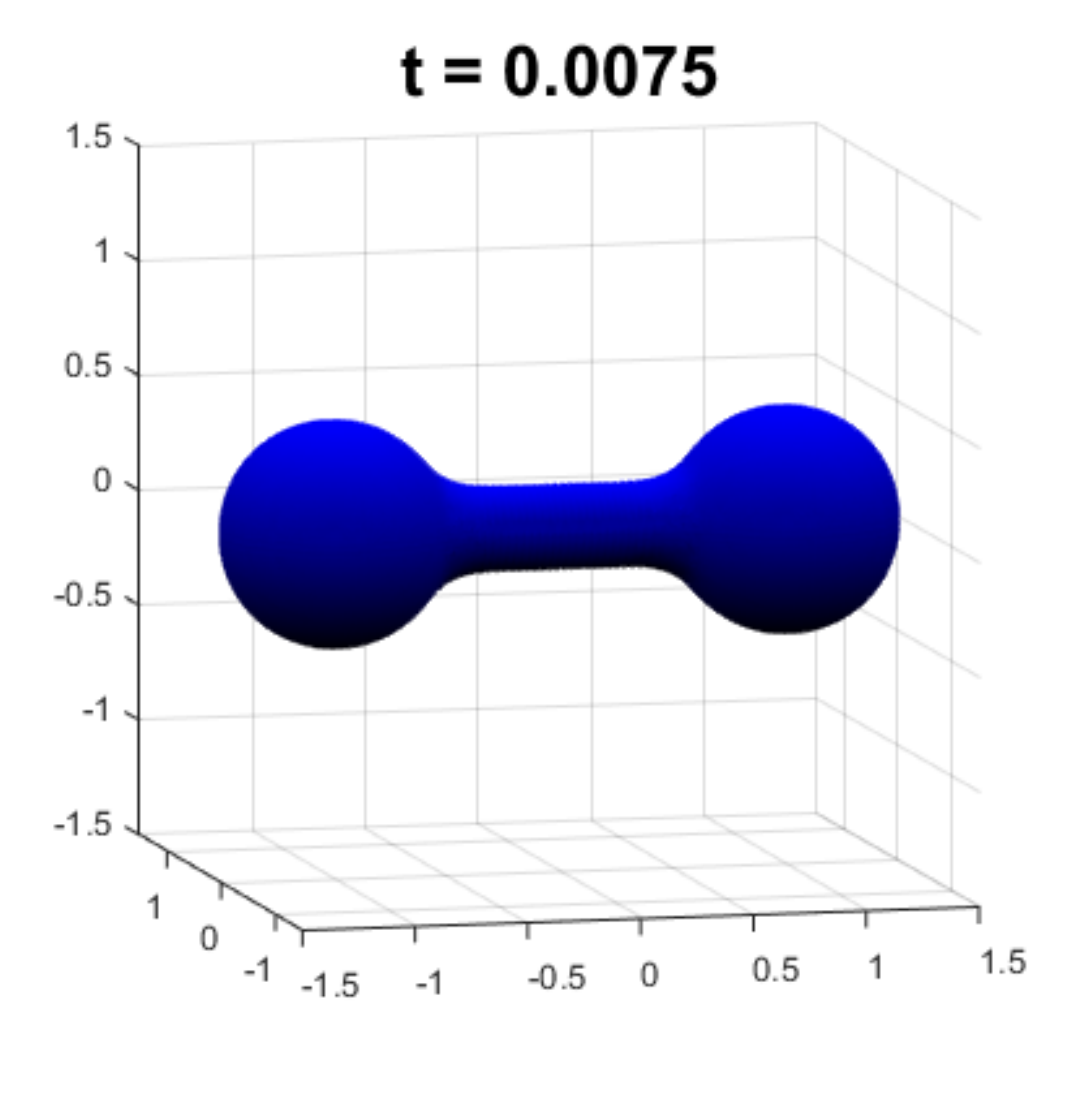}
    \includegraphics[width=0.49\textwidth]{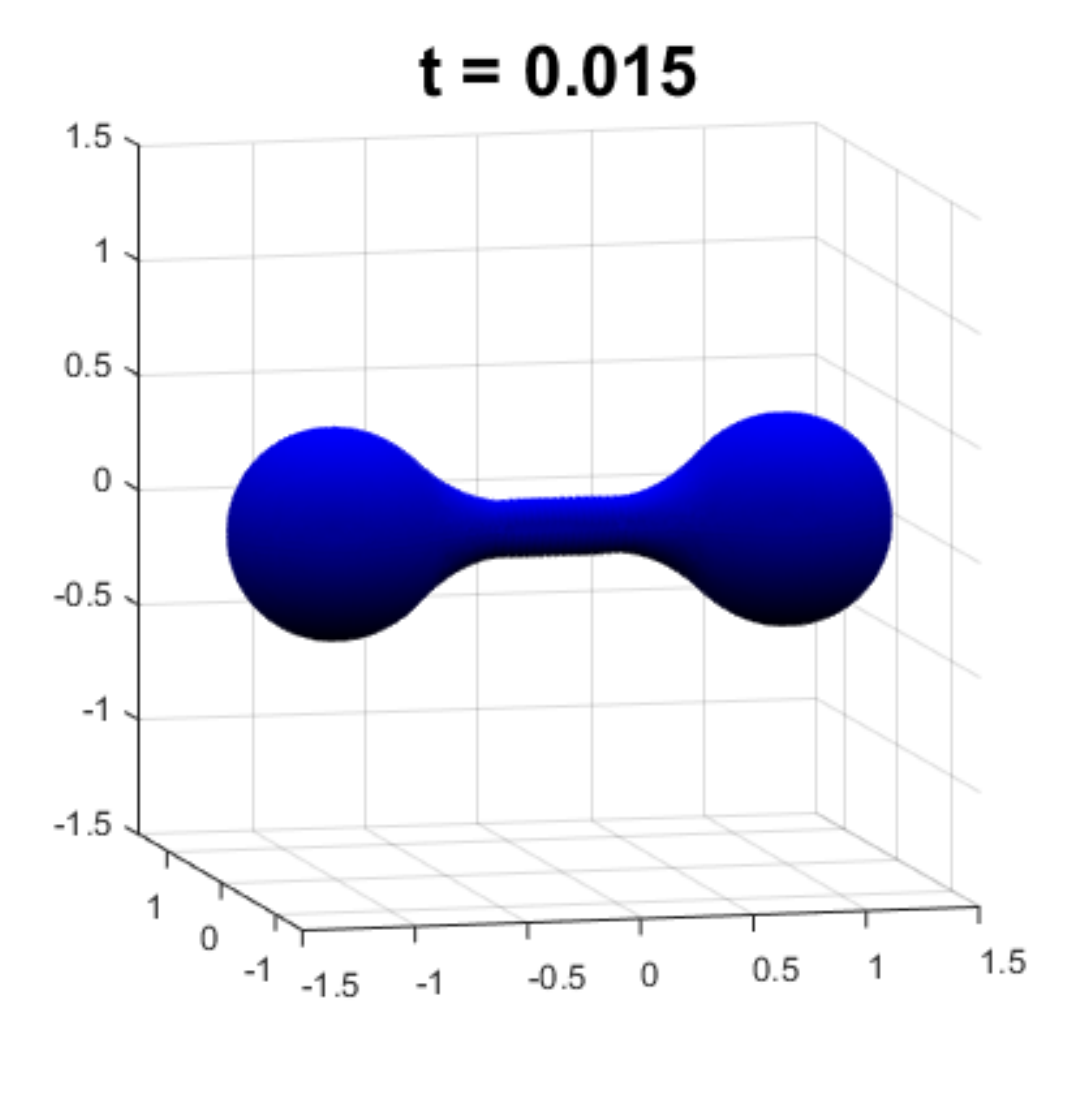}
    \includegraphics[width=0.49\textwidth]{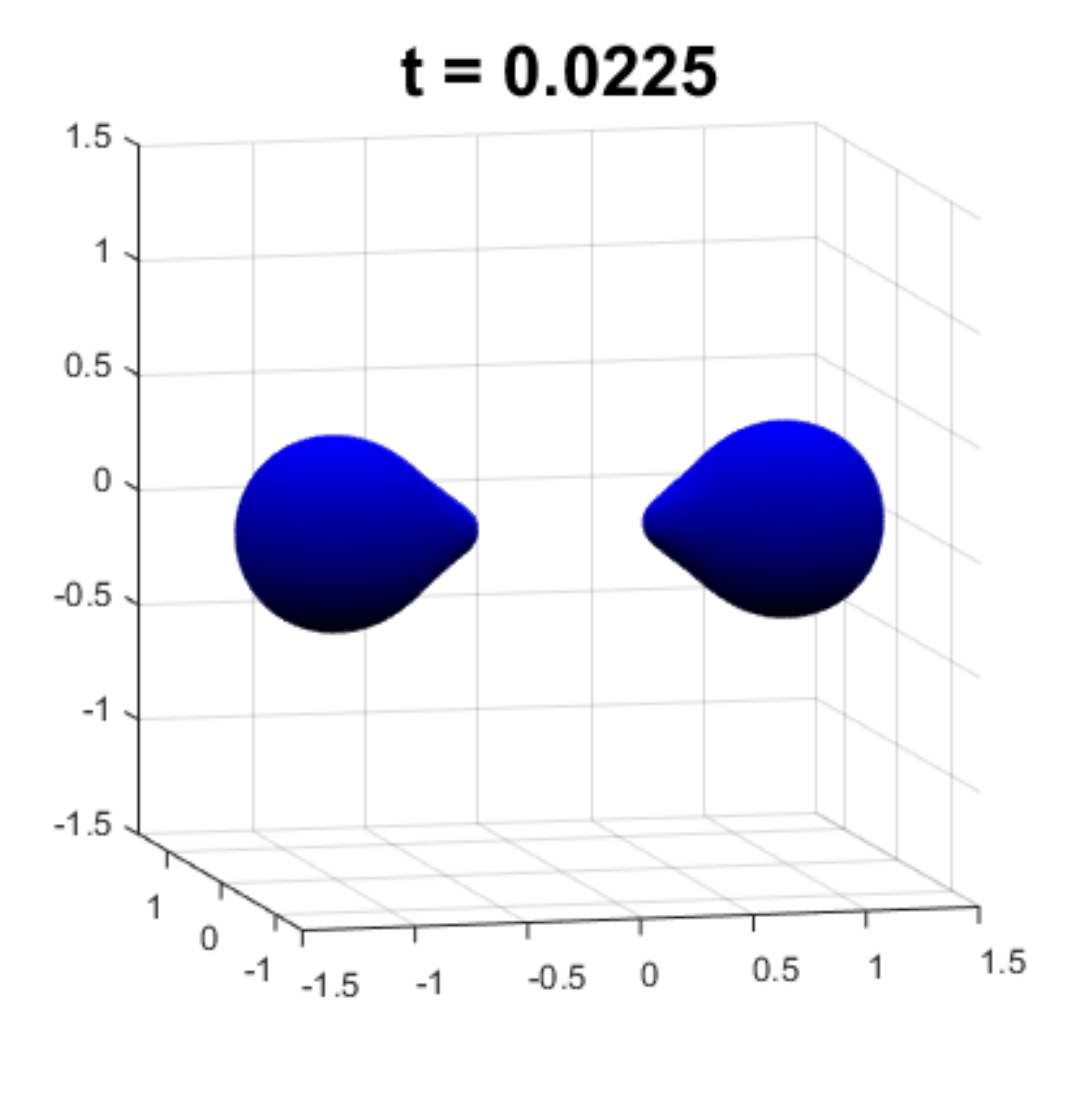}
    \includegraphics[width=0.49\textwidth]{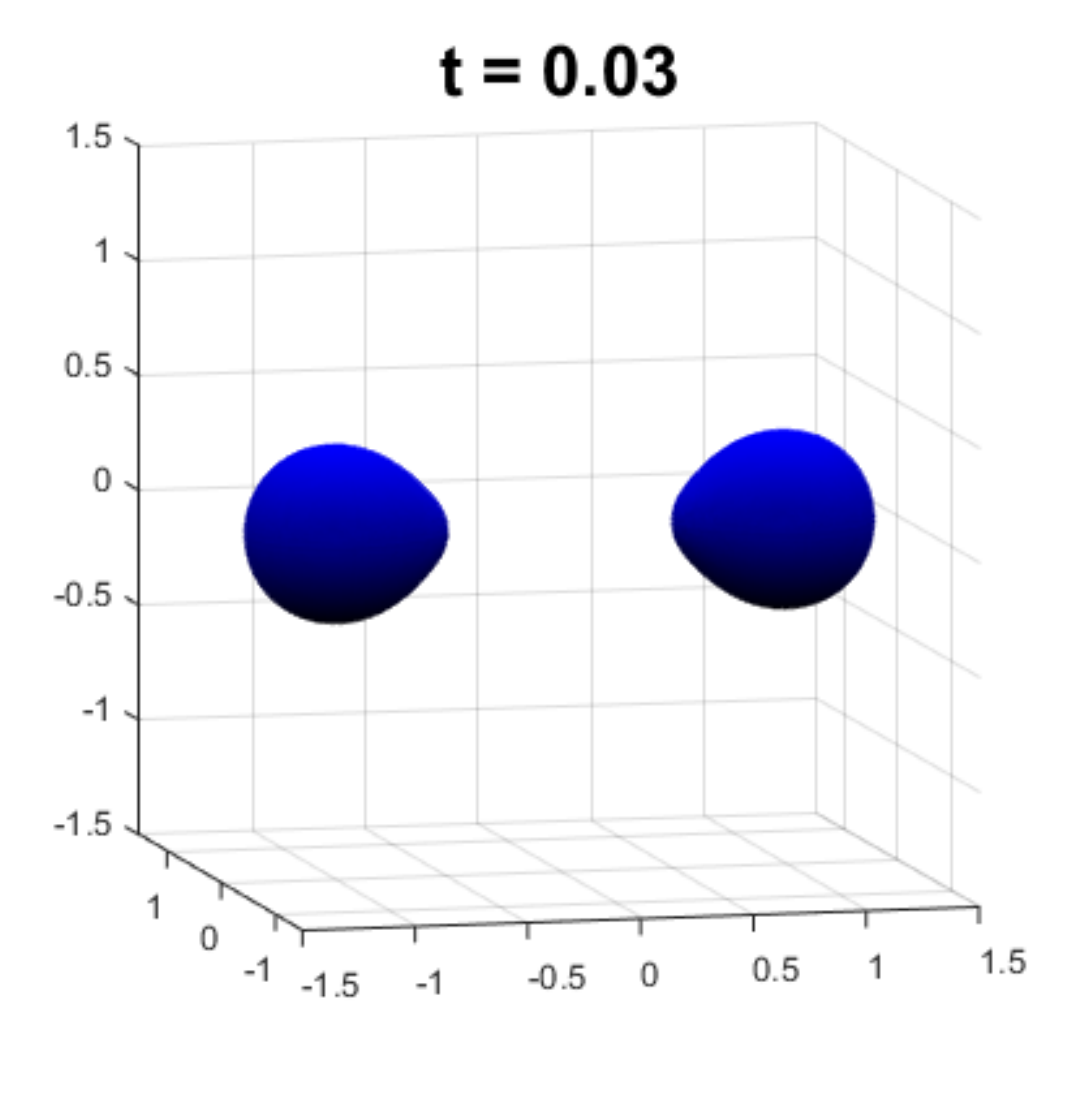}
    \caption{The evolution of a dumbbell under mean curvature motion at selected times $t$.}
    \label{DumbbellSplit}
\end{figure}

Start from an initial dumbbell shape.   Taking discretization step-sizes $\Delta x=0.0333$ and $\Delta t = 0.4\Delta x^2$, and a final time $t_{final}=0.03$, we obtain the results displayed in Figure~\ref{DumbbellSplit}. The modified GBPM captures the split that arises in this well-known example. An examination of the results indicates that the topological change occurs at approximately $t=0.021$.

In this example, control of topology is carried out as described in Section~\ref{GBPMSubsection}. Specifically, if the selected $m$ footpoints have inconsistent Lagrangian information (i.e., $\mathbf{n_0}\cdot\mathbf{n}<\cos\Big(\frac{3\pi}{4}\Big)$), we deactivate the reference grid point and its footpoint.

\section{A coupled method for PDEs on moving surfaces}\label{CoupledMethodSection}
In this section, we propose a coupled method to solve PDEs on surfaces moving by curvature-dependent flows and perform convergence tests. In all examples, we discretize using the forward Euler scheme in time and second-order centered differences in space.

\subsection{A coupled method}
We now introduce our coupled method. To begin, we apply the initialization step of the CPM using a tube radius $\gamma$ (which we shall define later).
This gives a closest point representation over the computational tube.
Nodal values are initialized using a closest point extension of initial surface values: $u(\mathbf{x})=u(cp(\mathbf{x}),0)$.

After initialization, the system is evolved in time. Each time step of size $\Delta t$ consists of two steps:

\begin{enumerate}
  \item \textbf{CPM evolution}: Steps 1-2 of the CPM are applied using a time step-size $\Delta t$
and a tube radius $\gamma$.  (See Section~\ref{CPMSection} for details on the CPM.)
 \item \textbf{Modified GBPM evolution}: We now think of the closest point representation on the computational tube as a collection of footpoints defined on a set of active grid points. Steps 1-3 of the modified GBPM are applied using a time step-size $\Delta t$
and a tube radius $\gamma$.  This gives a closest point representation of the time-evolved surface over a tube of radius $\gamma$.
(See Section~\ref{RemarksSection} for details on the modified GBPM.)
\end{enumerate}

We now discuss the tube radius $\gamma$.
Let $\Gamma_n$ be the computed surface after $n$ time steps,
and $S^\gamma_n$ be the corresponding discrete computational tube
(i.e., the set of grid points that are within a distance $\gamma$ of $\Gamma_n$).
The approximation of $u$ (from the CPM evolution)
will be defined on the set $S^\gamma_{n-1}$, while
the closest point representation (from the modified GBPM evolution)
will be defined on the set $S^\gamma_n$.  The CPM step uses both these
quantities over a tube of radius $\gamma_{CPM}$ (see Section~\ref{CPMSection}).
This suggests setting the tube radius
\begin{equation} \label{eq:safegamma}
\gamma = \gamma_{CPM}+\Delta t \cdot v_n^{max}
\end{equation}
where $v_n^{max}$ is a bound on the normal speed of a footpoint at time step $n$.
While (\ref{eq:safegamma}) does not guarantee a complete tube of nodal values in the CPM step,
it was not observed to fail in our numerical tests.

In practice, condition~(\ref{eq:safegamma}) is too conservative and we select $\gamma$ adaptively.
Specifically, we set $\gamma=\gamma_{CPM}$ and check for violations of the tube radius.
Such events are rare, and when they arise we simply apply a larger computational tube (e.g., we compute using (\ref{eq:safegamma})).

\subsection{Numerical experiments}
In this section, the coupled method is applied to a number of examples in two and three dimensions.
Unless stated otherwise, we solve Equation~(\ref{AdvDifPDE2}) with a flux $\mathbf{q}=\nabla_\Gamma u$.
This corresponds to diffusion with diffusivity parameter $\mathcal{D}=1$.
We note that the use of the CPM extension step in our coupled method leads to a vanishing normal derivative, $\partial u/\partial n=0$,  and a corresponding simplification in Equation~(\ref{AdvDifPDE2}).

\subsubsection{Diffusion on an expanding circle}\label{Examples2D}
In our first experiment, we approximate diffusion on an expanding circle (cf. \cite{elliott2011numerical}). The surface starts from a unit circle and evolves outwards with a constant velocity given by $\mathbf{v}=5\mathbf{n}$.

\begin{figure}
    \centering
    \includegraphics[width=0.24\textwidth]{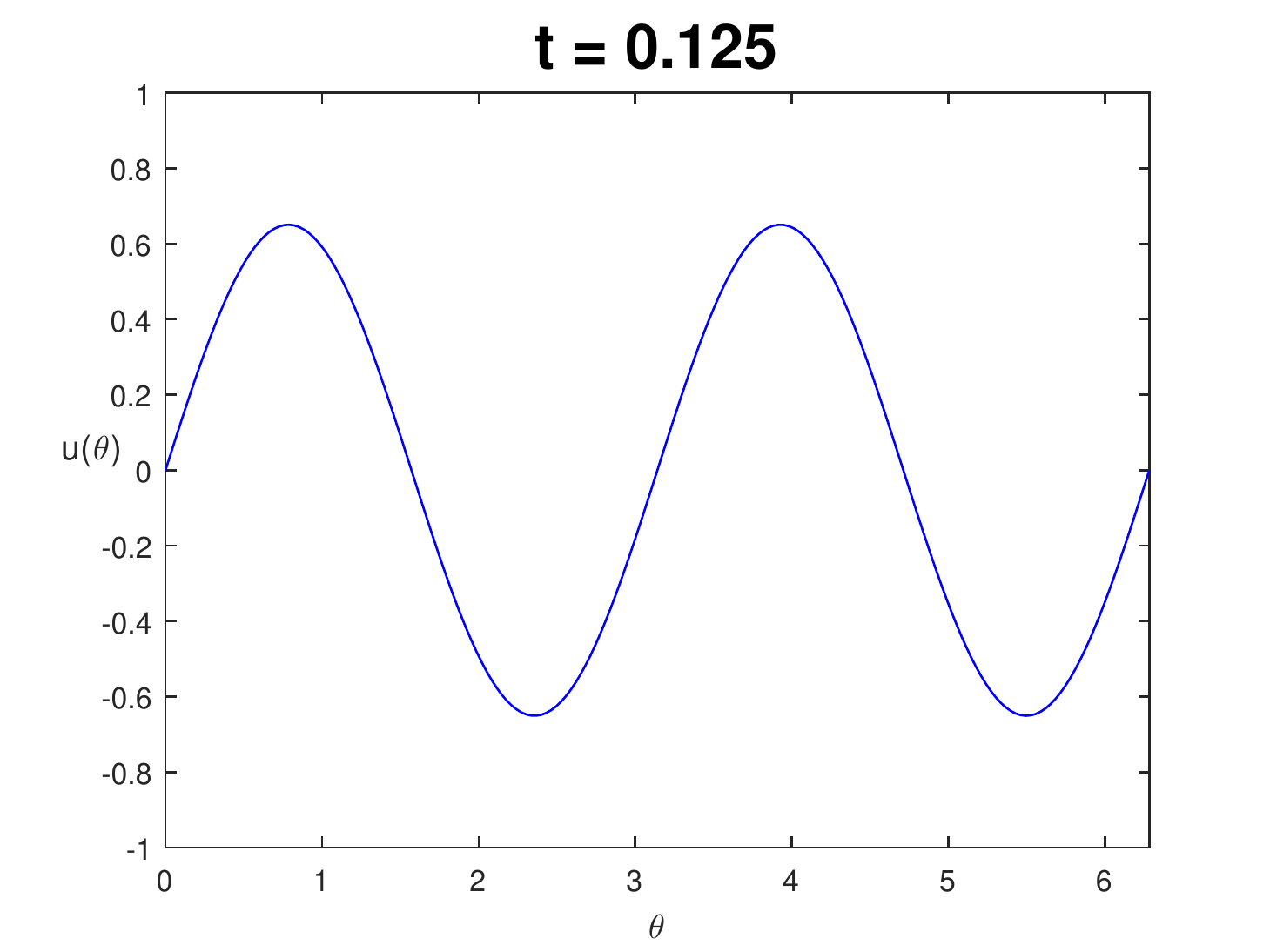}
    \includegraphics[width=0.24\textwidth]{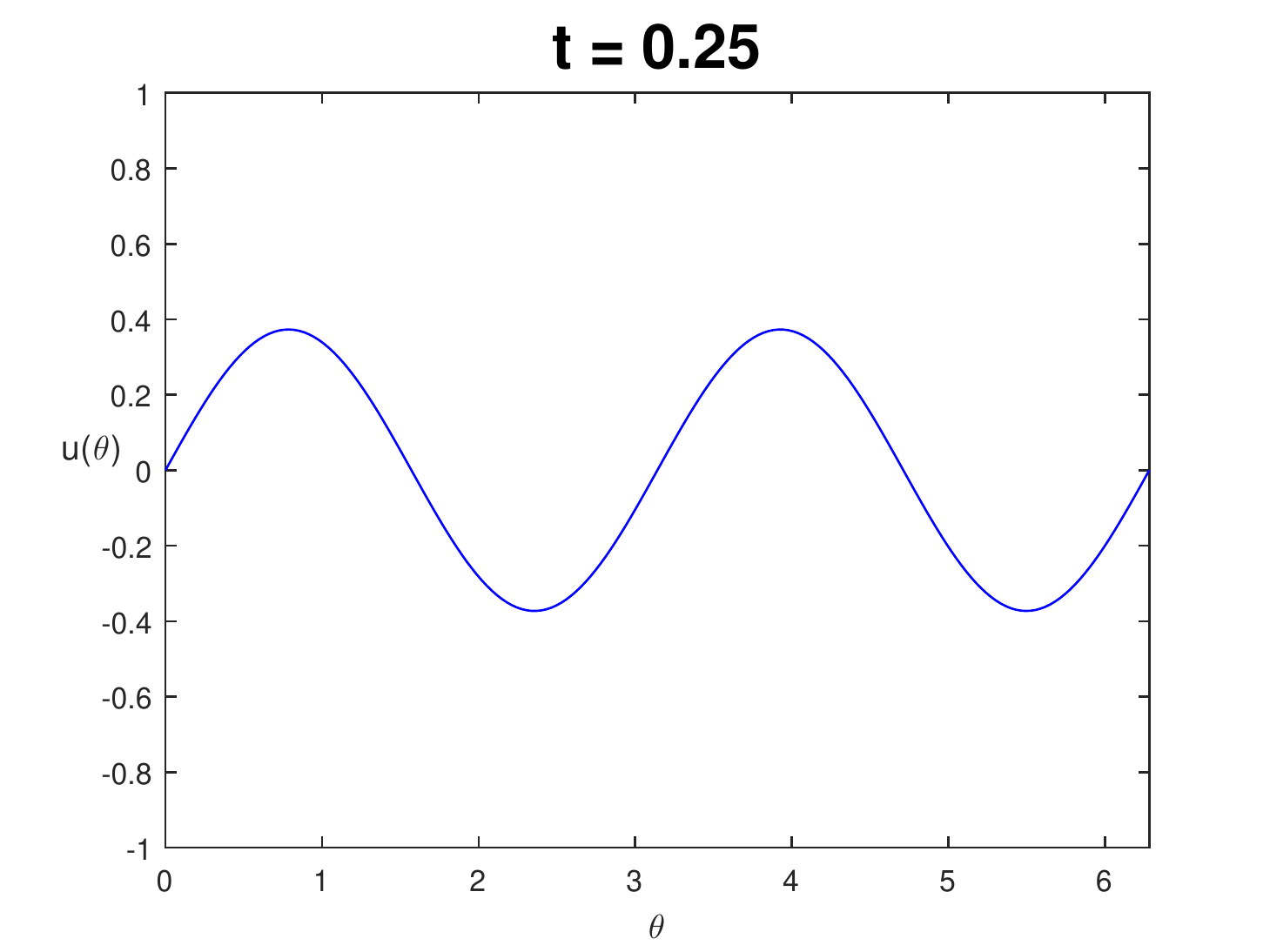}
    \includegraphics[width=0.24\textwidth]{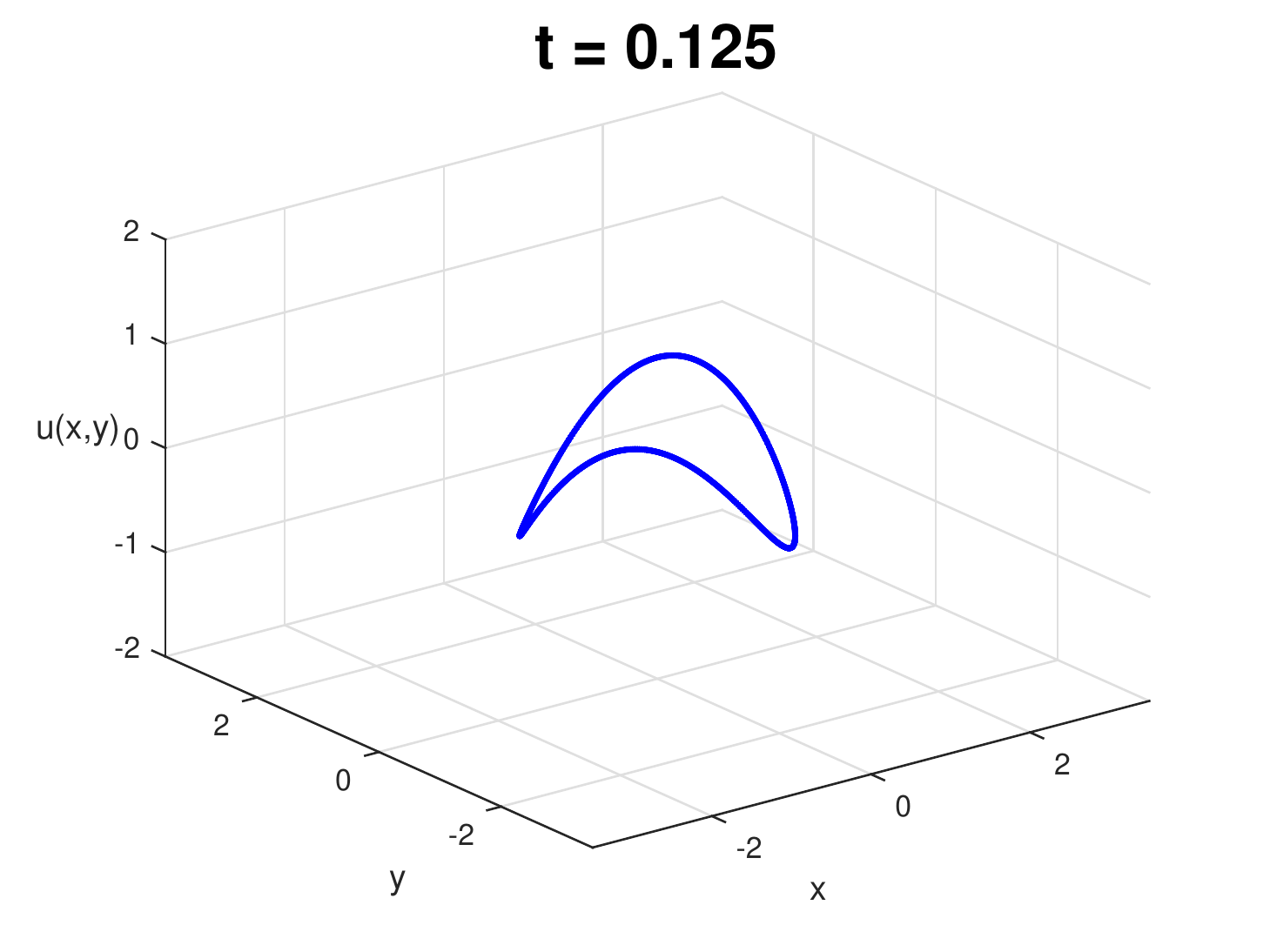}
    \includegraphics[width=0.24\textwidth]{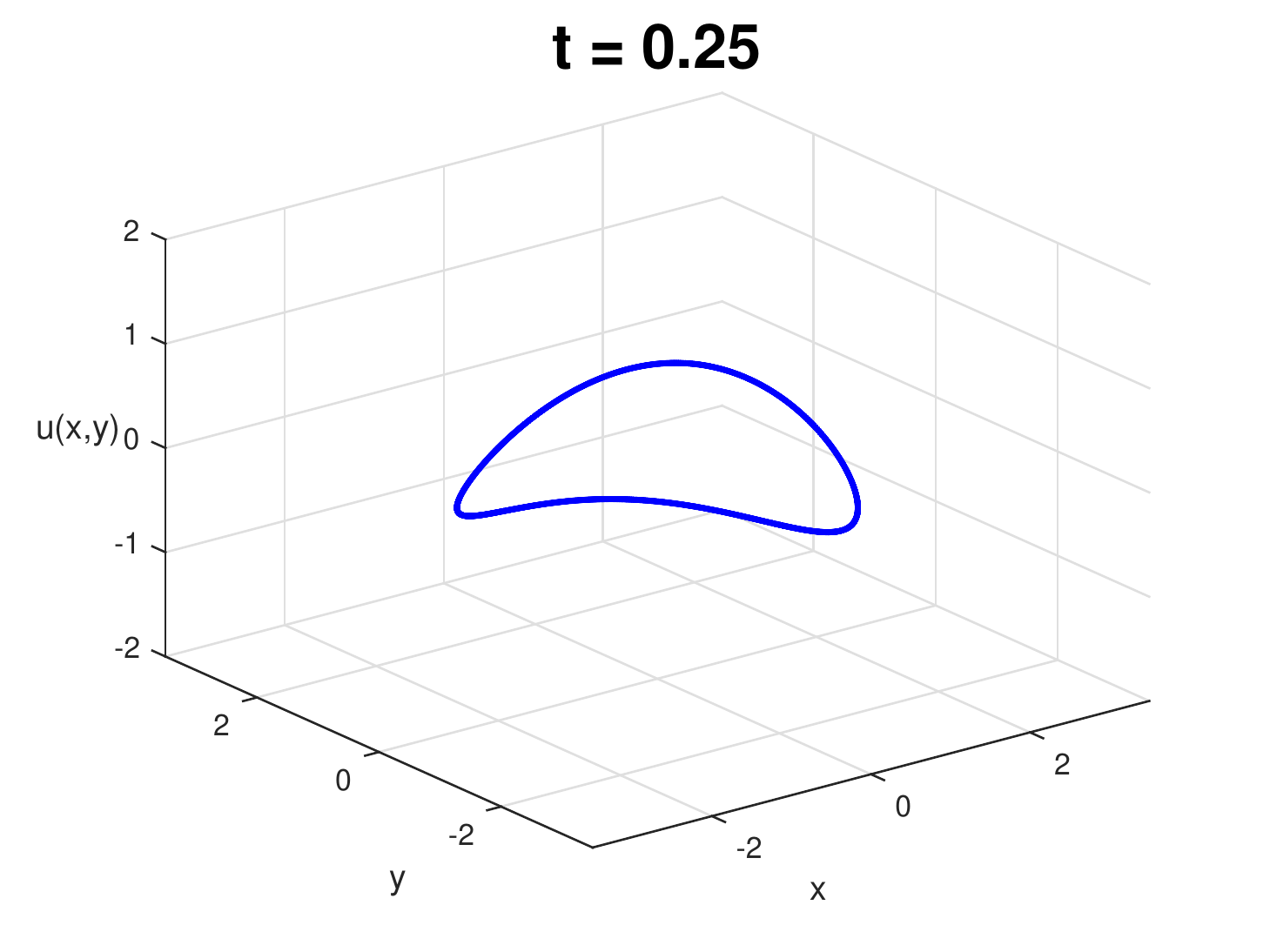}
    \includegraphics[width=0.24\textwidth]{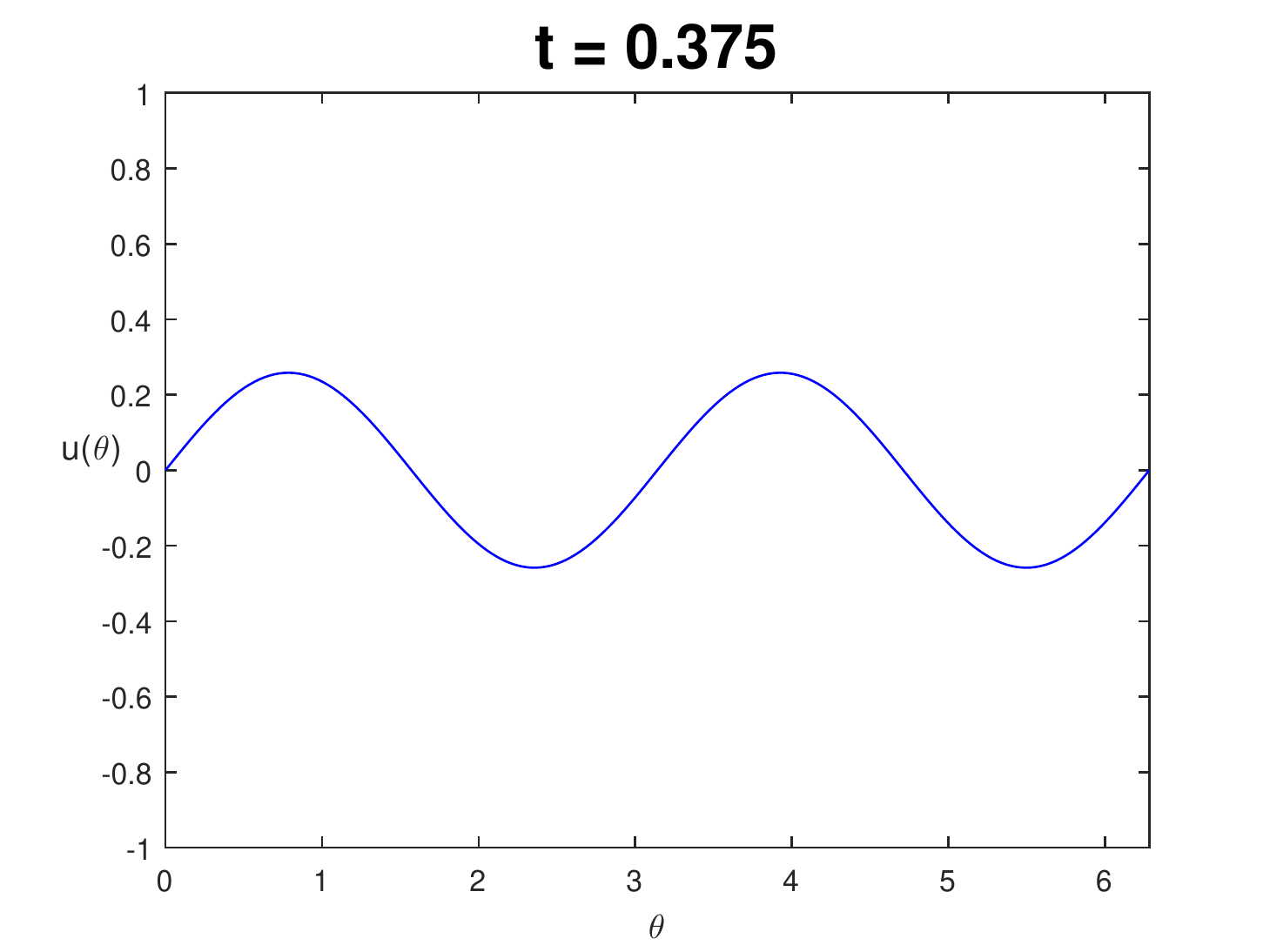}
    \includegraphics[width=0.24\textwidth]{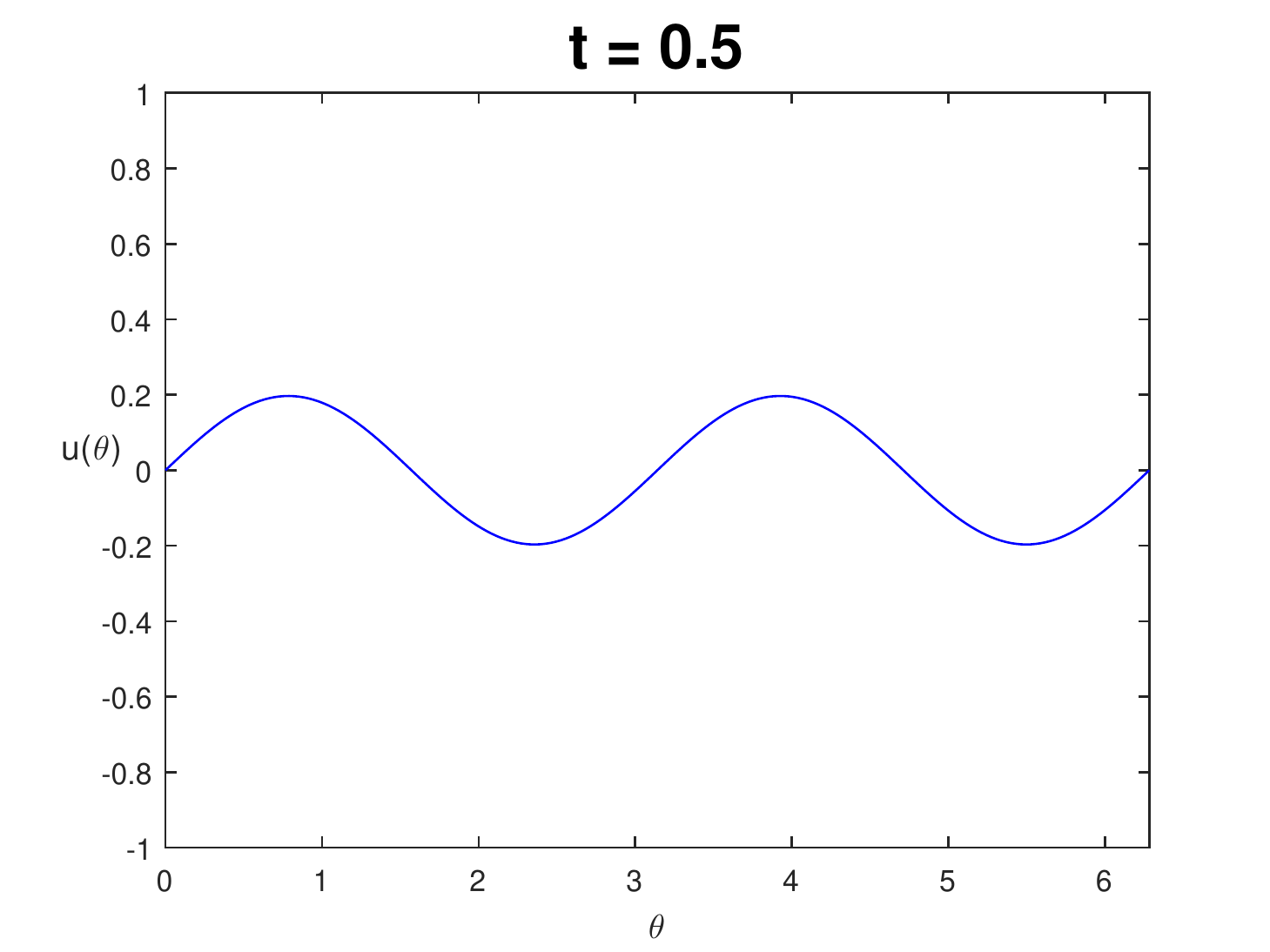}
    \includegraphics[width=0.24\textwidth]{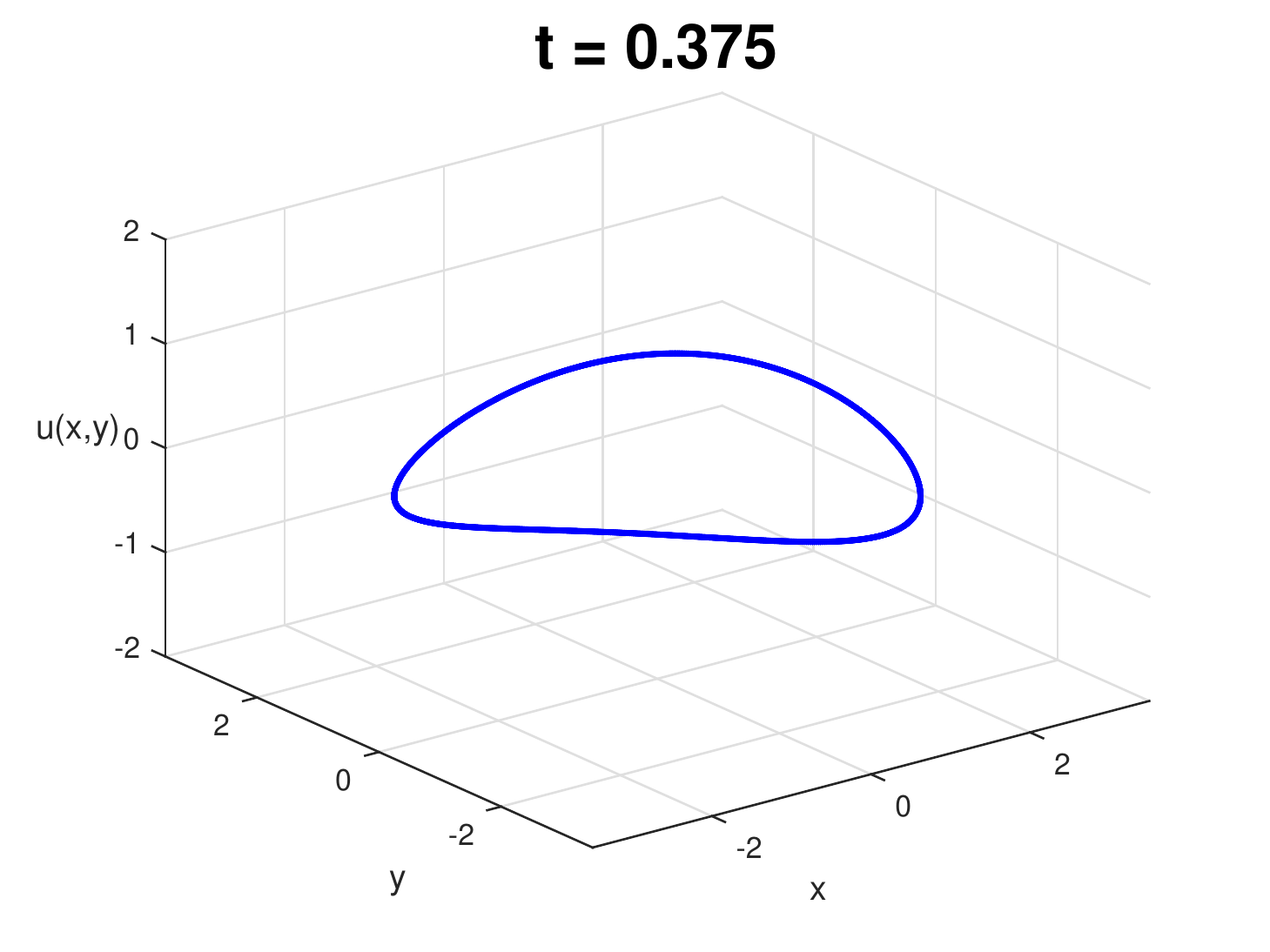}
    \includegraphics[width=0.24\textwidth]{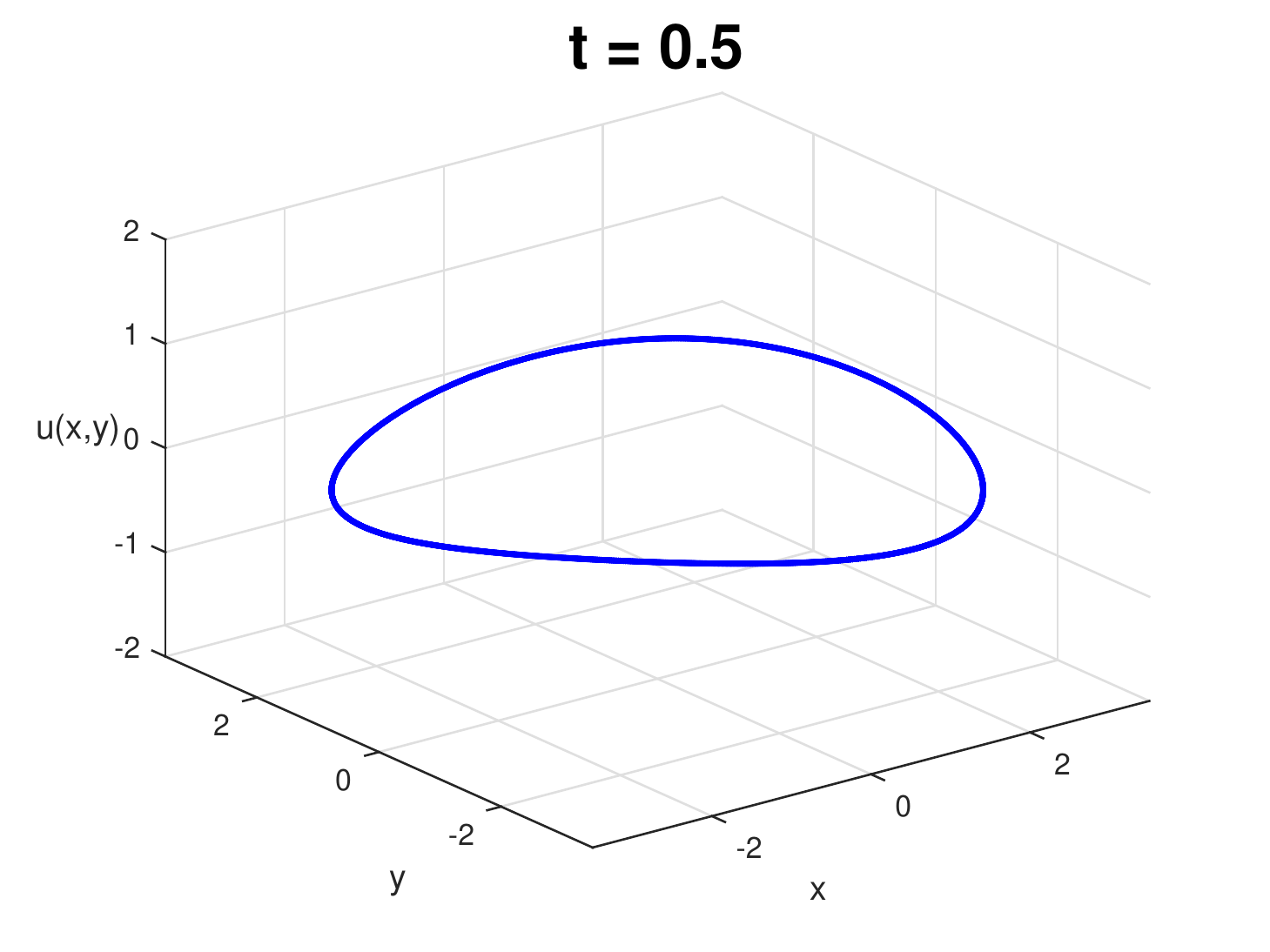}
    \caption{The graph of the solution of the diffusion model on the evolving circle at selected times $t$. Left:  The solution is depicted in relation to the surface parameter $\theta$. Right: The solution is shown in relation to the $x$ and $y$ coordinates.}
    \label{solCircle}
\end{figure}

The exact solution of the homogeneous PDE (\ref{AdvDifPDE2}) on our evolving circle is
$$u(\theta,t)=e^{4/(5r(t))}\frac{\cos\theta\;\sin\theta}{r(t)}$$
where $r(t)=1+5t$ is the radius of the circle at time $t$. The result may be verified by substituting $u$ into the PDE (\ref{AdvDifPDE}) and noting that $$\nabla_{\Gamma}=(-\sin\theta,\cos\theta)\frac{\partial_{\theta}}{|\gamma(\theta,t)|},$$
where $\Gamma=\Gamma(t)$ is the moving surface and $|\gamma(\theta,t)|$ is the length of the curve. Because there are no tangential terms in the velocity, diffusion dominates the flow. The evolution of the surface and the solution are depicted in Figure~\ref{solCircle}.

\begin{figure}
    \centering
    \includegraphics[width=0.49\textwidth]{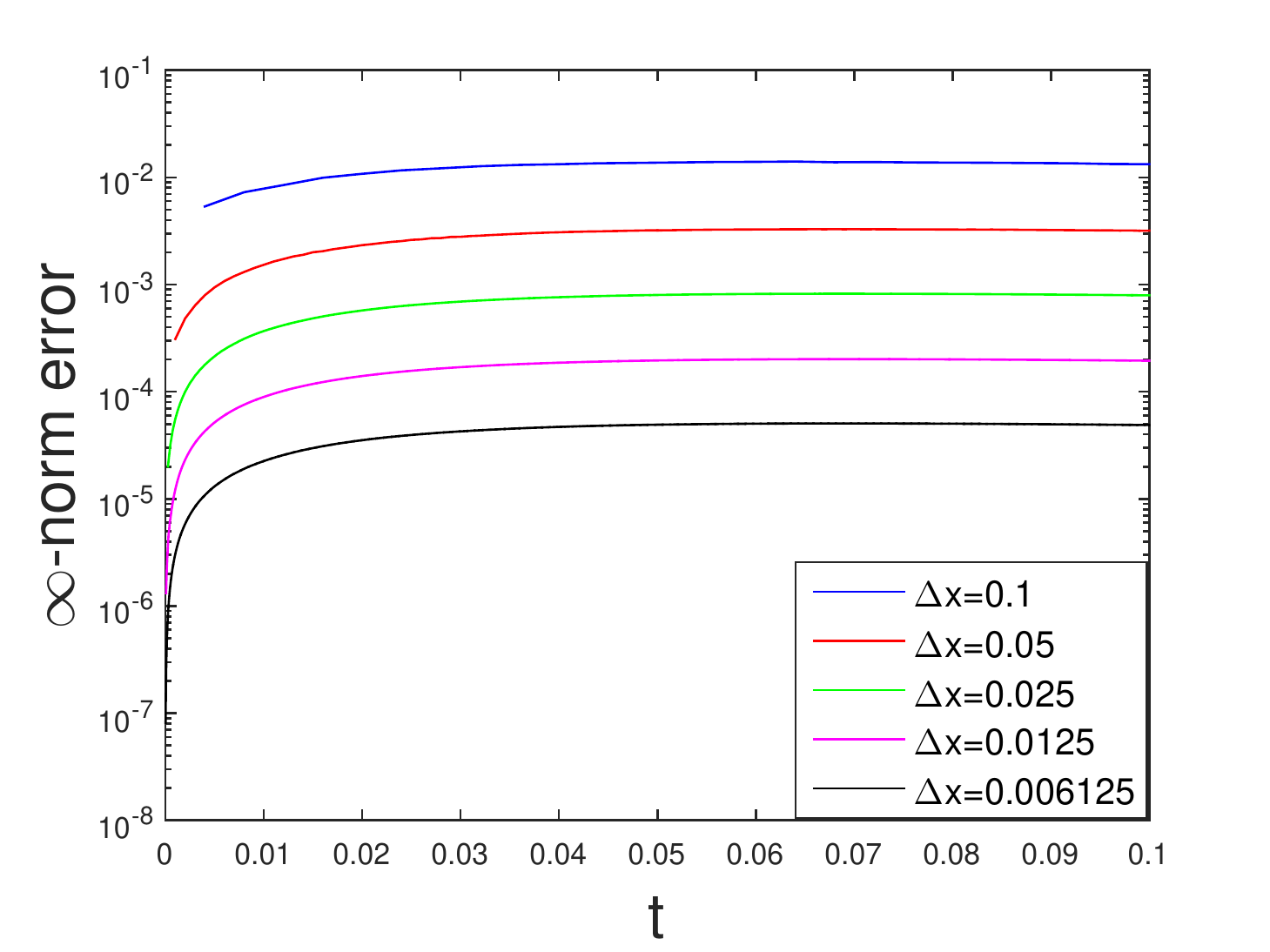}
    \caption{The absolute error of the numerical solution of the diffusion model on the expanding circle over time. Errors are computed using the analytical solution.}
    \label{errorCircle}
\end{figure}

We compute the absolute error as measured in the infinity norm at various times $k\Delta t$, i.e., we compute
$$\max_{i}|u(\mathbf{x}_i,k\Delta t)-U_i^k|$$
where $u(\mathbf{x}_i,k\Delta t)$ denotes the analytical solution at node $\mathbf{x}_i$ and $U_i^k$ denotes the corresponding numerical result. As shown in Figure~\ref{errorCircle}, errors accumulate in time and are reduced by mesh refinement.

To clarify the convergence rates, Table~\ref{errorCircleTable} provides the results and a numerical convergence study for selected times $t$. In all four cases, a clear second-order convergence is observed.

\begin{table}
\centering
  \begin{tabular}{|l|cc|cc|cc|cc|}
    \hline
    $\Delta x$ & $t =$ 0.025 & e.o.c. & $t =$ 0.05 & e.o.c. & $t =$ 0.075 & e.o.c. & $t =$ 0.1 & e.o.c.  \\ \hline
    0.1 & 1.16$\times 10^{-2}$ & - & 1.36$\times 10^{-2}$ & - & 1.37$\times 10^{-2}$ & - & 1.33$\times 10^{-2}$ & - \\
    0.05 & 2.61$\times 10^{-3}$ & 2.16 & 3.20$\times 10^{-3}$ & 2.09 & 3.27$\times 10^{-3}$ & 2.07 & 3.19$\times 10^{-3}$ & 2.06\\
    0.025 & 6.42$\times 10^{-4}$ & 2.03 & 8.00$\times 10^{-4}$ & 2.00 & 8.20$\times 10^{-4}$ & 2.00 & 7.93$\times 10^{-4}$ & 2.01\\
    0.0125 & 1.57$\times 10^{-4}$ & 2.03 & 1.96$\times 10^{-4}$ & 2.03 & 2.01$\times 10^{-4}$ & 2.03 & 1.95$\times 10^{-4}$ & 2.03 \\
    0.006125 & 3.95$\times 10^{-5}$ & 1.99 & 4.92$\times 10^{-5}$ & 1.99 & 5.05$\times 10^{-5}$ & 1.99 & 4.88$\times 10^{-5}$ & 1.99\\
    \hline
  \end{tabular}
  \caption{Absolute errors as measured in the infinity norm and the estimated order of convergence (e.o.c.) at various times $t$.}
  \label{errorCircleTable}
\end{table}

\subsubsection{Advection-diffusion on an oscillating sphere}
Moving to three dimensions, we apply the velocity
$$\mathbf{v} = \frac{a'(t)}{2a(t)}(x_1,0,0),$$
where
$$a(t)=1+\sin(2t)$$
to the unit sphere centered at the origin. Following \cite{elliott2011numerical}, the exact solution of the surface at all times $t\geq0$ is
$$\mathbf{x}(\theta,\phi,t) = (\sqrt{a(t)}\cos(\theta)\cos(\phi),\sin(\theta)\cos(\phi),\sin(\phi))$$
where $\theta$ and $\phi$ correspond to the azimuth and the elevation of the sphere respectively. Notice that this velocity causes the sphere to oscillate along the $x$-axis with period $\pi/2$.

To form our surface PDE, we insert
$$u(\mathbf{x},t) = e^{-6t}x_1x_2$$
into the non-homogeneous PDE (\ref{AdvDifPDE2}) and compute the right-hand side $f$. This gives
$$f = u(\mathbf{x},t) \times \Big(-6+\frac{a'(t)}{a(t)}\Big(1-\frac{x_1^2}{2N}\Big)+\frac{1+5a(t)+2a^2(t)}{N}-\frac{1+a(t)}{N^2}(x_1^2+a^3(t)(x_2^2+x_3^2))\Big),$$
with
$$N = x_1^2+a^2(t)(x_2^2+x_3^2).$$

\begin{figure}
    \centering
    \includegraphics[width=0.24\textwidth]{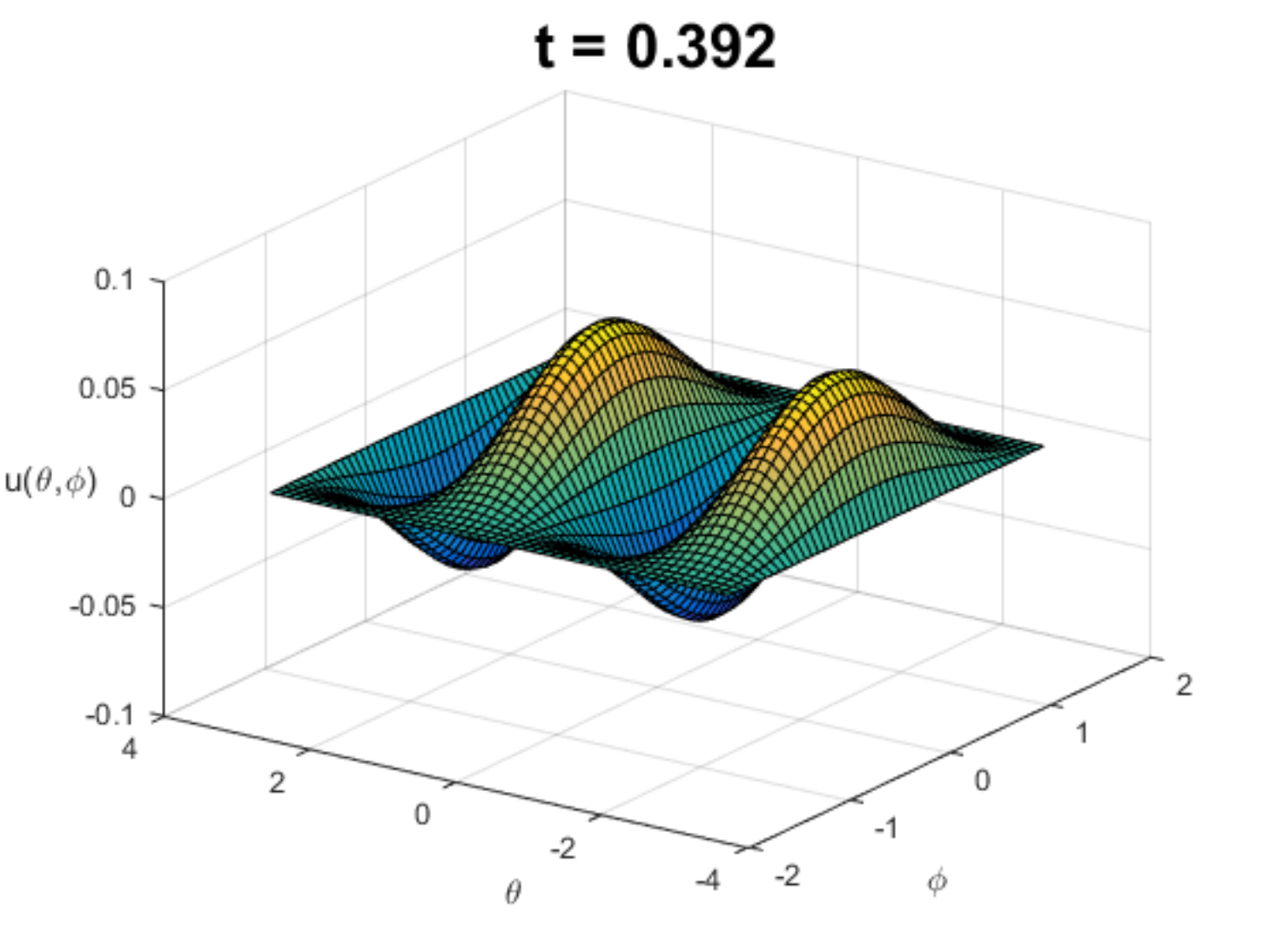}
    \includegraphics[width=0.24\textwidth]{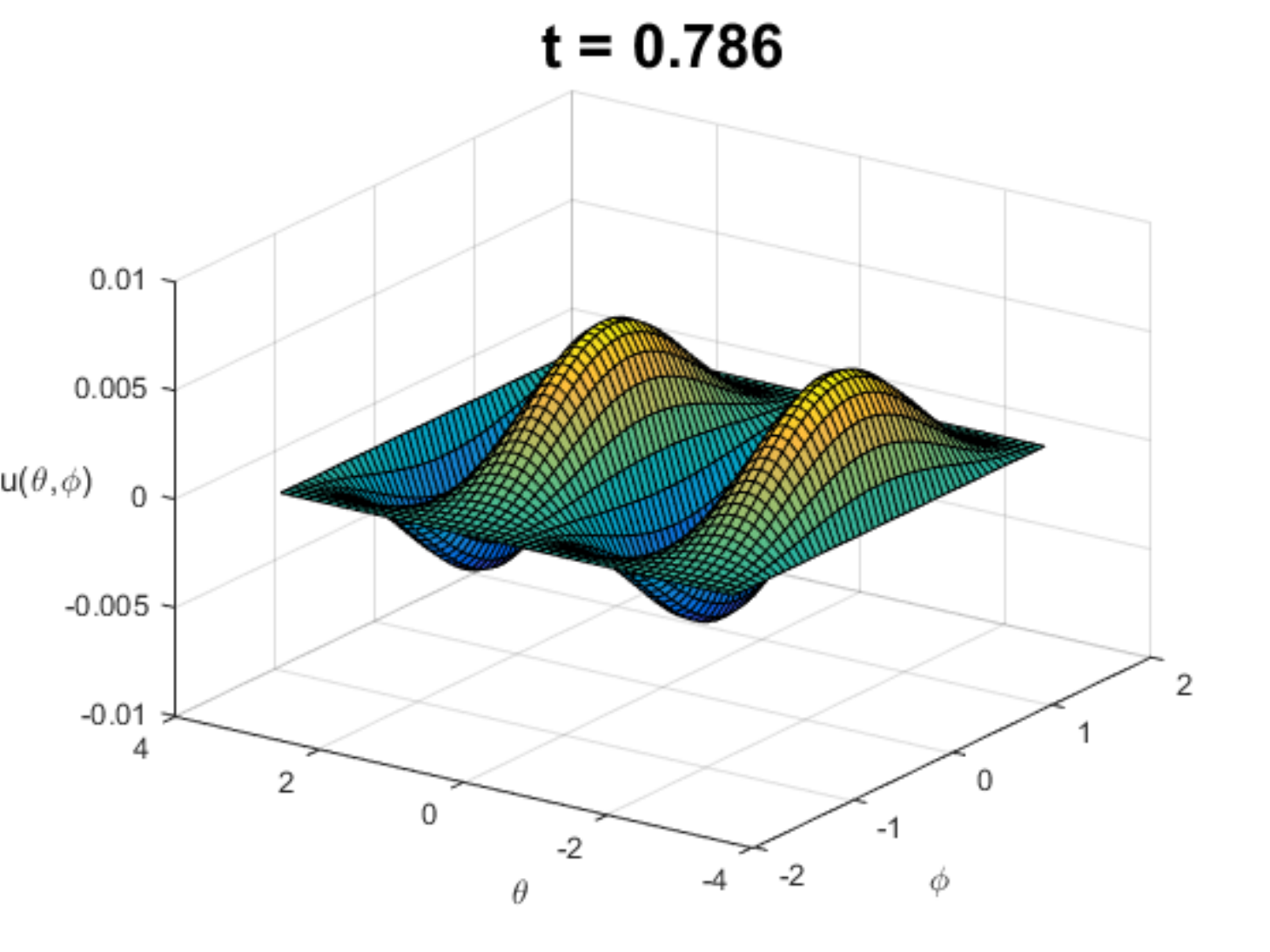}
    \includegraphics[width=0.24\textwidth]{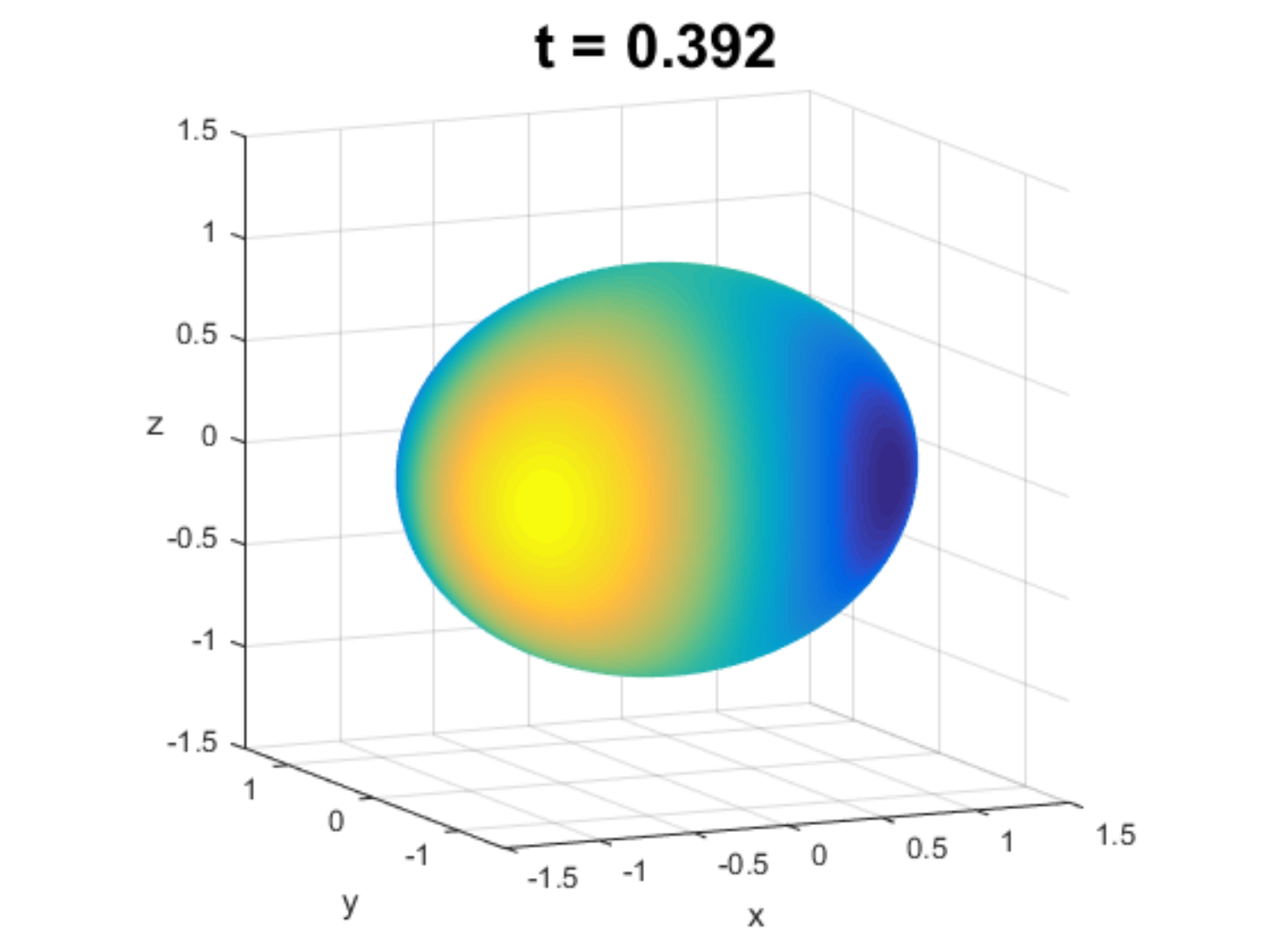}
    \includegraphics[width=0.24\textwidth]{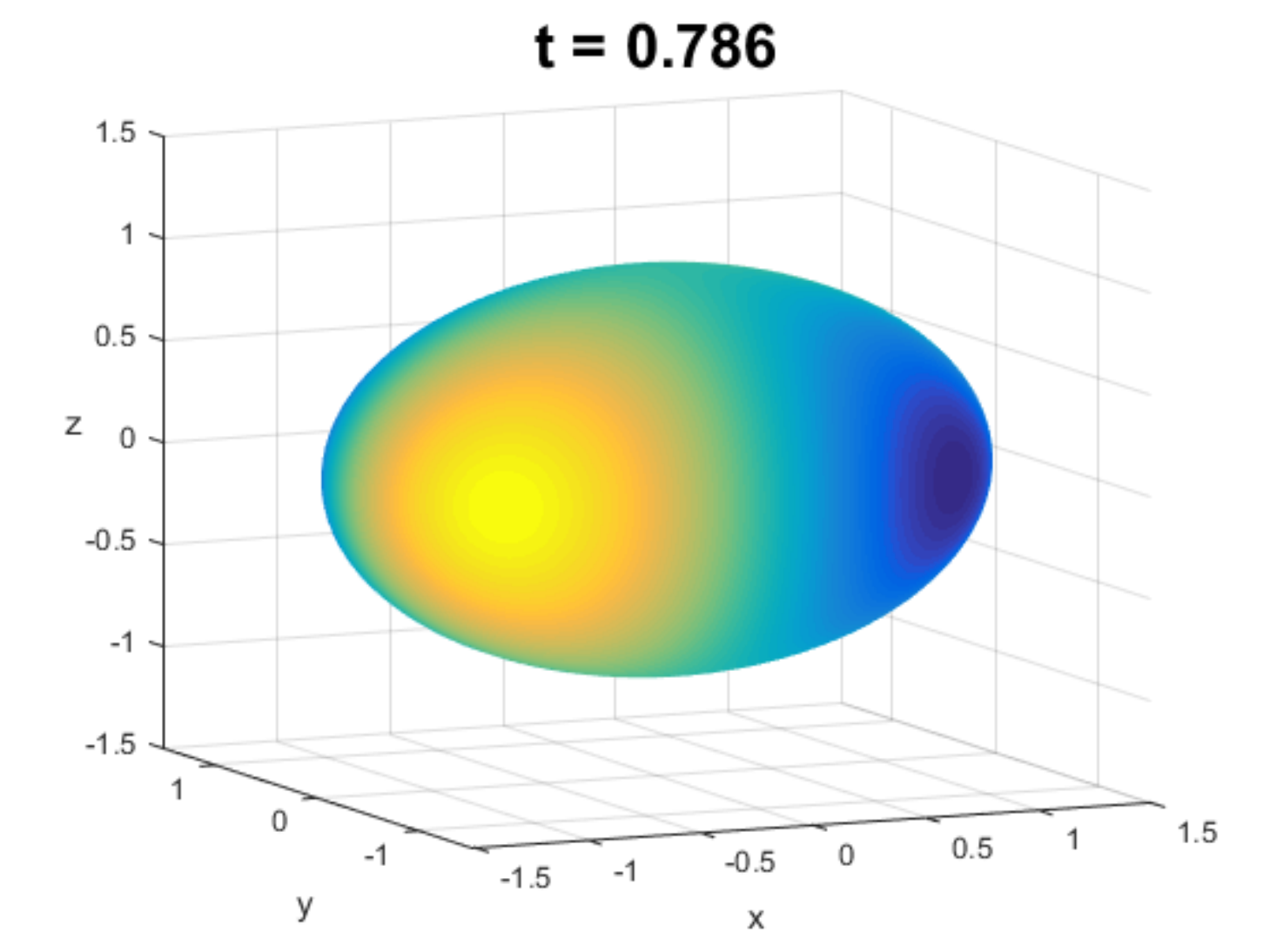}
    \includegraphics[width=0.24\textwidth]{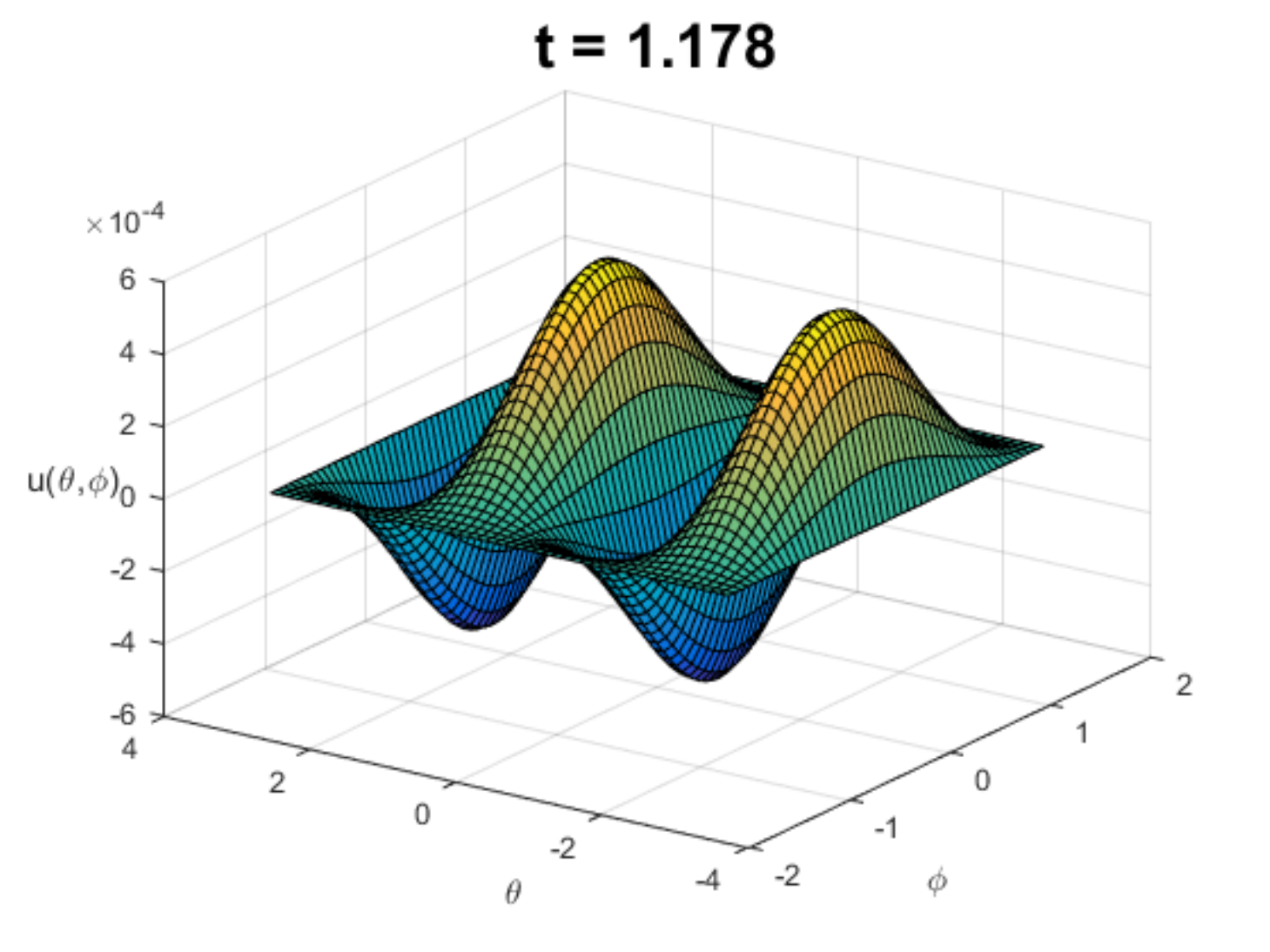}
    \includegraphics[width=0.24\textwidth]{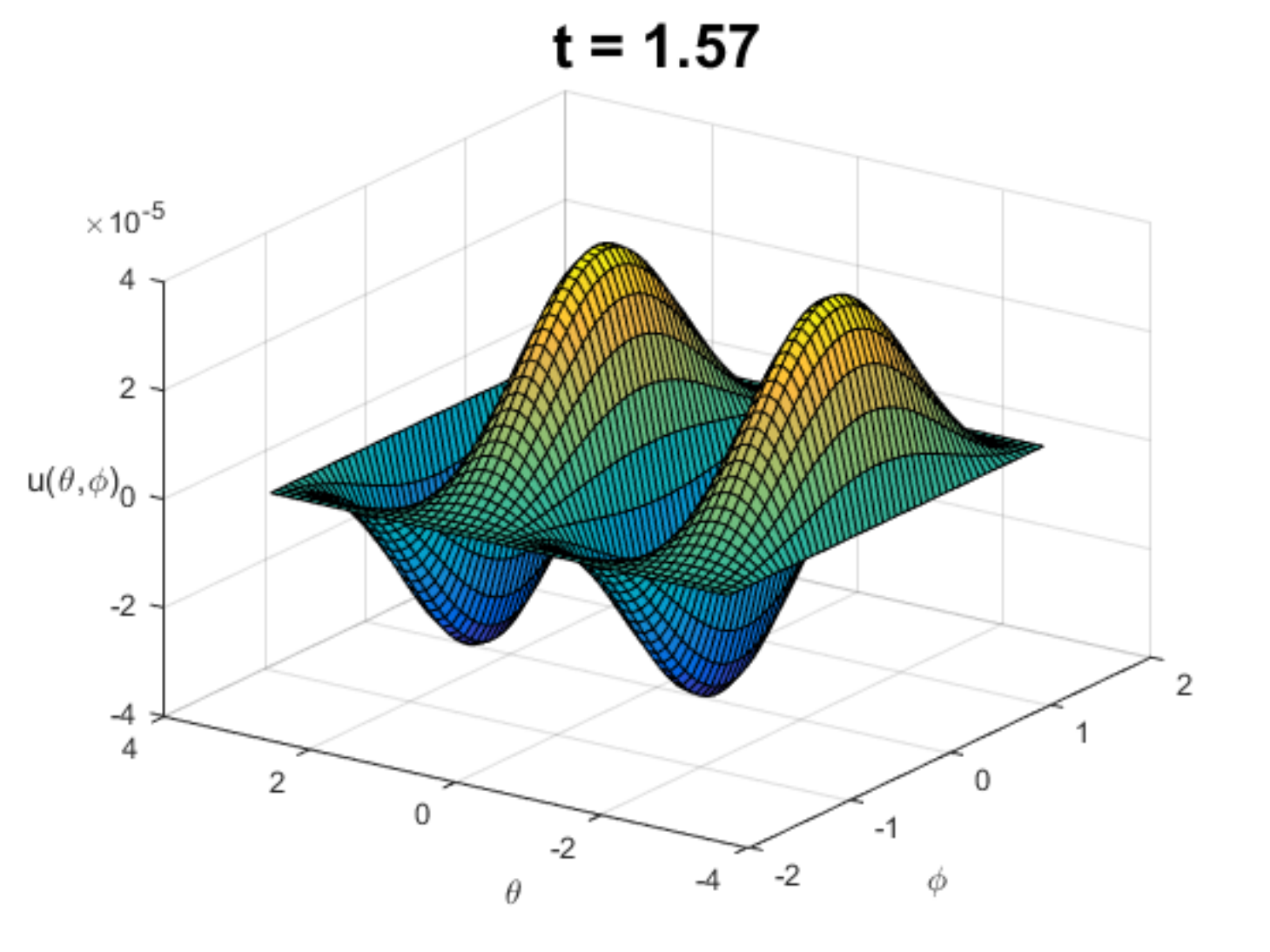}
    \includegraphics[width=0.24\textwidth]{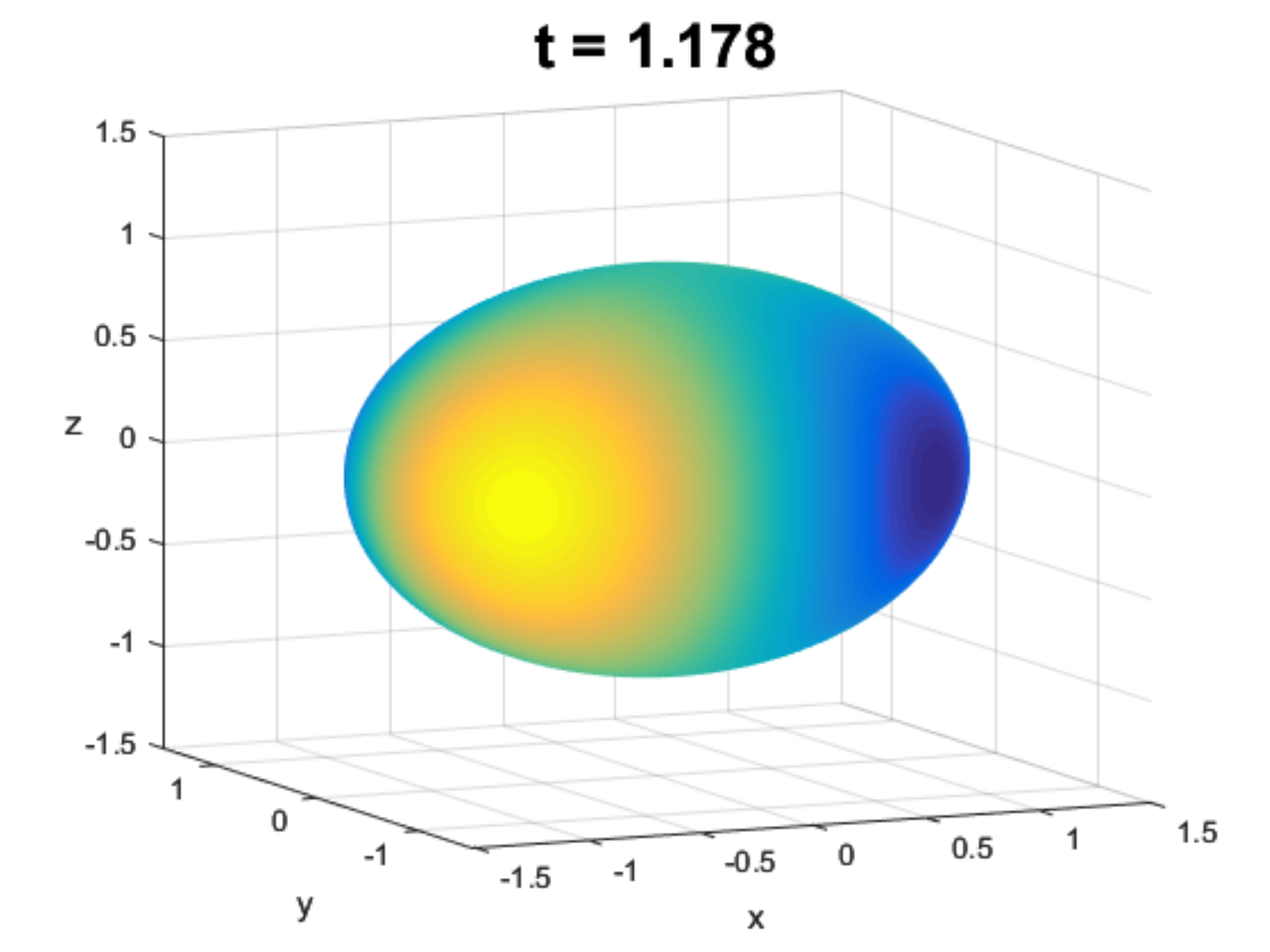}
    \includegraphics[width=0.24\textwidth]{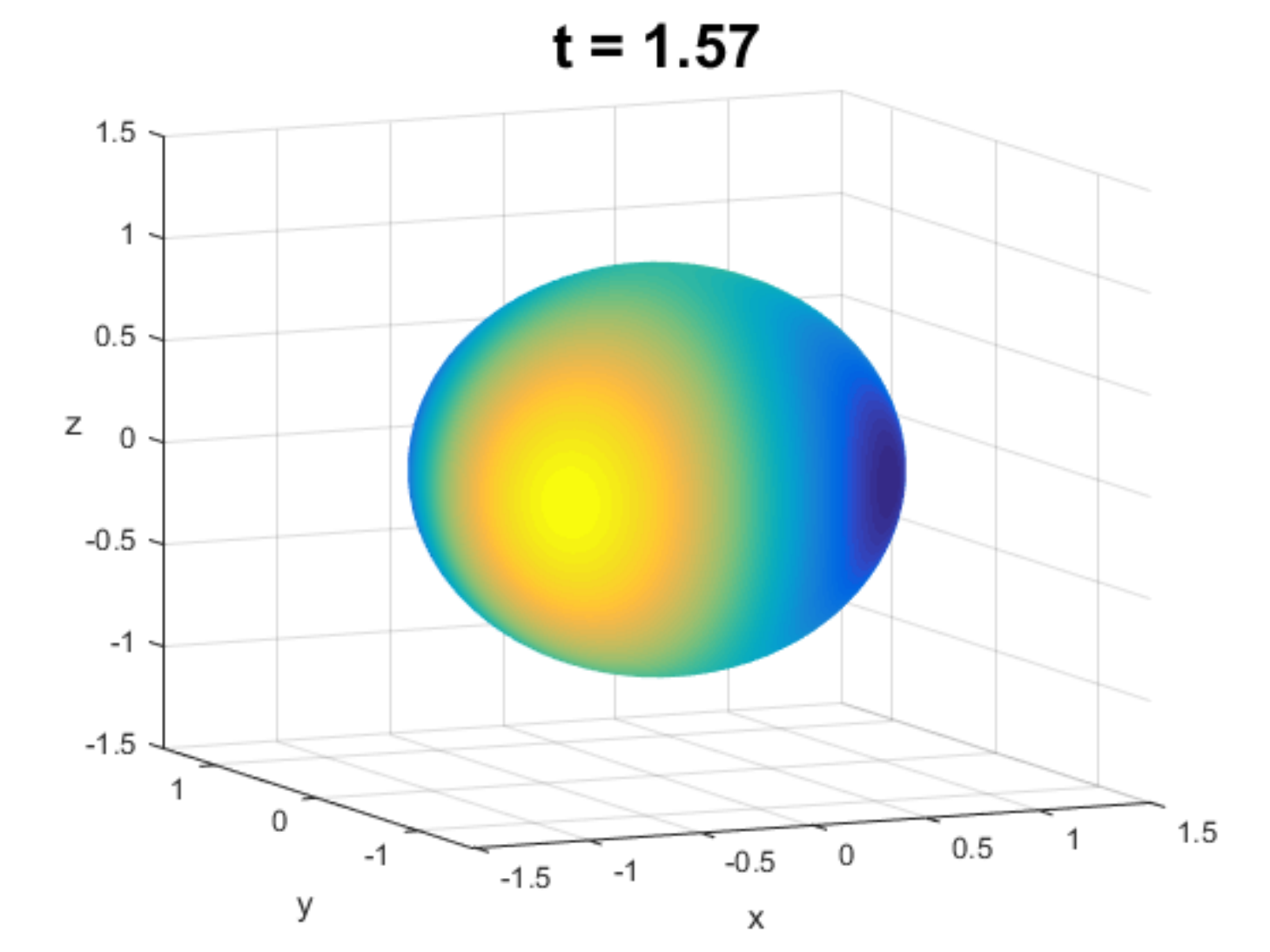}
    \caption{The graph of the numerical solution of an advection-diffusion model on an oscillating ellipsoid at various times $t$. Left: The solution is plotted as a graph of the surface parameters $\theta$ and $\phi$. Right: The solution is visualized on the evolving surface. Yellow corresponds to large solution values and dark blue corresponds to small values.}
    \label{solEllipsoid}
\end{figure}

Using a mesh spacing $\Delta x=0.1$ and a time step-size $\Delta t=0.1\Delta x^2$, we obtain the numerical solution displayed in Figure~\ref{solEllipsoid} at the final time $t_{final}=\pi/2$. The relative error at the final time is less than $5\%$ as measured in the infinity norm.

\begin{figure}
    \centering
    \includegraphics[width=0.49\textwidth]{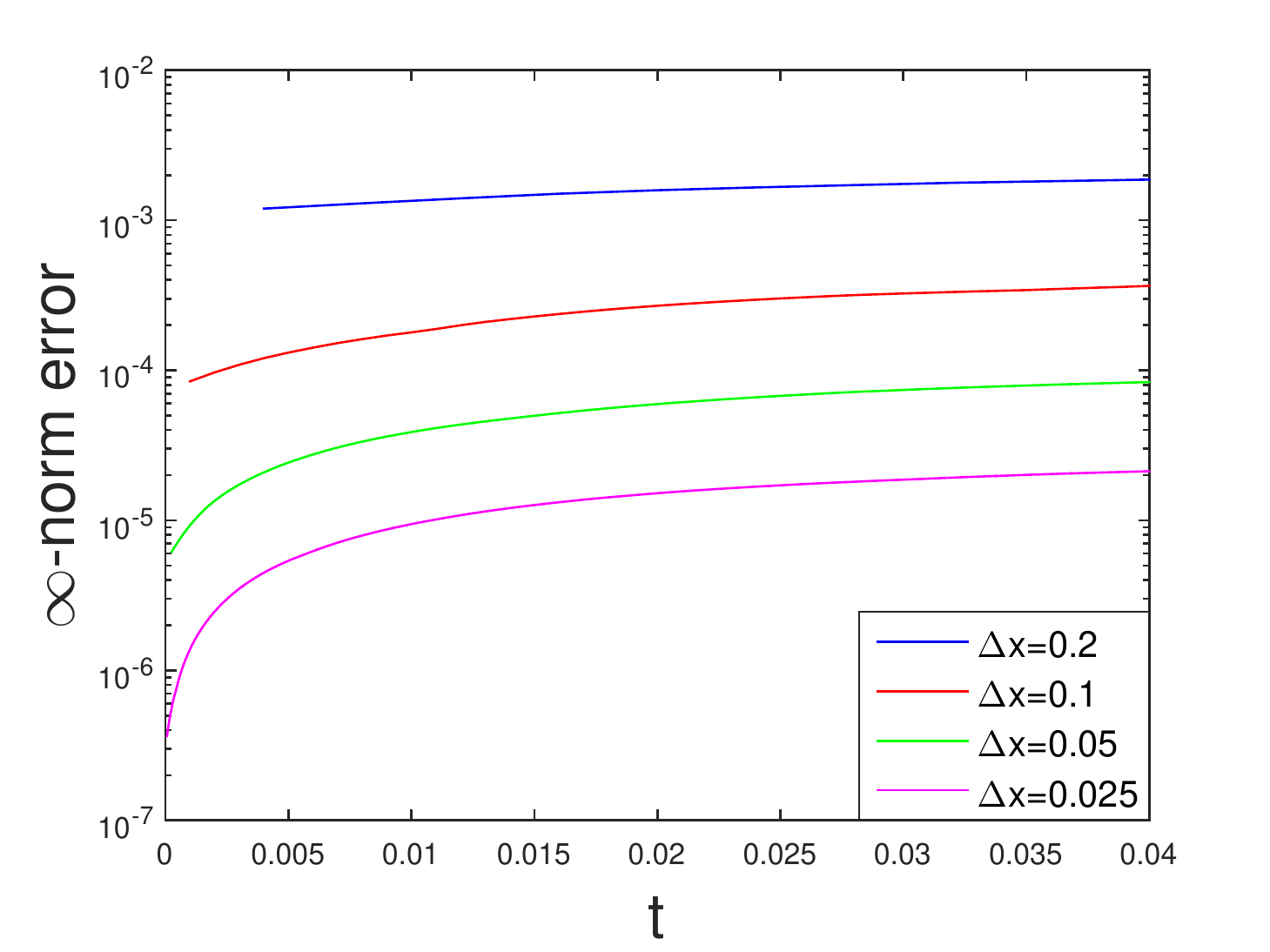}
    \caption{The absolute error of the numerical solution of the advection-diffusion model on the oscillating ellipsoid over time for various mesh spacings.}
    \label{errorEllipsoid}
\end{figure}

\begin{table}
\centering
  \begin{tabular}{|l|cc|cc|cc|cc|}
    \hline
    $\Delta x$ & $t =$ 0.01 & e.o.c. & $t =$ 0.02 & e.o.c. & $t =$ 0.03 & e.o.c. & $t =$ 0.04 & e.o.c.  \\ \hline
    0.2 & 1.30$\times 10^{-3}$ & - & 1.58$\times 10^{-3}$ & - & 1.72$\times 10^{-3}$ & - & 1.87$\times 10^{-3}$ & - \\
    0.1 & 1.79$\times 10^{-4}$ & 2.86 & 2.69$\times 10^{-4}$ & 2.56 & 3.25$\times 10^{-4}$ & 2.40 & 3.65$\times 10^{-4}$ & 2.35\\
    0.05 & 3.87$\times 10^{-5}$ & 2.21 & 5.95$\times 10^{-5}$ & 2.18 & 7.39$\times 10^{-5}$ & 2.14 & 8.34$\times 10^{-5}$ & 2.13\\
    0.025 & 9.43$\times 10^{-6}$ & 2.04 & 1.52$\times 10^{-5}$ & 1.97 & 1.87$\times 10^{-5}$ & 1.99 & 2.12$\times 10^{-5}$ & 1.97 \\
    \hline
  \end{tabular}
  \caption{Absolute error as measured in the infinity norm and the estimated order of convergence (e.o.c.) at times $t=0.01,0.02,0.03$ and $0.04$.}
  \label{errorEllipsoidTable}
\end{table}

Figure~\ref{errorEllipsoid} gives a plot of the error as a function of the elapsed time for various mesh sizes. We observe that the error accumulates in time, and that the error decreases with mesh refinement. Table~\ref{errorEllipsoidTable} gives a convergence study for times $t=0.01,0.02,0.03$ and $0.04$. In all four cases, second-order convergence is observed.

\subsubsection{Strongly coupled flow on an evolving torus}

Next, we evolve an initial torus according to a velocity
$$\mathbf{v} = (0.1\kappa+5u)\mathbf{n}$$
where $\kappa$ is the mean curvature, $u$ is the solution of the PDE (\ref{AdvDifPDE2}) and $\mathbf{n}$ is the unit normal vector. Because the surface velocity depends on the solution of the PDE, we call this system \emph{strongly coupled}. We select
$$u(\mathbf{x},0)=1+20x_1x_2x_3$$
for the initial condition.

\begin{figure}
    \centering
    \includegraphics[width=0.32\textwidth]{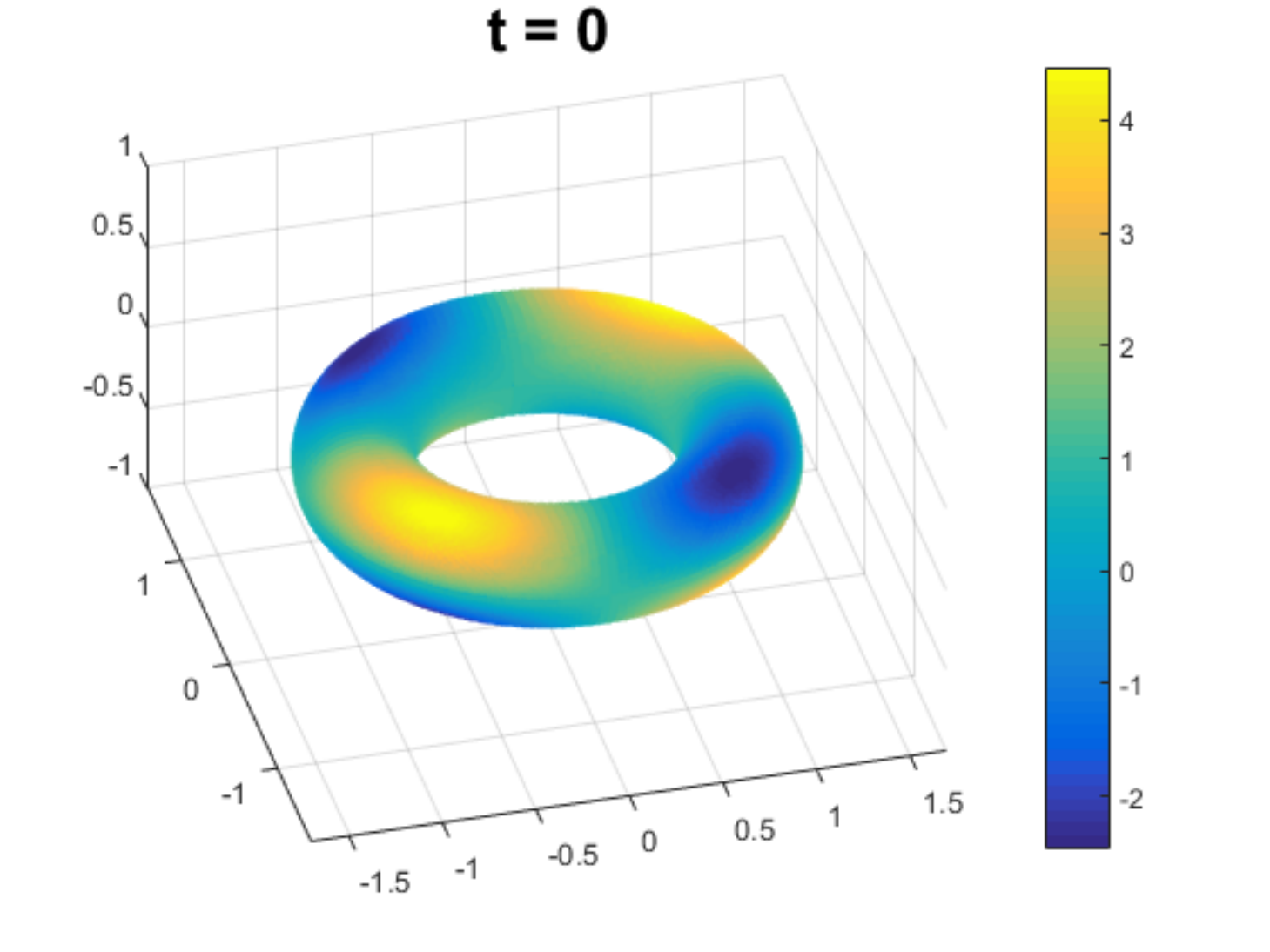}
    \includegraphics[width=0.32\textwidth]{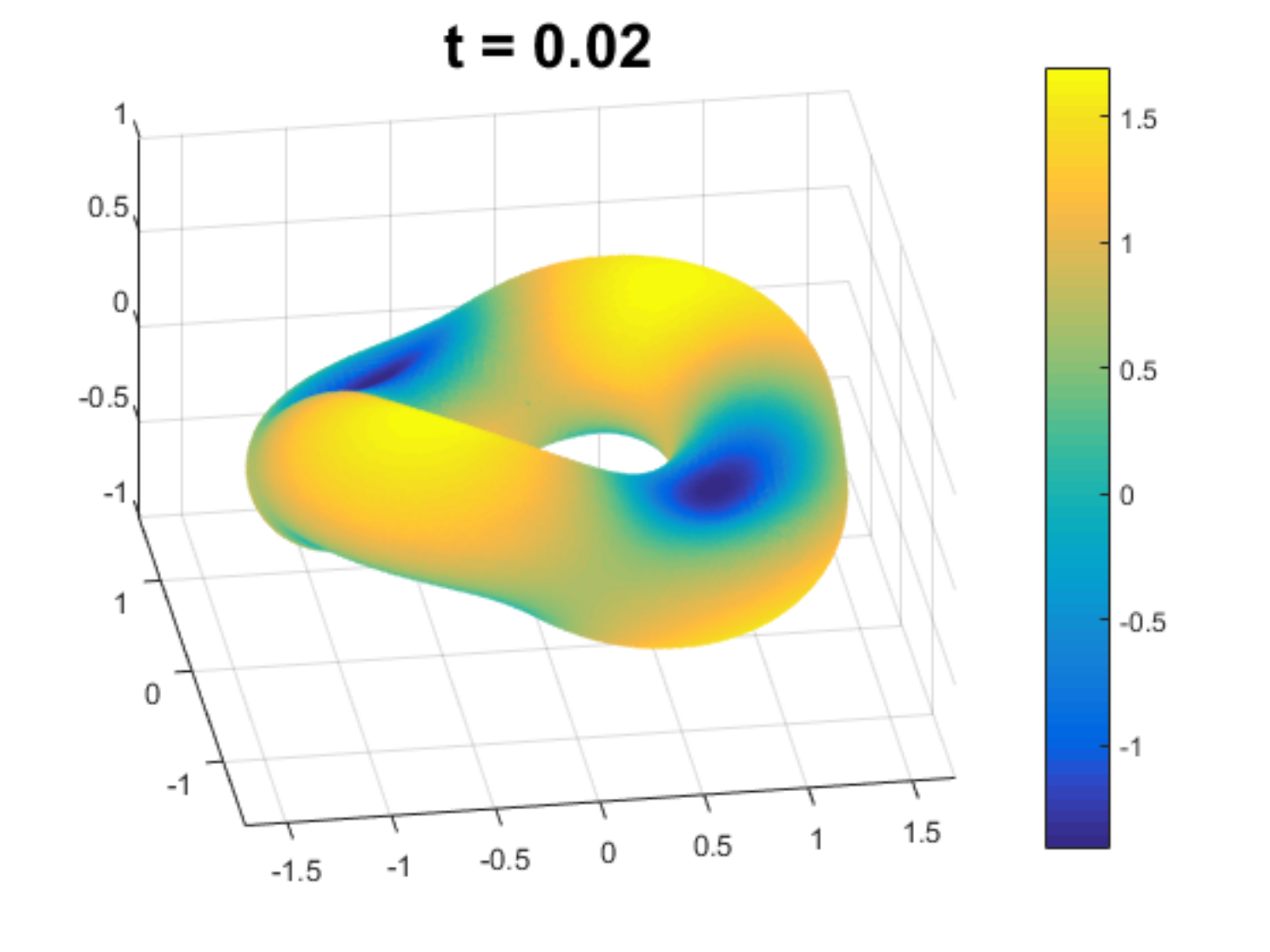}
    \includegraphics[width=0.32\textwidth]{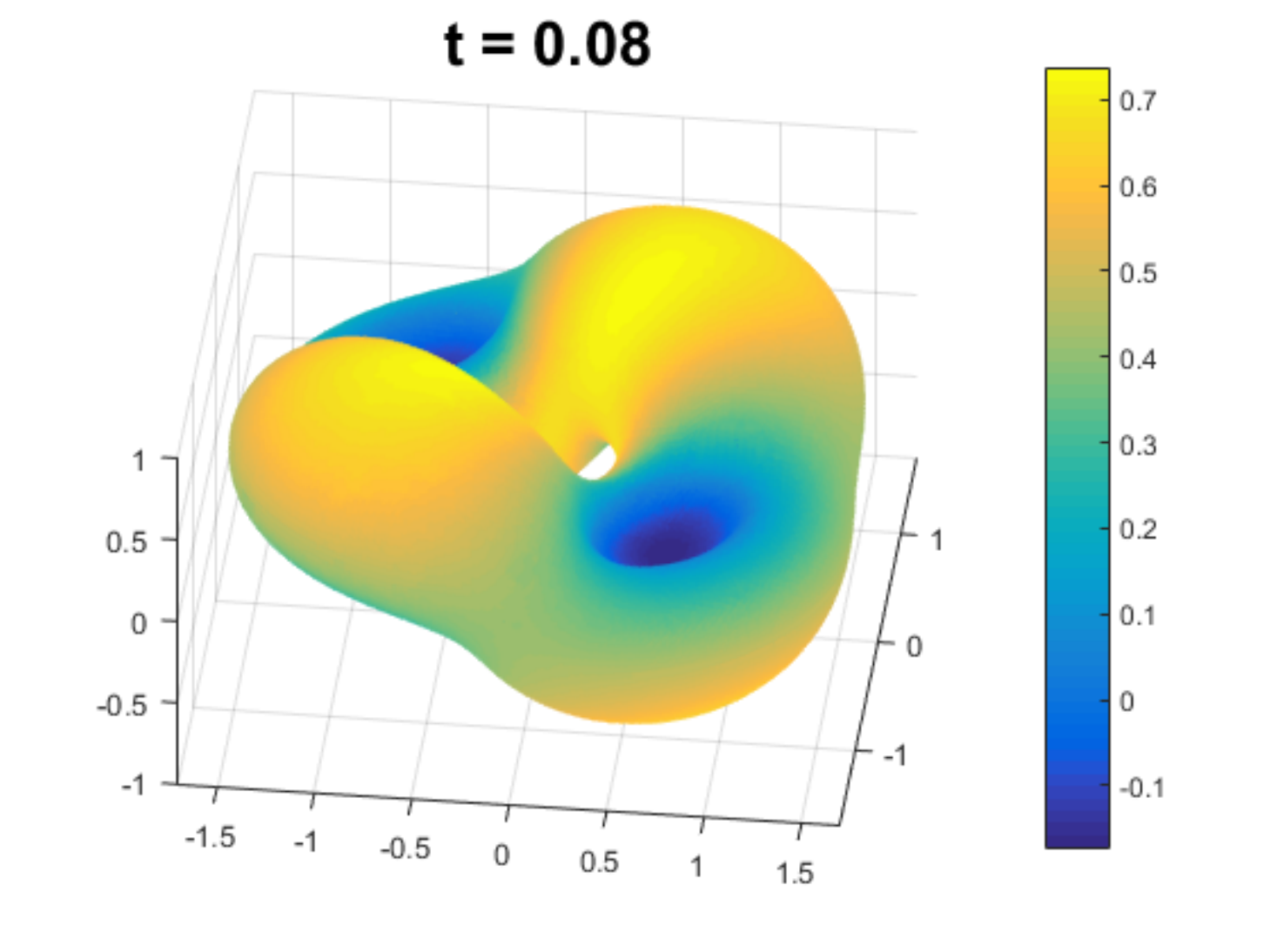}
    \caption{A visualization of the numerical solution of the strongly coupled flow at selected times $t$.  The initial surface is a torus.}
    \label{solTorus}
\end{figure}

Using a mesh spacing $\Delta x=0.05$ and a time step-size $\Delta t=0.2\Delta x^2$, we obtain the evolution displayed in Figure~\ref{solTorus}. In the initial flow, the torus experiences some rapid changes in the normal direction from the variation in the solution $u$. Over time, diffusion causes $u$ to become more uniform and the surface motion slows.

\subsubsection{A model for solid tumor growth}

Our final example considers solid tumor growth in the avascular growth phase, as presented in \cite{elliott2012ale}. The mathematical model for this example has the form
\begin{equation} \label{eq:tumormodel}
    \begin{array}{l}
        \displaystyle\frac{Du}{Dt} = \Delta_{\Gamma}u-u\nabla_{\Gamma}\cdot\mathbf{v}+f_1(u,w),\\
        \\
        \displaystyle\frac{Dw}{Dt} = \mathcal{D}\Delta_{\Gamma}w-w\nabla_{\Gamma}\cdot\mathbf{v}+f_2(u,w)\\
  \end{array}
\end{equation}
where $u$ and $w$ are scalar quantities on an evolving surface $\Gamma(t)$, $\frac{D}{Dt}$ is the material derivative, $\mathbf{v}$ is the velocity of $\Gamma(t)$, $\mathcal{D}$ is a positive constant and $f_1,f_2$ are functions that couple the solutions of the PDEs. For this example, we set
$$\begin{array}{l}
  \displaystyle f_1(u,w) = 100(0.1-u+u^2w),\\
  \\
  \displaystyle f_2(u,w) = 100(0.9-u^2w)
\end{array}$$
and
$$\mathbf{v}=(0.01\kappa+0.4u)\mathbf{n}$$
where $\mathbf{n}$ is the unit normal vector and $\kappa$ is the mean curvature. Finally, we set the diffusivity constant $\mathcal{D}=10$. Two identical initial conditions are chosen:
$$u(\mathbf{x},0) = 1+2x_1x_2x_3 = w(\mathbf{x},0).$$

\begin{figure}
    \centering
    \includegraphics[width=0.49\textwidth]{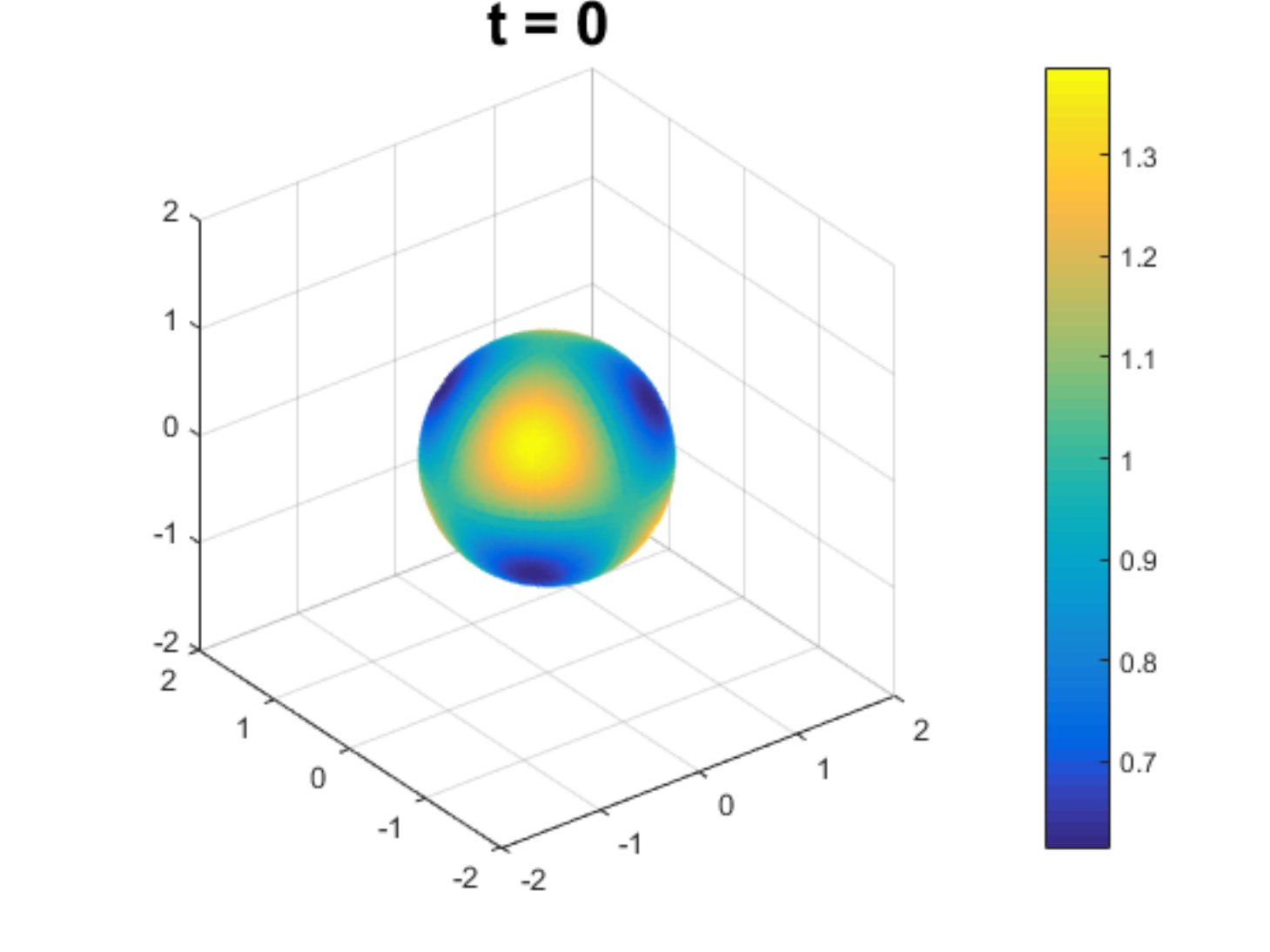}
    \includegraphics[width=0.49\textwidth]{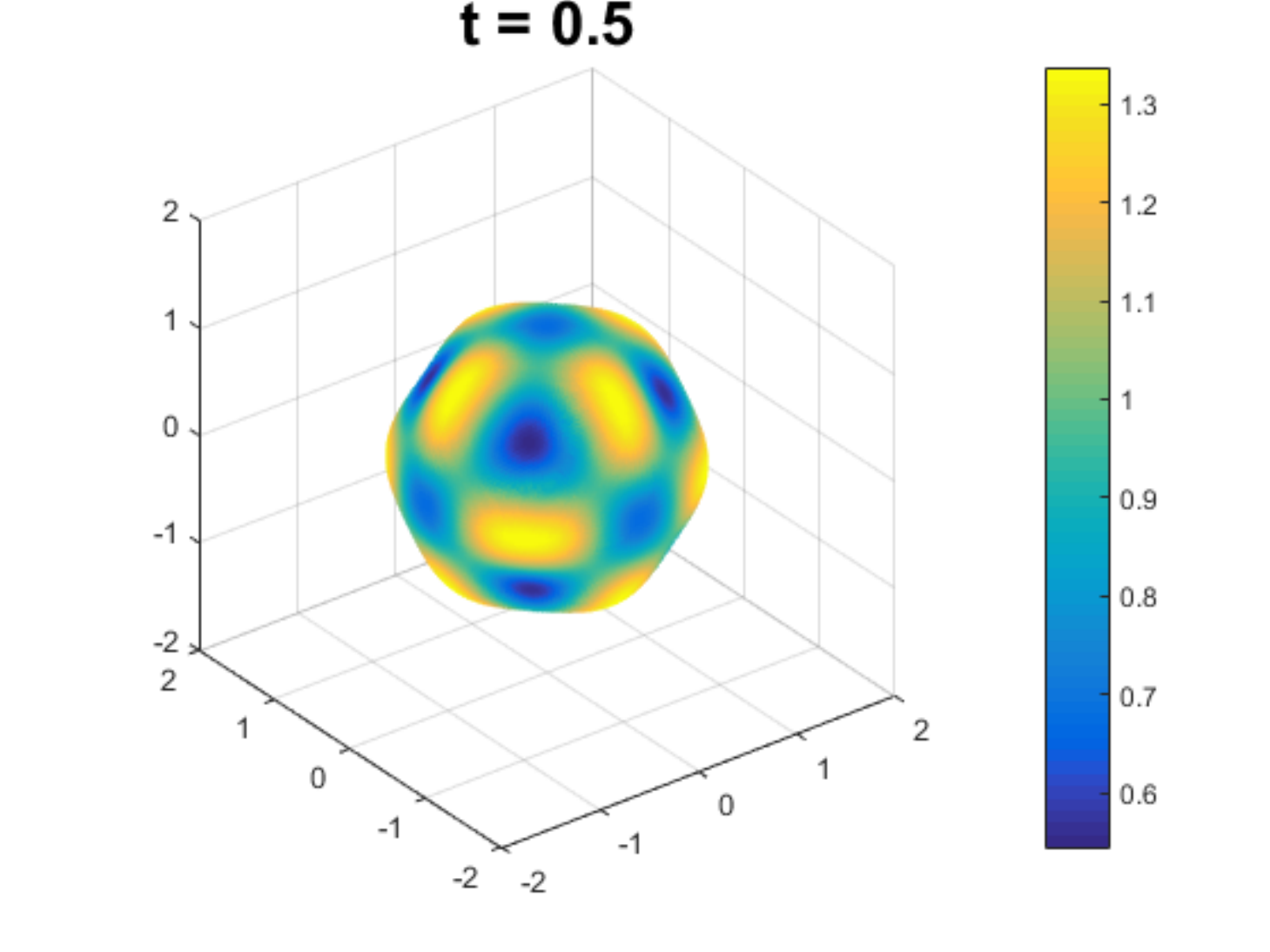}
    \includegraphics[width=0.49\textwidth]{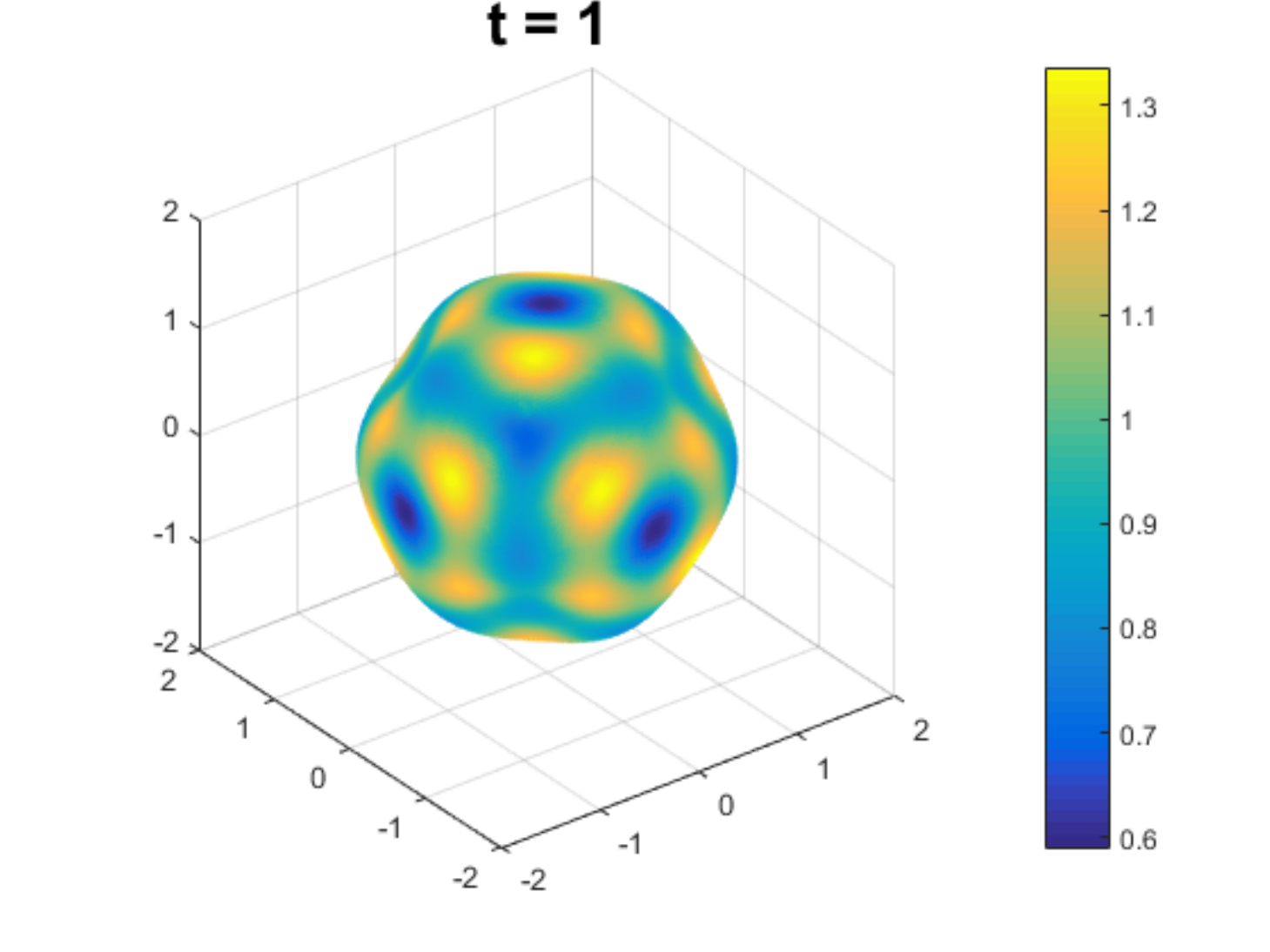}
    \includegraphics[width=0.49\textwidth]{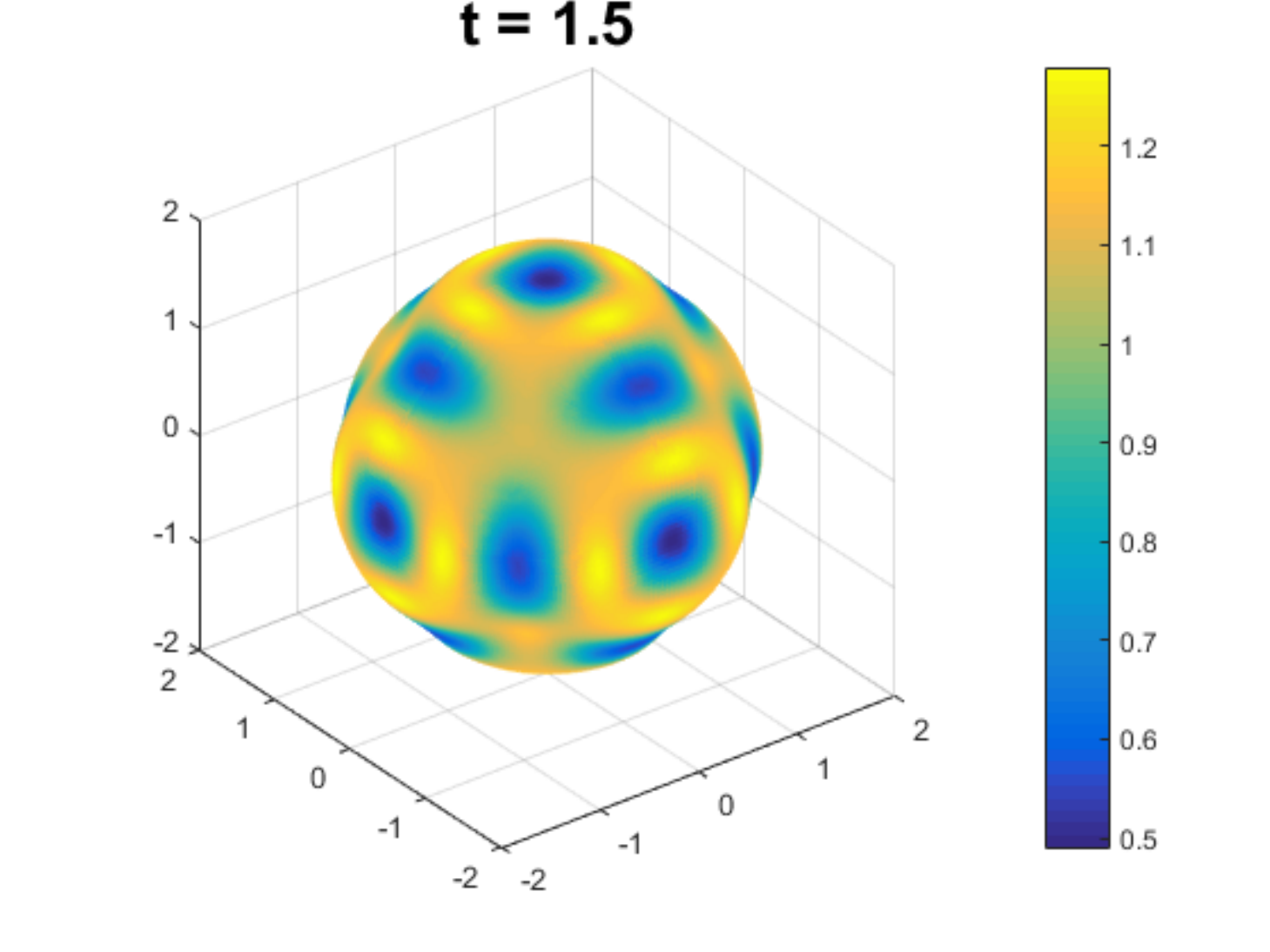}
    \caption{A visualization of the numerical solution $u$ of the tumor growth model~(\ref{eq:tumormodel}) at selected times $t$.  The initial surface is a sphere.}
    \label{solSphereU}
\end{figure}

\begin{figure}
    \centering
    \includegraphics[width=0.49\textwidth]{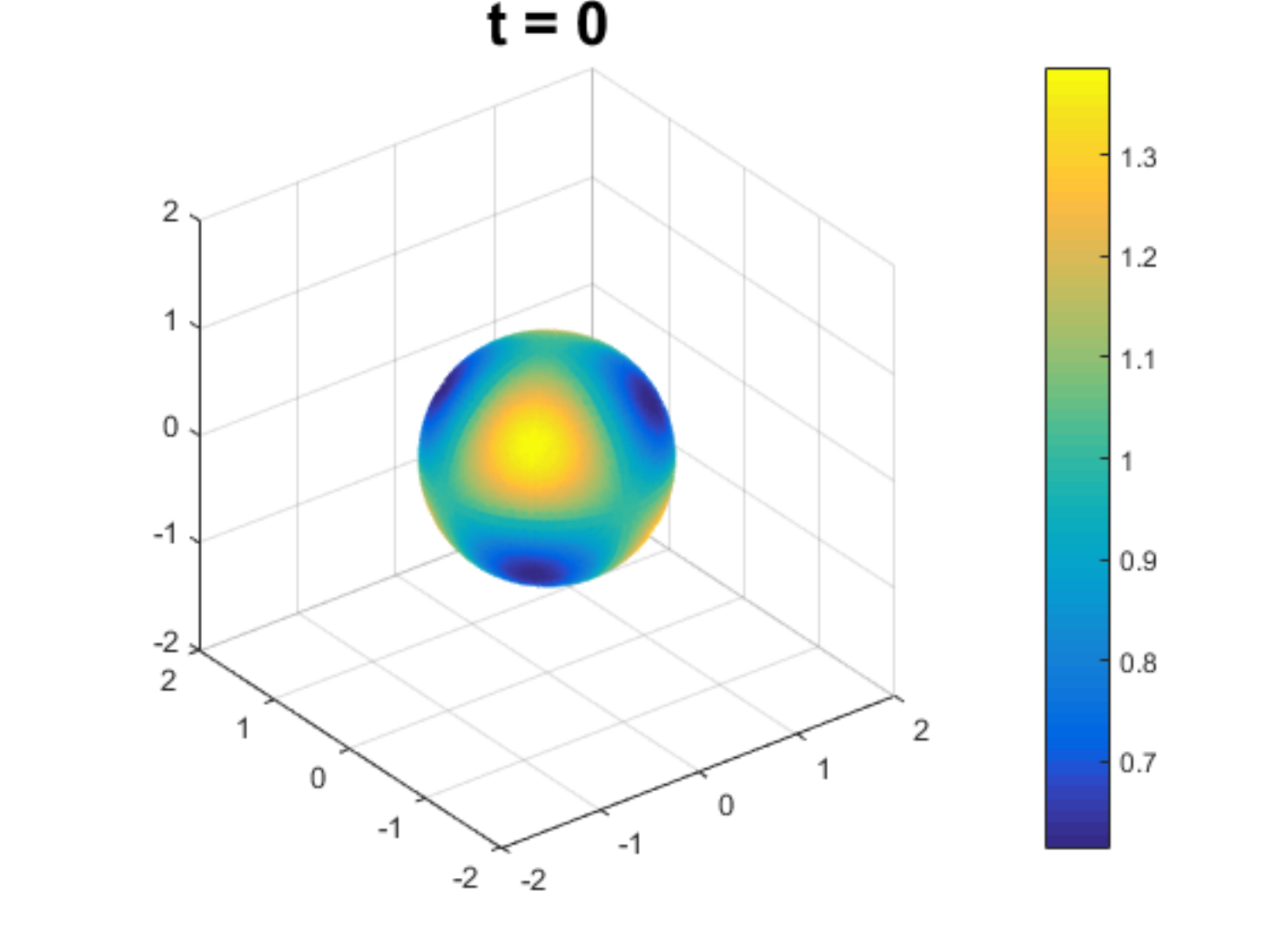}
    \includegraphics[width=0.49\textwidth]{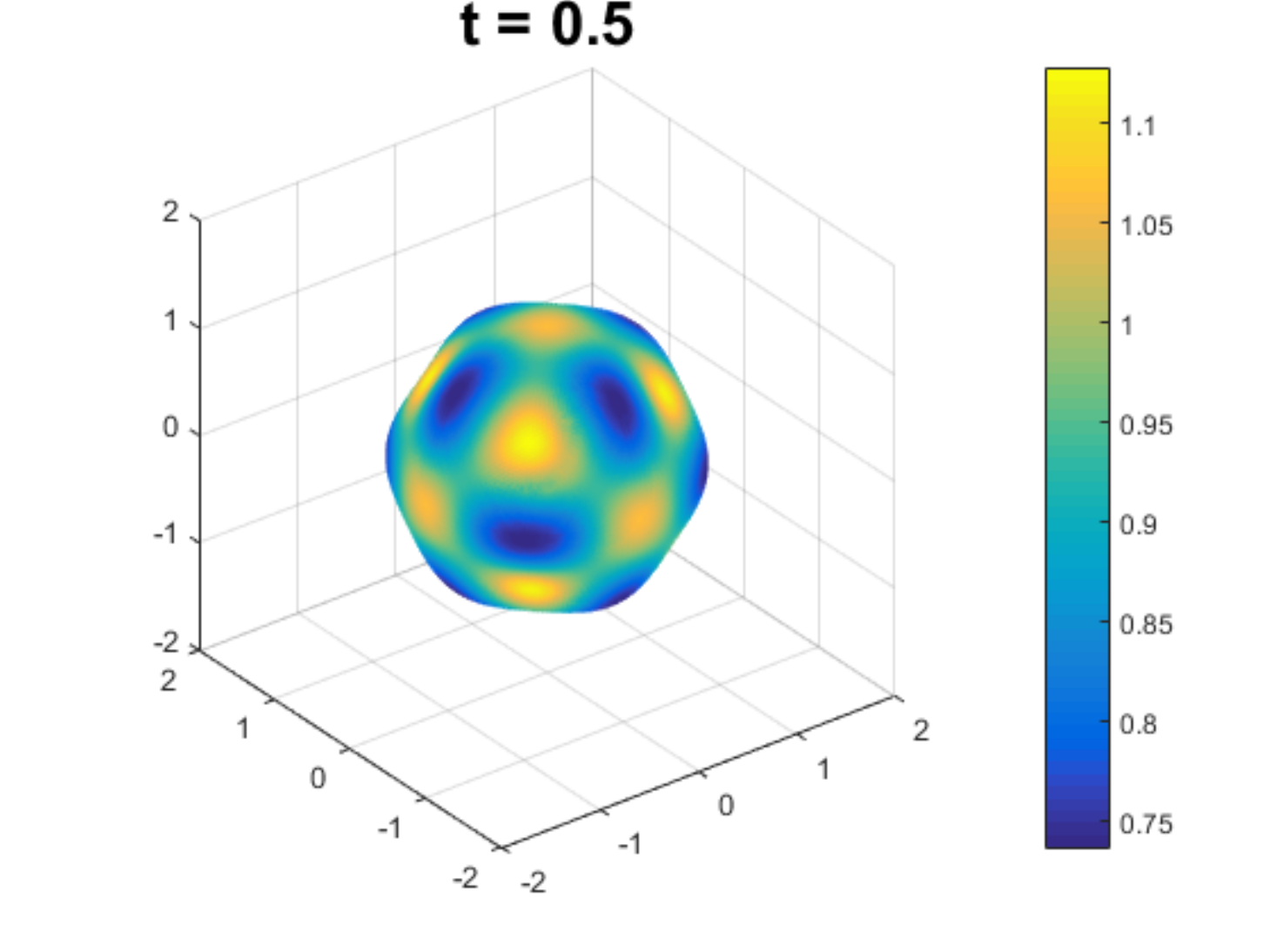}
    \includegraphics[width=0.49\textwidth]{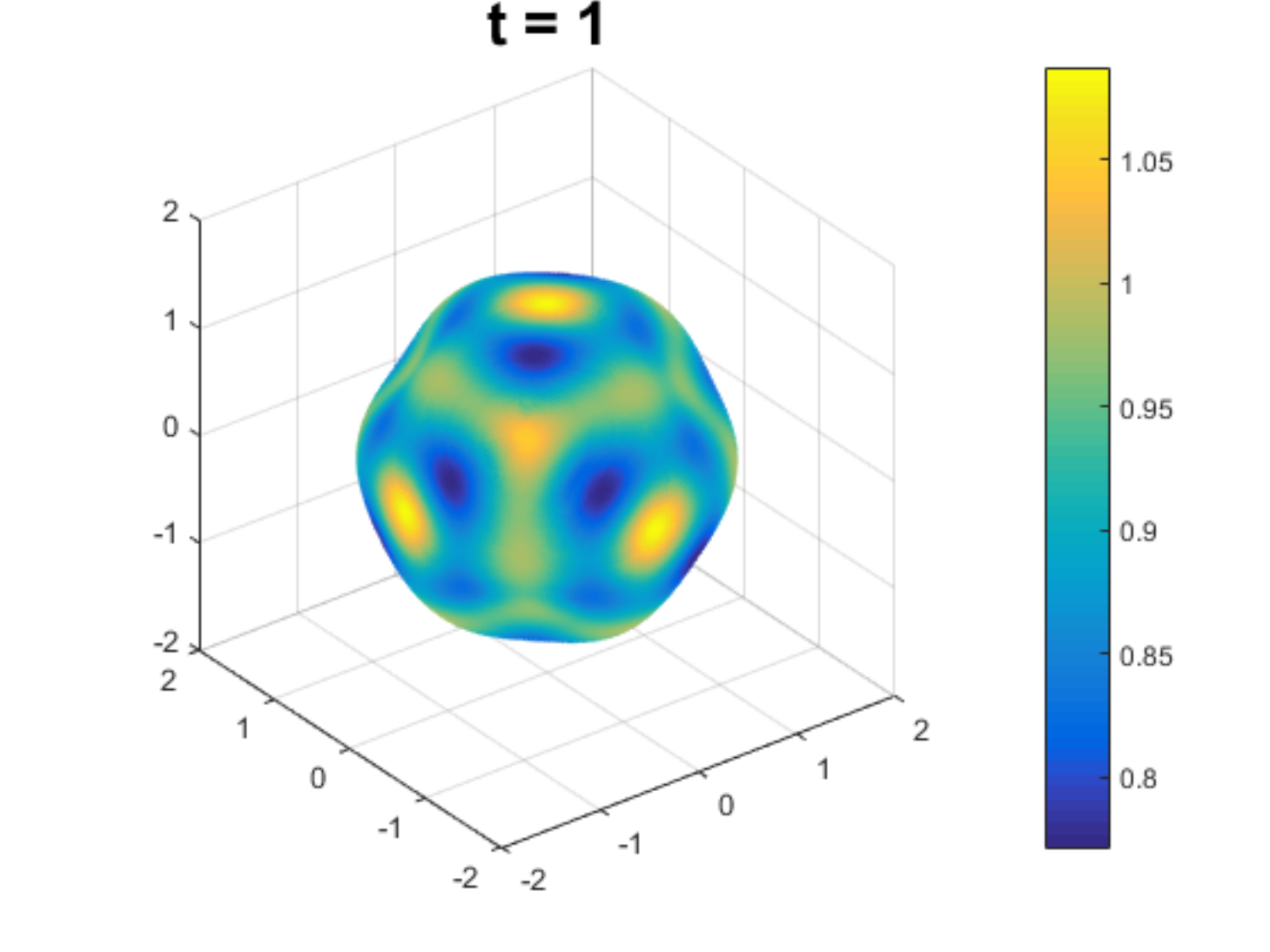}
    \includegraphics[width=0.49\textwidth]{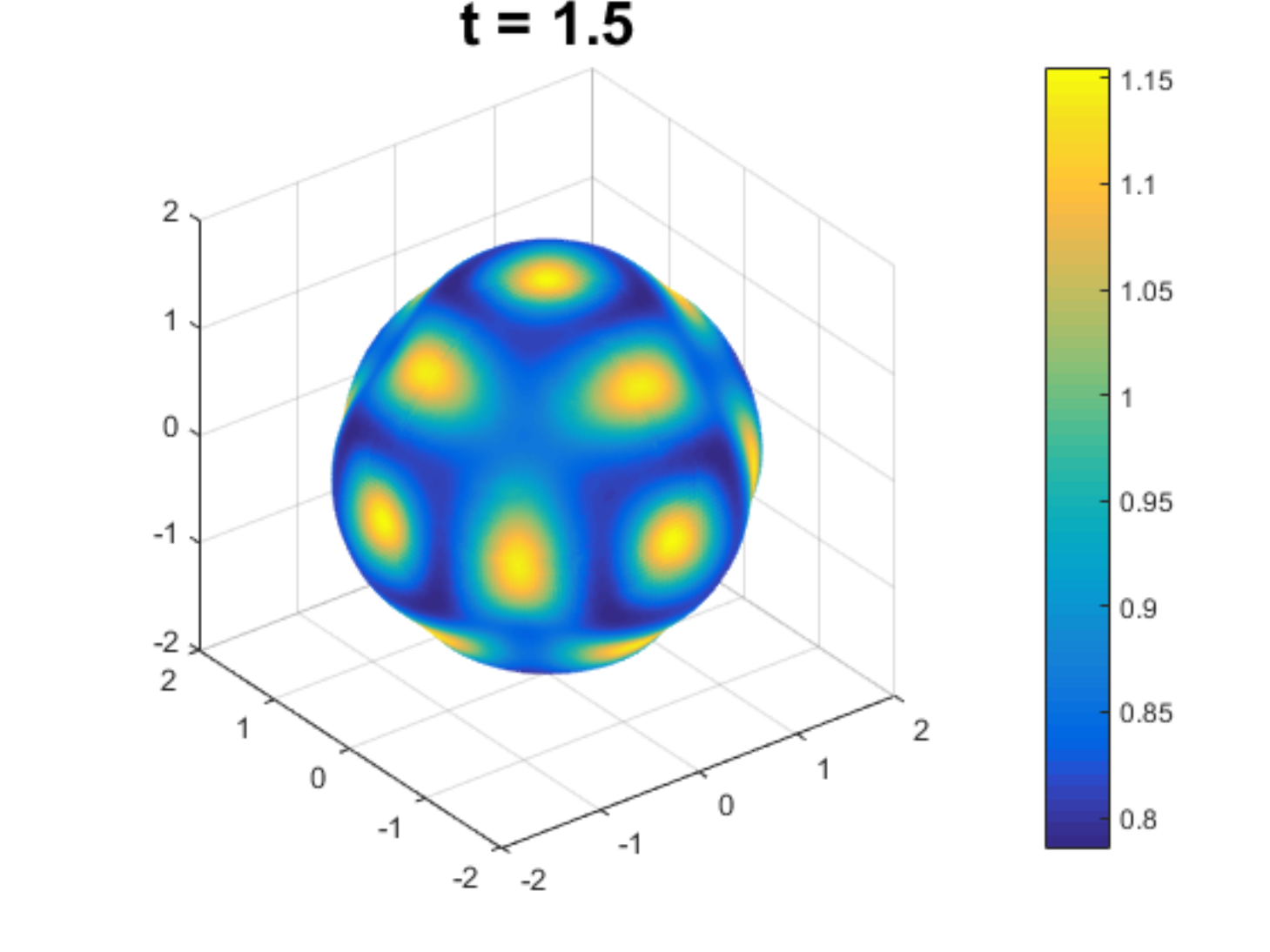}
    \caption{A visualization of the numerical solution $w$ of the tumor growth model~(\ref{eq:tumormodel}) at selected times $t$.  The initial surface is a sphere.}
    \label{solSphereW}
\end{figure}

Using a mesh spacing of $\Delta x=0.05$ and a time step-size $\Delta t=0.02\Delta x^2$ yields the results displayed in Figures~\ref{solSphereU} and \ref{solSphereW}. We observe tumor growth with preferred directions according to the solution $u$.

A smaller time step-size of $\Delta t=0.02\Delta x^2$ is used here to ensure stability for a relatively large diffusivity parameter $\mathcal{D}=10$. While this was sufficient for our purposes, stable, implicit methods are also available. See \cite{macdonald2009implicit} for details.

\section{Summary}\label{SummarySection}
In this paper, we combine the closest point method (CPM) and a (small) modification of the grid based particle method (GBPM)  in order to solve PDEs on moving surfaces. A modification of the GBPM was introduced to maintain all the grid points within the computational tube as active.
A coupled method that alternates between the CPM and the modified GBPM was proposed and applied to a number of examples of PDEs on surfaces. In numerical experiments, second-order convergence was observed as the grid was refined. The coupled method successfully approximates
the solution to curvature-dependent flows and systems with coupling between the PDEs and the surface motion.

Note that our proposed method should not be applied to problems involving
surfaces that are nonsmooth (at any time) since the original CPM is not designed for
nonsmooth surfaces.
Indeed, even on smooth surfaces, the computational domain must be sufficiently refined to avoid
the development of a multivalued closest point operator \cite{macdonald2008level}.
In recent work \cite{cheung2015localized}, we have found that a radial basis function discretization of the CPM can be
effective for solving PDEs on static surfaces with folds.
We are investigating combinations of such methods with the GBPM to solve PDEs on
nonsmooth moving surfaces, including surfaces undergoing topological change.

Other interesting directions for future work include the development of techniques for computing surface integrals, such as those appearing in volume conservation. Approximation methods developed for the cell based particle method \cite{hon2014cell} appear promising here. Note, however, that these methods do not directly apply to our situation, since they are based on $L^\infty$-distances rather than the $L^2$-distances we consider.

Other generalizations are also of strong interest. For example, it is natural to consider coupling the CPM and GBPM to solve PDEs on open surfaces and on surfaces moving by higher-order motion laws (cf., \cite{leung2009gridopen,leung2011grid}). Application of the coupled method to other PDE models is also of considerable interest. In particular, methods for the Cahn-Hilliard flux (\ref{CahnHilliardFlux}) are a natural target for future work.

\section*{Acknowledgements}
The work of the authors was partially supported by an NSERC Canada grant (RGPIN 227823).

\bibliographystyle{elsarticle-num-names}
\bibliography{draftBibliography}

\begin{thebibliography}{29}
\expandafter\ifx\csname natexlab\endcsname\relax\def\natexlab#1{#1}\fi
\providecommand{\url}[1]{\texttt{#1}}
\providecommand{\href}[2]{#2}
\providecommand{\path}[1]{#1}
\providecommand{\DOIprefix}{doi:}
\providecommand{\ArXivprefix}{arXiv:}
\providecommand{\URLprefix}{URL: }
\providecommand{\Pubmedprefix}{pmid:}
\providecommand{\doi}[1]{\href{http://dx.doi.org/#1}{\path{#1}}}
\providecommand{\Pubmed}[1]{\href{pmid:#1}{\path{#1}}}
\providecommand{\bibinfo}[2]{#2}
\ifx\xfnm\relax \def\xfnm[#1]{\unskip,\space#1}\fi
\bibitem[{Elliott and Stinner(2010)}]{elliott2010modeling}
\bibinfo{author}{C.~M. Elliott}, \bibinfo{author}{B.~Stinner},
\newblock \bibinfo{title}{Modeling and computation of two phase geometric
  biomembranes using surface finite elements},
\newblock \bibinfo{journal}{Journal of Computational Physics}
  \bibinfo{volume}{229} (\bibinfo{year}{2010}) \bibinfo{pages}{6585--6612}.
\bibitem[{Elliott et~al.(2012)Elliott, Stinner, and
  Venkataraman}]{elliott2012modelling}
\bibinfo{author}{C.~M. Elliott}, \bibinfo{author}{B.~Stinner},
  \bibinfo{author}{C.~Venkataraman},
\newblock \bibinfo{title}{Modelling cell motility and chemotaxis with evolving
  surface finite elements},
\newblock \bibinfo{journal}{Journal of The Royal Society Interface}
  (\bibinfo{year}{2012}) \bibinfo{pages}{rsif20120276}.
\bibitem[{Barreira et~al.(2011)Barreira, Elliott, and
  Madzvamuse}]{barreira2011surface}
\bibinfo{author}{R.~Barreira}, \bibinfo{author}{C.~M. Elliott},
  \bibinfo{author}{A.~Madzvamuse},
\newblock \bibinfo{title}{The surface finite element method for pattern
  formation on evolving biological surfaces},
\newblock \bibinfo{journal}{Journal of Mathematical Biology}
  \bibinfo{volume}{63} (\bibinfo{year}{2011}) \bibinfo{pages}{1095--1119}.
\bibitem[{Venkataraman et~al.(2011)Venkataraman, Sekimura, Gaffney, Maini, and
  Madzvamuse}]{venkataraman2011modeling}
\bibinfo{author}{C.~Venkataraman}, \bibinfo{author}{T.~Sekimura},
  \bibinfo{author}{E.~A. Gaffney}, \bibinfo{author}{P.~K. Maini},
  \bibinfo{author}{A.~Madzvamuse},
\newblock \bibinfo{title}{Modeling parr-mark pattern formation during the early
  development of amago trout},
\newblock \bibinfo{journal}{Physical Review E} \bibinfo{volume}{84}
  (\bibinfo{year}{2011}) \bibinfo{pages}{041923}.
\bibitem[{Eilks and Elliott(2008)}]{eilks2008numerical}
\bibinfo{author}{C.~Eilks}, \bibinfo{author}{C.~M. Elliott},
\newblock \bibinfo{title}{Numerical simulation of dealloying by surface
  dissolution via the evolving surface finite element method},
\newblock \bibinfo{journal}{Journal of Computational Physics}
  \bibinfo{volume}{227} (\bibinfo{year}{2008}) \bibinfo{pages}{9727--9741}.
\bibitem[{Adalsteinsson and Sethian(2003)}]{adalsteinsson2003transport}
\bibinfo{author}{D.~Adalsteinsson}, \bibinfo{author}{J.~A. Sethian},
\newblock \bibinfo{title}{Transport and diffusion of material quantities on
  propagating interfaces via level set methods},
\newblock \bibinfo{journal}{Journal of Computational Physics}
  \bibinfo{volume}{185} (\bibinfo{year}{2003}) \bibinfo{pages}{271--288}.
\bibitem[{James and Lowengrub(2004)}]{james2004surfactant}
\bibinfo{author}{A.~J. James}, \bibinfo{author}{J.~Lowengrub},
\newblock \bibinfo{title}{A surfactant-conserving volume-of-fluid method for
  interfacial flows with insoluble surfactant},
\newblock \bibinfo{journal}{Journal of Computational Physics}
  \bibinfo{volume}{201} (\bibinfo{year}{2004}) \bibinfo{pages}{685--722}.
\bibitem[{Auer and Westermann(2013)}]{auer2013semi}
\bibinfo{author}{S.~Auer}, \bibinfo{author}{R.~Westermann},
\newblock \bibinfo{title}{A semi-lagrangian closest point method for deforming
  surfaces},
\newblock in: \bibinfo{booktitle}{Computer Graphics Forum},
  volume~\bibinfo{volume}{32}, \bibinfo{organization}{Wiley Online Library},
  \bibinfo{year}{2013}, pp. \bibinfo{pages}{207--214}.
\bibitem[{Dziuk and Elliott(2007)}]{dziuk2007finite}
\bibinfo{author}{G.~Dziuk}, \bibinfo{author}{C.~M. Elliott},
\newblock \bibinfo{title}{Finite elements on evolving surfaces},
\newblock \bibinfo{journal}{IMA journal of numerical analysis}
  \bibinfo{volume}{27} (\bibinfo{year}{2007}) \bibinfo{pages}{262--292}.
\bibitem[{Dziuk and Elliott(2013)}]{dziuk2013finite}
\bibinfo{author}{G.~Dziuk}, \bibinfo{author}{C.~M. Elliott},
\newblock \bibinfo{title}{Finite element methods for surface pdes},
\newblock \bibinfo{journal}{Acta Numerica} \bibinfo{volume}{22}
  (\bibinfo{year}{2013}) \bibinfo{pages}{289--396}.
\bibitem[{Nemadjieu(2012)}]{nemadjieu2012finite}
\bibinfo{author}{S.~F. Nemadjieu}, \bibinfo{title}{Finite volume methods for
  advection diffusion on moving interfaces and application on surfactant driven
  thin film flow}, Ph.D. thesis, Faculty of Mathematics and Natural Sciences,
  University of Bonn, \bibinfo{year}{2012}.
\bibitem[{Leung and Zhao(2009)}]{leung2009grid}
\bibinfo{author}{S.~Leung}, \bibinfo{author}{H.~Zhao},
\newblock \bibinfo{title}{A grid based particle method for moving interface
  problems},
\newblock \bibinfo{journal}{Journal of Computational Physics}
  \bibinfo{volume}{228} (\bibinfo{year}{2009}) \bibinfo{pages}{2993--3024}.
\bibitem[{Leung et~al.(2011)Leung, Lowengrub, and Zhao}]{leung2011grid}
\bibinfo{author}{S.~Leung}, \bibinfo{author}{J.~Lowengrub},
  \bibinfo{author}{H.~Zhao},
\newblock \bibinfo{title}{A grid based particle method for solving partial
  differential equations on evolving surfaces and modeling high order
  geometrical motion},
\newblock \bibinfo{journal}{Journal of Computational Physics}
  \bibinfo{volume}{230} (\bibinfo{year}{2011}) \bibinfo{pages}{2540--2561}.
\bibitem[{Leung and Zhao(2010)}]{leung2010gaussian}
\bibinfo{author}{S.~Leung}, \bibinfo{author}{H.~Zhao},
\newblock \bibinfo{title}{Gaussian beam summation for diffraction in
  inhomogeneous media based on the grid based particle method},
\newblock \bibinfo{journal}{Communications in Computational Physics}
  \bibinfo{volume}{8} (\bibinfo{year}{2010}) \bibinfo{pages}{758}.
\bibitem[{Sethian(1999)}]{sethian1999level}
\bibinfo{author}{J.~A. Sethian}, \bibinfo{title}{Level set methods and fast
  marching methods: evolving interfaces in computational geometry, fluid
  mechanics, computer vision, and materials science},
  volume~\bibinfo{volume}{3}, \bibinfo{publisher}{Cambridge University Press},
  \bibinfo{year}{1999}.
\bibitem[{Osher and Fedkiw(2006)}]{osher2006level}
\bibinfo{author}{S.~Osher}, \bibinfo{author}{R.~Fedkiw}, \bibinfo{title}{Level
  set methods and dynamic implicit surfaces}, volume \bibinfo{volume}{153},
  \bibinfo{publisher}{Springer Science \& Business Media},
  \bibinfo{year}{2006}.
\bibitem[{Xu and Zhao(2003)}]{xu2003eulerian}
\bibinfo{author}{J.-J. Xu}, \bibinfo{author}{H.-K. Zhao},
\newblock \bibinfo{title}{An eulerian formulation for solving partial
  differential equations along a moving interface},
\newblock \bibinfo{journal}{Journal of Scientific Computing}
  \bibinfo{volume}{19} (\bibinfo{year}{2003}) \bibinfo{pages}{573--594}.
\bibitem[{Dziuk and Elliott(2010)}]{dziuk2010eulerian}
\bibinfo{author}{G.~Dziuk}, \bibinfo{author}{C.~M. Elliott},
\newblock \bibinfo{title}{An eulerian approach to transport and diffusion on
  evolving implicit surfaces},
\newblock \bibinfo{journal}{Computing and Visualization in Science}
  \bibinfo{volume}{13} (\bibinfo{year}{2010}) \bibinfo{pages}{17--28}.
\bibitem[{Kim et~al.(2013)Kim, Tessendorf, and Thuerey}]{kim2013closest}
\bibinfo{author}{T.~Kim}, \bibinfo{author}{J.~Tessendorf},
  \bibinfo{author}{N.~Thuerey},
\newblock \bibinfo{title}{Closest point turbulence for liquid surfaces},
\newblock \bibinfo{journal}{ACM Transactions on Graphics (TOG)}
  \bibinfo{volume}{32} (\bibinfo{year}{2013}) \bibinfo{pages}{15}.
\bibitem[{Ruuth and Merriman(2008)}]{ruuth2008simple}
\bibinfo{author}{S.~J. Ruuth}, \bibinfo{author}{B.~Merriman},
\newblock \bibinfo{title}{A simple embedding method for solving partial
  differential equations on surfaces},
\newblock \bibinfo{journal}{Journal of Computational Physics}
  \bibinfo{volume}{227} (\bibinfo{year}{2008}) \bibinfo{pages}{1943--1961}.
\bibitem[{Macdonald and Ruuth(2009)}]{macdonald2009implicit}
\bibinfo{author}{C.~B. Macdonald}, \bibinfo{author}{S.~J. Ruuth},
\newblock \bibinfo{title}{The implicit closest point method for the numerical
  solution of partial differential equations on surfaces},
\newblock \bibinfo{journal}{SIAM Journal on Scientific Computing}
  \bibinfo{volume}{31} (\bibinfo{year}{2009}) \bibinfo{pages}{4330--4350}.
\bibitem[{Macdonald et~al.(2011)Macdonald, Brandman, and
  Ruuth}]{Macdonald20117944}
\bibinfo{author}{C.~B. Macdonald}, \bibinfo{author}{J.~Brandman},
  \bibinfo{author}{S.~J. Ruuth},
\newblock \bibinfo{title}{Solving eigenvalue problems on curved surfaces using
  the closest point method},
\newblock \bibinfo{journal}{Journal of Computational Physics}
  \bibinfo{volume}{230} (\bibinfo{year}{2011}) \bibinfo{pages}{7944 -- 7956}.
\bibitem[{Leung and Zhao(2009)}]{leung2009gridopen}
\bibinfo{author}{S.~Leung}, \bibinfo{author}{H.~Zhao},
\newblock \bibinfo{title}{A grid based particle method for evolution of open
  curves and surfaces},
\newblock \bibinfo{journal}{Journal of Computational Physics}
  \bibinfo{volume}{228} (\bibinfo{year}{2009}) \bibinfo{pages}{7706--7728}.
\bibitem[{Hon et~al.(2014)Hon, Leung, and Zhao}]{hon2014cell}
\bibinfo{author}{S.~Y. Hon}, \bibinfo{author}{S.~Leung},
  \bibinfo{author}{H.~Zhao},
\newblock \bibinfo{title}{A cell based particle method for modeling dynamic
  interfaces},
\newblock \bibinfo{journal}{Journal of Computational Physics}
  \bibinfo{volume}{272} (\bibinfo{year}{2014}) \bibinfo{pages}{279--306}.
\bibitem[{LeVeque(1996)}]{leveque1996high}
\bibinfo{author}{R.~J. LeVeque},
\newblock \bibinfo{title}{High-resolution conservative algorithms for advection
  in incompressible flow},
\newblock \bibinfo{journal}{SIAM Journal on Numerical Analysis}
  \bibinfo{volume}{33} (\bibinfo{year}{1996}) \bibinfo{pages}{627--665}.
\bibitem[{Elliott et~al.(2011)Elliott, Stinner, Styles, and
  Welford}]{elliott2011numerical}
\bibinfo{author}{C.~M. Elliott}, \bibinfo{author}{B.~Stinner},
  \bibinfo{author}{V.~Styles}, \bibinfo{author}{R.~Welford},
\newblock \bibinfo{title}{Numerical computation of advection and diffusion on
  evolving diffuse interfaces},
\newblock \bibinfo{journal}{IMA Journal of Numerical Analysis}
  \bibinfo{volume}{31} (\bibinfo{year}{2011}) \bibinfo{pages}{786--812}.
\bibitem[{Elliott and Styles(2012)}]{elliott2012ale}
\bibinfo{author}{C.~M. Elliott}, \bibinfo{author}{V.~Styles},
\newblock \bibinfo{title}{An ale esfem for solving pdes on evolving surfaces},
\newblock \bibinfo{journal}{Milan Journal of Mathematics}
  (\bibinfo{year}{2012}) \bibinfo{pages}{1--33}.
\bibitem[{Macdonald and Ruuth(2008)}]{macdonald2008level}
\bibinfo{author}{C.~B. Macdonald}, \bibinfo{author}{S.~J. Ruuth},
\newblock \bibinfo{title}{Level set equations on surfaces via the closest point
  method},
\newblock \bibinfo{journal}{Journal of Scientific Computing}
  \bibinfo{volume}{35} (\bibinfo{year}{2008}) \bibinfo{pages}{219--240}.
\bibitem[{Cheung et~al.(2015)Cheung, Ling, and Ruuth}]{cheung2015localized}
\bibinfo{author}{K.~C. Cheung}, \bibinfo{author}{L.~Ling},
  \bibinfo{author}{S.~Ruuth},
\newblock \bibinfo{title}{A localized meshless method for diffusion on folded
  surfaces},
\newblock \bibinfo{journal}{Journal of Computational Physics}
  \bibinfo{volume}{297} (\bibinfo{year}{2015}) \bibinfo{pages}{194--206}.

\end{thebibliography}

\end{document}